\documentclass{article}
\usepackage{amssymb,amsmath,amsthm,amsfonts,latexsym}
\usepackage{newtxtext}

\usepackage{bbm}
\usepackage[toc,page]{appendix}
\usepackage{tikz}
\usepackage{stackrel}
\usepackage{fancyvrb}
\usepackage{graphicx, color, subfigure}
\usepackage{siunitx}
\usepgflibrary{arrows}
\newtheorem{theorem}{Theorem}[section]
\newtheorem{lemma}[theorem]{Lemma}

\newtheorem{remark}{Remark}[section]
\usepackage{caption}
\usepackage{afterpage}
\usepackage{fullpage}
\usepackage{mathtools}
\usepackage[mathscr]{eucal}
\usepackage{pdfpages}
\usepackage[T1]{fontenc}
\usepackage[english]{babel}
\usepackage{xspace}
\usepackage{amscd,bm}
\usepackage{xcolor,soul}
\usepackage{verbatim}
\usepackage{csvsimple}
\numberwithin{equation}{section}

\newtheorem{ex}{Example}[section]

\usepackage{hyperref}
\hypersetup{
    colorlinks=true,
    linkcolor=black,
    filecolor=black,
    urlcolor=black,
    citecolor=black,
}
\title{The hybrid dimensional representation of permeability tensor: a reinterpretation of the discrete fracture model and its extension on nonconforming meshes\footnote{\baselineskip 1.2pc Supported by the NSF grant DMS-1818467}}
\author{Ziyao Xu\footnote{Department of Mathematical Sciences,
Michigan Technological University, Houghton, MI 49931. E-mail: ziyaox@mtu.edu},\quad Yang Yang\footnote{Department of Mathematical Sciences,
Michigan Technological University, Houghton, MI 49931. E-mail: yyang7@mtu.edu}}
\date{}
\begin{document}
\baselineskip 1.2pc
\maketitle
\begin{abstract}
\baselineskip 1.2pc
The discrete fracture model (DFM) has been widely used in the simulation of fluid flow in fractured porous media. Traditional DFM uses the so-called hybrid-dimensional approach to treat fractures explicitly as low-dimensional entries (e.g. line entries in 2D media and face entries in 3D media) on the interfaces of matrix cells and then couple the matrix and fracture flow systems together based on the principle of superposition with the fracture thickness used as the dimensional homogeneity factor. Because of this methodology, DFM is considered to be limited on conforming meshes and thus may raise difficulties in generating high quality unstructured meshes due to the complexity of fracture's geometrical morphology. In this paper, we clarify that the DFM actually can be extended to non-conforming meshes without any essential changes. To show it clearly, we provide another perspective for DFM based on hybrid-dimensional representation of permeability tensor to describe fractures as one-dimensional line Dirac delta functions contained in permeability tensor. A finite element DFM scheme for single-phase flow on non-conforming meshes is then derived by applying Galerkin finite element method to it. Analytical analysis and numerical experiments show that our DFM automatically degenerates to the classical finite element DFM when the mesh is conforming with fractures. Moreover, the accuracy and efficiency of the model on non-conforming meshes are demonstrated by testing several benchmark problems. This model is also applicable to curved fracture with variable thickness.
\end{abstract}
\section{Introduction}
{As an important model problem arising from a variety of applications including enhanced oil recovery in naturally fractured reservoirs, contaminant transport in fractured rocks, and radioactive waste repository in subsurface, the study on fluid flow in fractured porous media is of high interest and has engaged a large number of researchers in the past half-century. Moreover, due to the recent prevalence of hydraulic fracturing techniques developed for unconventional reservoirs, efficient and accurate simulators for flow and transport in fractured media are increasingly desired.}

{The fractured porous media is composed of the highly conductive narrow fractures and the surrounding low-permeability rock matrix. Essentially, the fractured porous media is a special porous media with extreme heterogeneity contributed by the fractured regions.}

{Though with tiny thickness, the fractures have non-negligible effect on the flow in fractured media because of its high conductivity. Therefore, to efficiently and accurately build in the effect of fractures on the flow in fractured media is crucial to the simulation but remains to be challenging due to the extreme contrast of scale and permeability between porous matrix and fracture, as well as the complexity of the fracture network.}

{There are various remarkable models developed for the simulation of flow in fractured porous media in the past decades, among which the most widely used methods are dual-porosity model, single-porosity model, discrete fracture model (DFM), and embedded discrete fracture model (EDFM), etc. Moreover, the XFEM-class methods, mortar-type approaches, Lagrange multiplier methods and interface models are also developed in recent years. These approaches can be roughly divided into two categories: the continuum model and discrete fracture-matrix model.}

{The dual-porosity model, see e.g. \cite{DualPoro1,DualPoro2,DualPoro3,DualPoro4,DualPoro5,benchmark3}, is a widely practiced continuum model in fractured reservoir simulations. In dual-porosity model, the matrix and fracture network are both treated as continuum systems governed by Darcy's law with different permeability. The mass transfer between two systems is given by the matrix–fracture mass transfer function, which is determined by the pressure difference between matrix and fractures systems and the characteristic properties of rock and fluid. By simplifying the fracture network to continuum media, the dual-porosity model gained great advantage in efficiency and becomes a major approach used in the field-scale simulation of naturally fractured reservoirs. However, as a typical continuum model, the dual porosity model has severe limitations in accuracy when simulating disconnected fractured media, especially when the porous media contains several large discrete fractures that dominates the flow.} {In order to accurately account for the effect of individual fractures, a number of discrete fracture-matrix models was developed.}

\begin{figure}[!htbp]
\subfigure[Locally refined mesh]{\includegraphics[width = 3in]{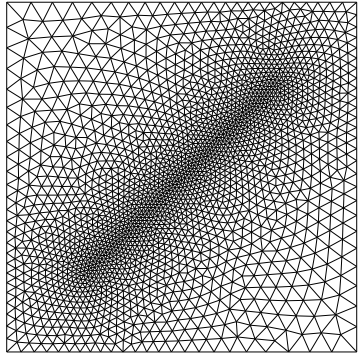}}
\subfigure[Conforming mesh]{\includegraphics[width = 3in]{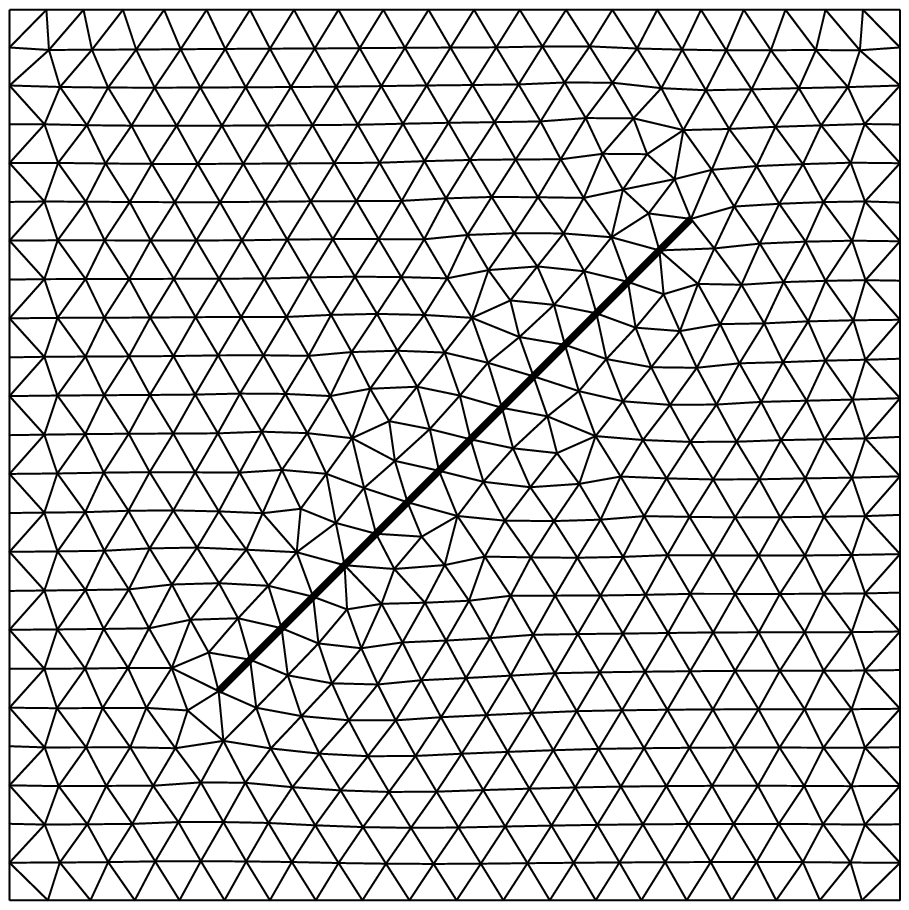}}\\
\subfigure[Non-conforming mesh 1]{\includegraphics[width = 3in]{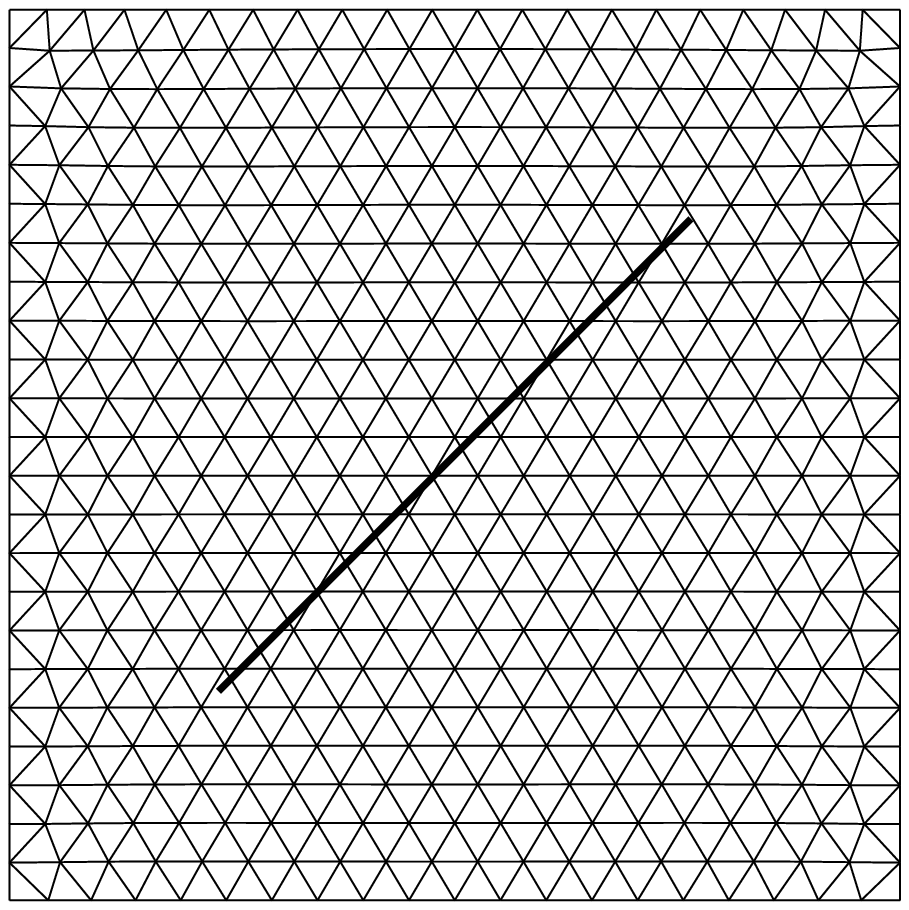}}
\subfigure[Non-conforming mesh 2]{\includegraphics[width = 3in]{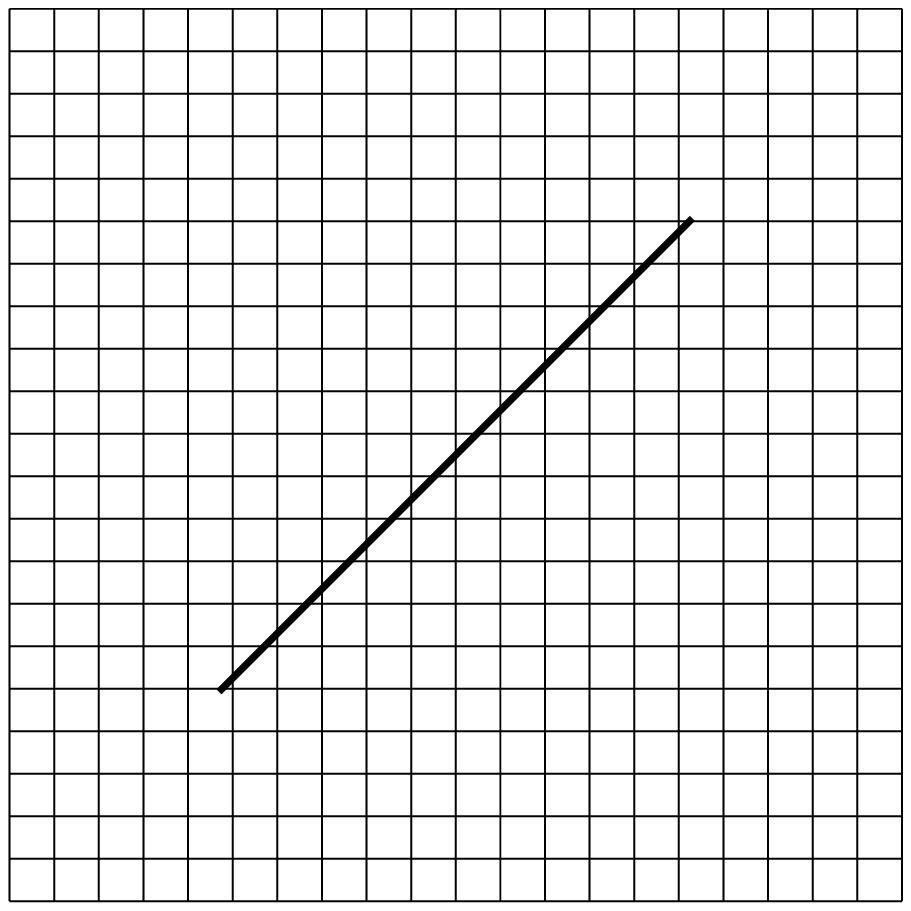}}
\caption{Different type of meshes on fractured media, triangulated by DistMesh\cite{DistMesh} \label{fig:meshes}}
\end{figure}

{The single-porosity model, see e.g. \cite{SinglePoro1}, came with the idea that fractured media is a particular case of heterogeneous media. It precisely capture and describe the fractures by means of local grid refinement in fractured region, see Figure \ref{fig:meshes} (a). However, though providing enough accuracy, the single-porosity model is not practical in real reservoir simulations because of the computational cost resulted from the enormous number of grids due to the scale difference between matrix and fracture thickness.
Therefore, the single-porosity model is only employed to give a reference solution to compare with other algorithms in most of literature, in which it also called equi-dimensional model sometimes.}

To accurately account for the effect of individual fractures without loss of efficiency, the DFM was proposed which treats fractures explicitly as low dimensional entries on the interfaces of high dimensional matrix cells.
The DFM uses the so-called hybrid-dimensional approach to process the flow in 1D fractures and 2D matrix respectively and then couple them together based on the principle of superposition with the fracture system multiplied by the fracture thickness as the dimensional homogeneity factor. By doing so, the DFM avoid the grid refinement in fractured regions to save the efficiency but keeps the accuracy meanwhile.
In 1982, Noorishad and Mehran \cite{FEMDFM1} proposed the first DFM approach in a convection-diffusion problem for single-phase flow. In their work, the conforming mesh was aligned with fractures, in which matrix was discretized by quadrilateral bilinear isoparametric elements and fractures were discretized by one-dimensional line elements as the edge of quadrilateral elements. An upstream weighted residual finite element method was adopted to discretize the transport equations on matrix and fractures, respectively, and then the two equation systems were coupled together using the aforementioned technique.
Later, Baca et al. \cite{FEMDFM2} considered the heat and solute transport in fractured media on conforming meshes where the matrix was discretized by isoparametric elements and the fractures were discretized by line elements along the sides of isoparametric elements, and used the principle of superposition to couple the governing equations for each element type together.
Since then, the DFM has been rapidly developed.
Kim and Deo \cite{SuperposeStiffnessMat,FEMDFM} employed the Galerkin finite element method to discretize the multi-phase flow in matrix on triangular meshes and fractures on its interfaces respectively, and superposed the fracture stiffness matrices on rock stiffness matrices. Then the authors applied an inexact Newton’s method for solving the resulting fully-implicit scheme.
Karimi-Fard and Firoozabadi \cite{DFMpaper1} also used Galerkin method in DFM to solve the two-phase flow problem on a conforming triangular mesh with implicit pressure–explicit saturation (IMPES) time discretization coupled with adaptive time step. The results have shown great agreement with the reference solution of single-porosity model on fine meshes. It is worthy mentioning that the authors explained the methodology of DFM as a decomposition of integration regions for matrix and fracture based on conforming mesh, i.e. $\displaystyle\int_{\Omega}FEQ~d\Omega = \int_{\Omega_m}FEQ~d\Omega_m + \epsilon\int_{\Omega_f}FEQ~d\Omega_f$, where $\epsilon$ is the fracture thickness, $FEQ$ is the flow equation, and $\Omega_m,\Omega_f$ are the regions of 2D matrix and 1D fractures, i.e. $\Omega = \Omega_m + \epsilon\Omega_f$. More DFM based on finite element methods can be find in \cite{FEMDFM3,FEMDFM4}, etc.
In addition to finite element methods, researchers also adopted the idea of DFM on finite volume methods for the purpose of local mass conservative.
Based on conforming Delaunay triangulation where fractures lay on the edges of triangles, the vertex-centered finite volume DFM (Box-DFM) \cite{BoxDFM1,BoxDFM2,BoxDFM3,BoxDFM4,BoxDFM5,BoxDFM6} associates unknowns to each vertex-centered control volumes (CV) which are the dual cells of Delaunay mesh formed by connecting the barycenters of neighboring Delaunay triangles and midpoints of edge around each vertices. In this method, the flux on each CV faces is composed of the rock matrix flux and fracture flux, where the flux on fractures is obtained by multiplying the fracture thickness with the flux on 1D fracture entries.
Another variation of finite volume DFM is the cell-centered finite volume DFM (CC-DFM) \cite{IntersectingFractures,CCDFM2,CCDFM3,CCDFM4,CCDFM5}, where the degrees of freedom are assigned to each triangular matrix elements and the line fracture elements on the edge of triangles. In this model, the transmissibility of different element-adjacent types are determined by the harmonic average of transmissibility of different elements. Depends on the anisotropy and heterogeneity of porous matrix, the fluxes can be approximated by the information of two points (TPFA) or multiple points (MPFA). In \cite{IntersectingFractures}, the authors proposed a widely recognized simplification for flux exchange in intersecting fractures in CC-DFM to remove the small element caused by the intersection.
Moreover, Firoozabadi et al. \cite{MFEMDGDFM1,MFEMDGDFM2,MFEMDGDFM3,MFEMDGDFM4,MFEMDGDFM6,MFEMDGDFM7,MFEMDGDFM8,MFEMDGDFM9}. combined the mixed finite element (MFE) and the discontinuous Galerkin (DG) methods in discrete fracture model to attain the local-conservativity of mass and high accuracy of spices in flow and transports in fractured media. The slope limiter is often cooperated with DG to gained a better stability. There are also conforming approaches based on DG methods\cite{DGinterfaces,DGinterfaces2} that model the fracture as a low dimensional interface imposed with suitable jump conditions of pressure and flux. Besides the above, the mortar-type methods (Mortar-DFM) \cite{MortarDFM1,MortarDFM2} and mimetic finite difference methods (MFD-DFM) \cite{MFDDFM1} are also well investigated in recent researches.

\begin{figure}[!htbp]
\subfigure[Poor mesh 1]{\includegraphics[width = 3in]{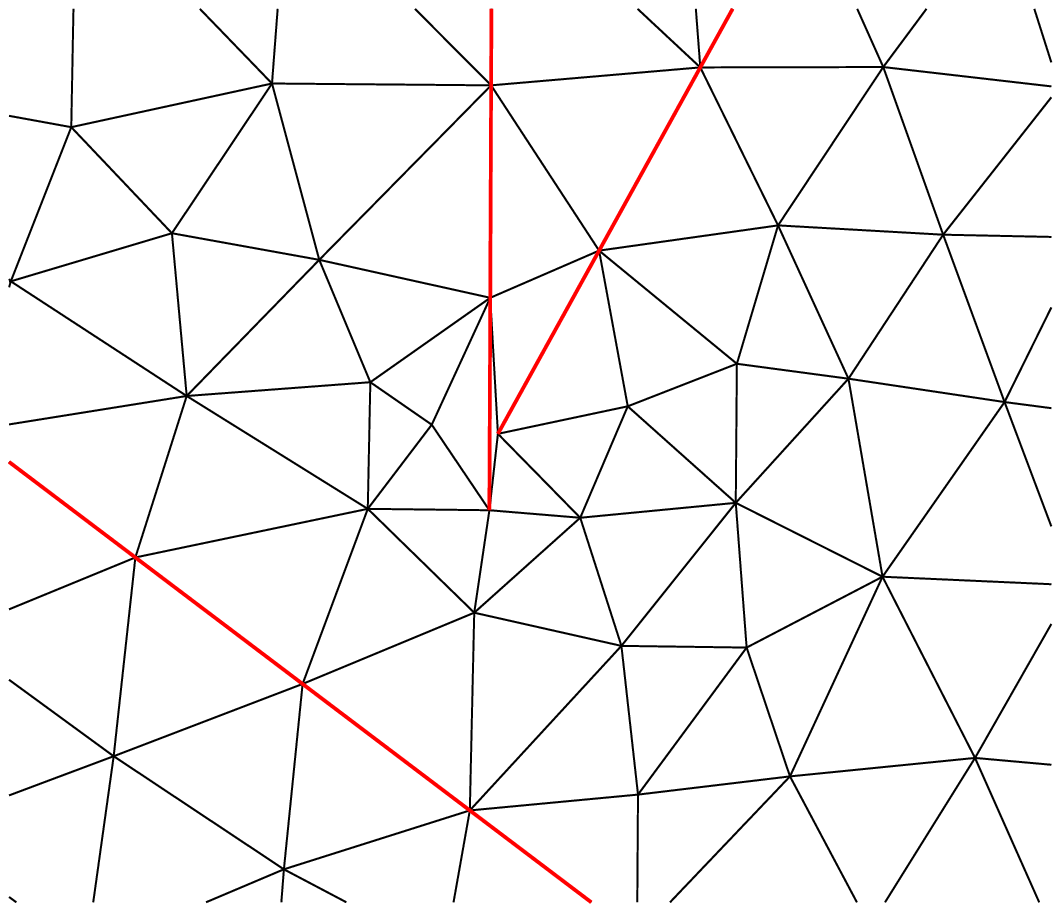}}
\subfigure[Poor mesh 2]{\includegraphics[width = 3in]{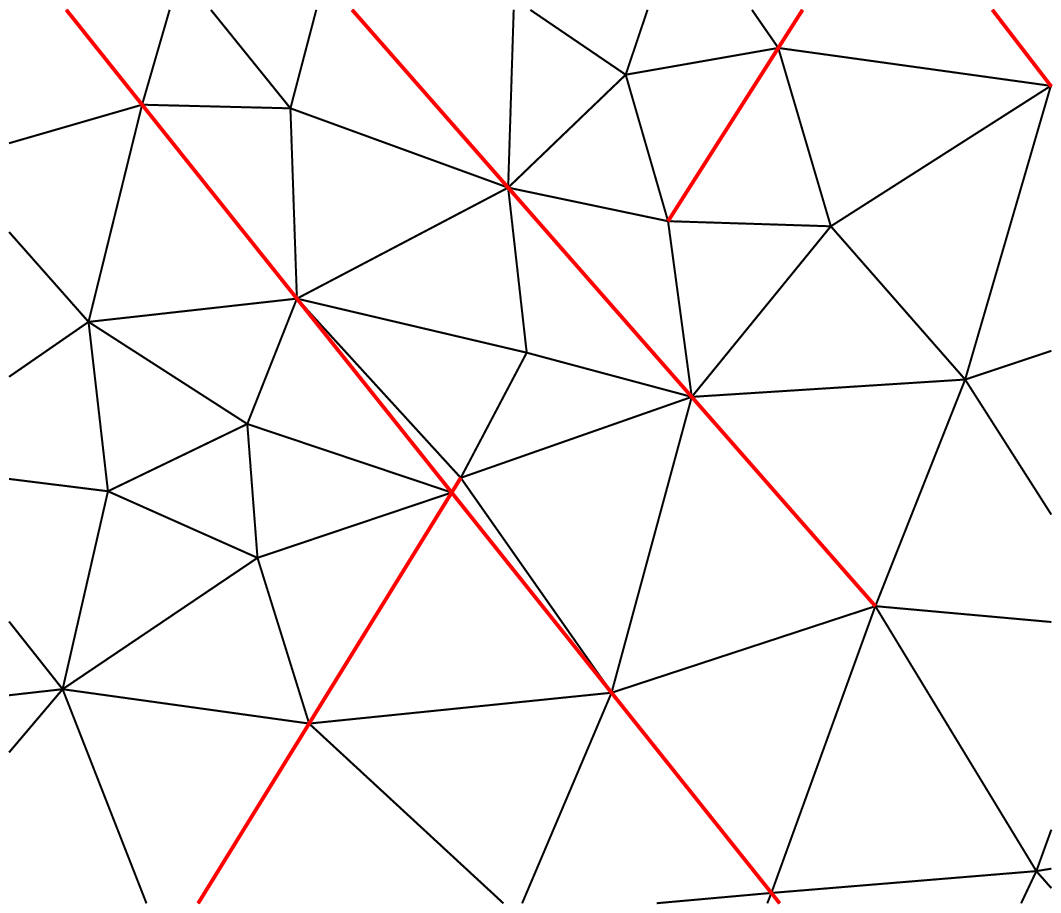}}\\
\subfigure[Poor mesh 3]{\includegraphics[width = 3in]{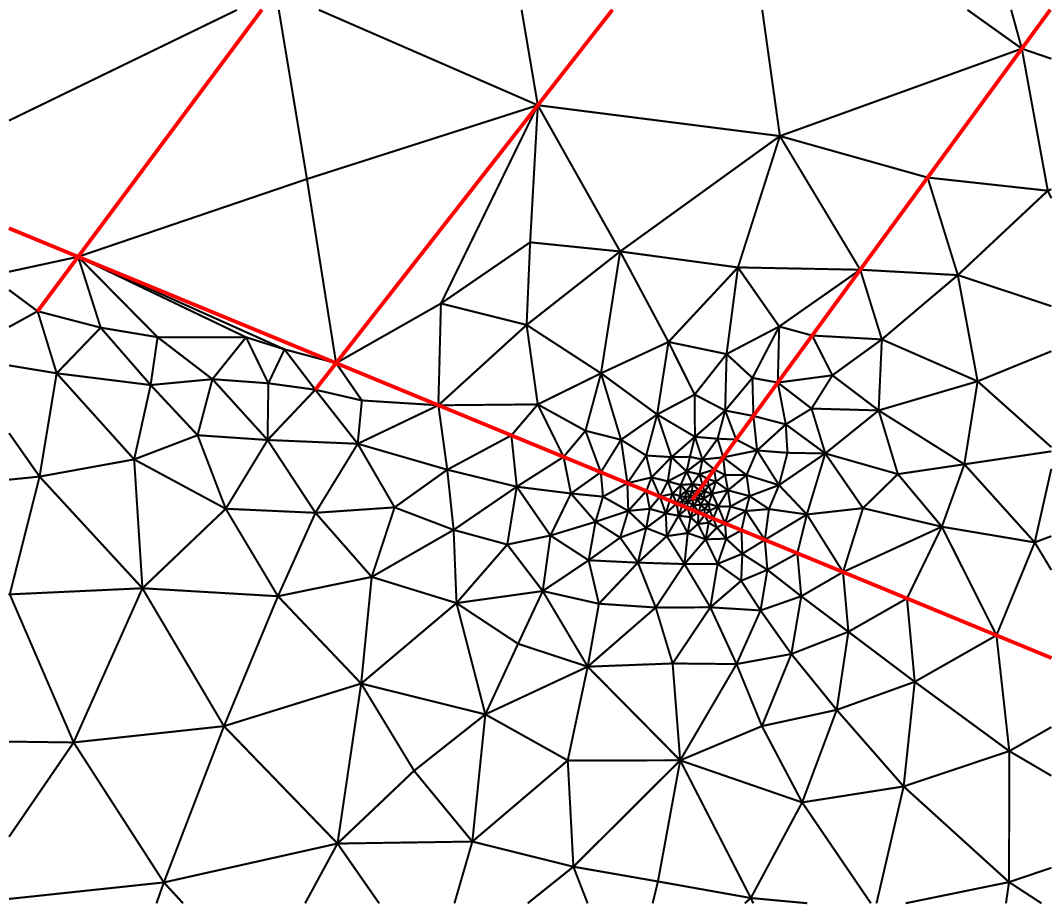}}
\subfigure[Poor mesh 4]{\includegraphics[width = 3in]{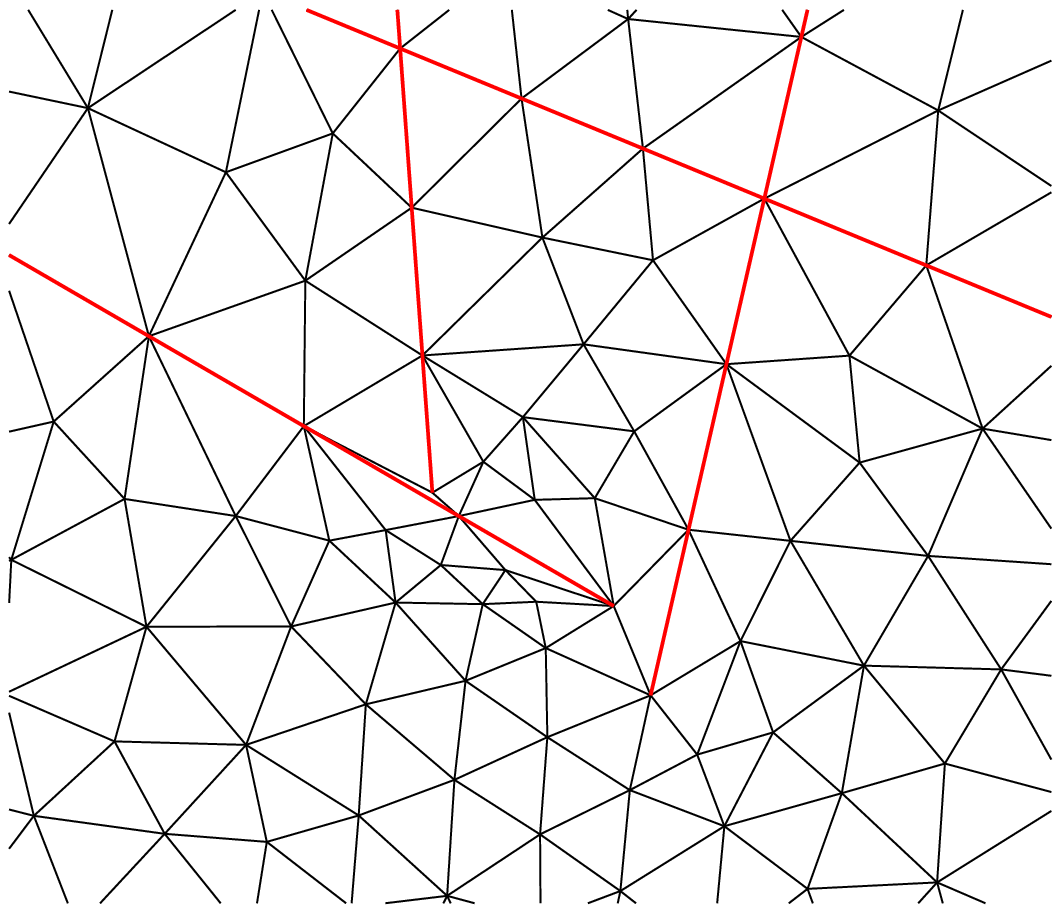}}
\caption[Caption for LOF]{Conforming mesh with poor quality, data from \cite{Website} \label{fig:poormeshes}}
\end{figure}

{However, all the aforementioned DFMs suffer from the limitation of conforming meshes, which is usually composed of triangular element for matrix and line entries as the edge of triangles for fractures, see Figure \ref{fig:meshes} (b). The resulting difficulty is the generation of meshes with high quality, especially when fracture network is complex and the distance or angle between fractures is small, see Figure \ref{fig:poormeshes}. To overcome this shortcoming, lots of efforts were made on finding alternative methods in non-conforming meshes.}

The embedded discrete fracture model (EDFM) is a successful alternative model that relieves the constraint on meshes but keeps explicit description for individual fractures meanwhile. In 2008, Li and Lee \cite{firstEDFM} proposed the first EDFM for black oil model in fractured media. Later it was adopted by Moinfar et al. in \cite{secondEDFM} as well as lots of other follow up works \cite{EDFM3,pEDFM,EDFM4,CrossShaped,EDFM5}.
Borrowing the idea of mass transfer from dual-porosity model, the EDFM computes the flow and transport in 2D matrix and 1D fractures network respectively with accounting for matrix-fracture and fracture-fracture mass transfer based on their pressure differences and geometric and physical information. A typical grids of EDFM with fracture is shown in Figure \ref{fig:meshes} (d), where the fractures are cut into small pieces by matrix grids with each matrix cell and fracture pieces associated by one degree of freedom. Although EDFM only needs structured grid, this method requires computing the average normal distance between matrix grid block and embedded fracture piece to calculate the fluid transport between them and considering the mass transfer between intersecting fractures.

Recently, a non-conforming finite element method\cite{LMFEM2D,LMFEM2D2, LMFEM3D} based on Lagrange multipliers is proposed. Originating from the fictitious domain method\cite{FictitiousDomain}, this approach couples the flow in $(n-1)$-dimensional fractures with that in the $n$-dimensional matrix by means of Lagrange multiplier.
The Lagrange multiplier is also used to impose the continuity of pressure across fractures in this approach. Since the integration forms of the matrix flow, fracture flow and Lagrange multiplier show up in the formulation separately, the meshes for these three terms can be mutually independent.

Other non-conforming methods like extended finite element discrete fracture model (XFEM-DFM) \cite{XFEMDFM1,XFEMDFM2,ThesisXFEM,XFEMDFM3,XFEMDFM4} based on interfaces models \cite{Interfaces1,Interfaces2,Interfaces3,Interfaces4,Interfaces5,Interfaces6,benchmark2} were also proposed in recent years. However, these methods are not as widely practiced as the EDFM in industry because of the difficulty of implementation issue when fracture network is of high geometrical complexity \cite{Benchmark}. Moreover, the CutFEM\cite{CutFEM}, which couples the fluid flow in all lower dimensional manifolds, is another alternative non-conforming method. But this method requires the fractures to cut the domain into completely disjoint subdomains, thus it's not applicable for all fractured media.

In this paper, instead of proposing a new non-conforming alternative of DFM, we clarify that the primal finite element discrete fracture model indeed is not really restricted on conforming meshes. Actually, we can reinterpret the DFM in a different way then the DFM can be applied to nonconforming meshes. The inspiration comes from the comb model \cite{DeltaPaper}, which uses a diffusion tensor containing Dirac-$\delta$ function to describe a special diffusion process that occurs only on $x$-axis in X-direction but on whole plane in Y-direction. Drawing from the comb model, we propose a hybrid-dimensional representation of permeability tensor of fractured media and obtain the non-conforming DFM scheme by directly applying Galerkin method on this model. The derivation is quite simple but helps us get rid of the restriction on conforming meshes. Typically applicable meshes for this scheme are shown in Figure \ref{fig:meshes} (c), (d).

Though both address the matrix-fracture mesh non-conformity, our model essentially differentiates from the EDFM. The EDFM borrows the idea from dual-porosity model and focus on the calculation of the transmissibility factor of different non-neighboring connections (NNCs) thus the fractures and matrix flows are two different systems, while our model based on the idea of representing fractured media as hybrid-dimensional permeability tensor thus the fracture and matrix flows are one system. Due to their differences, our model have some advantages compared with EDFM. First, the degrees of freedom and complexity of our model won't increase as the number of intersecting fractures increases, while in EDFM, the system will become more complicated if there are more intersecting fractures (resulting in more NNCs), especially when they intersect at one point. Second, our model can handle curved fractures naturally, which is shown in later sections. Moreover, since we don't need to compute the transmissibility factor used in EDFM, the costs on geometric computation of the average normal distance from matrix to fractures and from fractures to fractures are relieved. We also would like to point out some shortcomings and limitations of the method in this paper. Not like some adapted EDFM\cite{pEDFM}, the model is unable to treat the barriers directly. The pressure jump across the fracture is not representable in the model as well. In addition, due to the property of the finite element method, the approach proposed in the paper is not locally mass conservative. However, it is possible to apply the idea introduced in \cite{FEMconservation} to obtain the local mass conservation.

It's notable that the out look of the formulation of our method is more or less similar to that of Lagrange multiplier approach\cite{LMFEM2D,LMFEM2D2, LMFEM3D} but they are in totally different theoretical frameworks. One may refer to that approach if of interest.

{The rest of this paper is organized as follows. In Section 2, we introduce the equi-dimensional model problem of the steady-state single-phase flow in fractured porous media. In Section 3, we adopt the comb model to give a hybrid-dimensional representation for permeability tensor of fractured media. In Section 4, we apply the standard Galerkin finite element method to the model proposed in Section 3 to attain the DFM scheme on non-conforming meshes. The effectiveness, accuracy and consistency with traditional finite element DFM on conforming meshes of this scheme are demonstrated in Section 5 with plenty of numerical tests. Finally, we end in Section 6 with some concluding remarks.}

\section{Equi-dimensional Model for single-phase flow in fractured media}
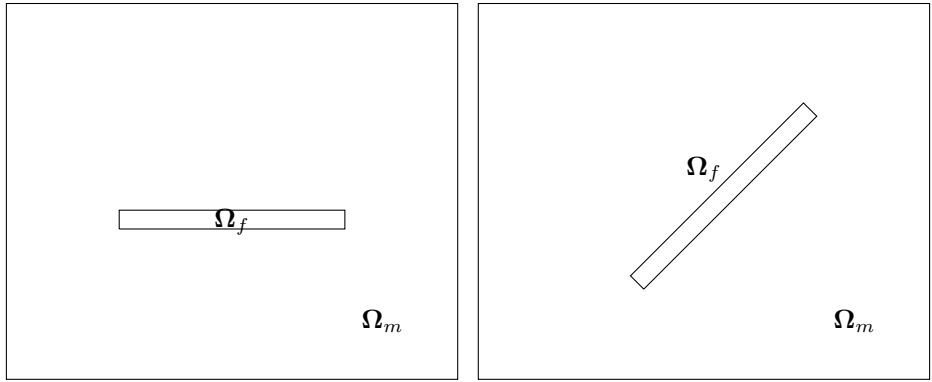
\begin{figure}[!htbp]
\centering
\subfigure[horizontal fracture]
{
\begin{tikzpicture}[scale=1]
\draw (-3,-2.5) -- (3,-2.5) -- (3,2.5) --(-3,2.5) -- (-3,-2.5) -- cycle;
\coordinate [label=$\bm{\Omega}_m$] (h) at(2,-2) ;

\draw (-1.5,-0.5) -- (1.5,-0.5) -- (1.5,-0.25) -- (-1.5,-0.25) -- (-1.5,-0.5);
\coordinate [label=$\bm{\Omega}_f$] (h) at(-0,-0.7) ;

\end{tikzpicture}
}
\subfigure[oblique fracture]
{
\begin{tikzpicture}[scale=1]
\draw (-3,-2.5) -- (3,-2.5) -- (3,2.5) --(-3,2.5) -- (-3,-2.5) -- cycle;
\coordinate [label=$\bm{\Omega}_m$] (h) at(2,-2) ;
\draw (-0.8,-1.3) -- (1.5,1.0) -- (1.32322,1.17678) -- (-0.97678,-1.12322) -- (-0.8,-1.3);
\coordinate [label=$\bm{\Omega}_f$] (h) at(0,-0) ;
\end{tikzpicture}
}
\caption{Fractured media in model problem}
\label{fig:Fractures2D}
\end{figure}

For the steady-state single-phase flow, the distribution of pressure $p$ in heterogeneous porous media $\Omega$ is governed by the Poisson's equation:
\begin{equation}\label{PoissonEq}
    -\nabla\cdot(\textbf{K}\nabla p)=f,\quad x\in\Omega,
\end{equation}
where ${\bf{K}}$ is the permeability tensor of porous media and $f$ is the source term. For fractured porous media, $\bf{K}$ can be expressed as follows: \\
\begin{equation}\label{EquiDimensionK}
{\bf{K}} =
 \begin{cases}
      {\bf K}_m & x\in \Omega_m \\
      {\bf K}_f & x\in \Omega_f  ,
   \end{cases}
\end{equation}
where $\Omega_m$ and $\Omega_f$ are the regions of porous matrix and fractures respectively, which compose the whole domain $\Omega$. See Figure \ref{fig:Fractures2D} as an illustration, where the thickness of fractures is exaggerated for the sake of visibility.

We consider the mixed boundary condition
\begin{equation}\label{BVP}
    p = p_D,\quad \text{on}~ \Gamma_D\in\partial\Omega,\quad \text{and} \quad -({\bf K}\nabla p) \cdot {\bf n}=q_N, \quad \text{on}~ \Gamma_N = \partial\Omega\setminus\Gamma_D,
\end{equation}
where ${\bf n}$ is the unit outer normal vector of the boundary $\partial \Omega$.

The model problem \eqref{PoissonEq}, \eqref{EquiDimensionK} and \eqref{BVP} is usually called the equi-dimensional model for single-phase flow in fractured media and the expression $\textbf{K}$ in \eqref{EquiDimensionK} is the equi-dimensional representation of the permeability tensor.

\section{Hybrid-dimensional representation of permeability tensor $\textbf{K}$}
Drawing from the comb model \cite{DeltaPaper} which used a diffusion tensor containing Dirac-$\delta$ function to describe a special diffusion process that
occurs only on $x$-axis in X-direction but on whole plane in Y-direction, this section establishes the hybrid-dimensional representation of permeability tensor of the fractured media.

\begin{figure}[!htbp]
\centering
\subfigure[horizontal fracture]
{
\begin{tikzpicture}[scale=1]
\draw (-3,-2.5) -- (3,-2.5) -- (3,2.5) --(-3,2.5) -- (-3,-2.5) -- cycle;
\draw[black, ultra thick] (-1.5,-0.5) -- (1.5,-0.5);
\coordinate [label=$\bm{l}$] (h) at(0,-0.4) ;

\draw[-latex] (0,0.5) -- (1,0.5);
\coordinate [label=$\bm{\nu}$] (h) at(1,0.5) ;
\draw[-latex] (0,0.5) -- (0,1.5);
\coordinate [label=$\bm{\sigma}$] (h) at(0,1.5) ;

\draw [thin,dash dot] (-1.5,-0.5) -- (-1.5,-3);
\coordinate [label=$\xi_1$] (h) at(-1.7,-3) ;
\draw [thin,dash dot] (1.5,-0.5) -- (1.5,-3);
\coordinate [label=$\xi_2$] (h) at(1.7,-3) ;
\draw [thin,dash dot] (-1.5,-0.5) -- (-2.2,-0.5);
\coordinate [label=$\eta_0$] (h) at(-2,-0.5) ;

\draw[-latex] (-2.2,-3) -- (3,-3);
\coordinate [label=$\mathbf{x/\xi}$] (h) at(3,-3) ;
\draw[-latex] (-2.2,-3) -- (-2.2,0.8);
\coordinate [label=$\mathbf{y/\eta}$] (h) at(-2,0.8) ;
\coordinate [label=$\bm{\Omega}$] (h) at(2,-2) ;

\coordinate [label=$O$] (h) at(-2.2,-3.5) ;

\end{tikzpicture}
}
\subfigure[oblique fracture]
{
\begin{tikzpicture}[scale=1]
\draw (-3,-2.5) -- (3,-2.5) -- (3,2.5) --(-3,2.5) -- (-3,-2.5) -- cycle;

\coordinate [label=$\bm{\Omega}$] (h) at(2,-2) ;

\draw[black, ultra thick] (-0.8,-1.3) -- (1.5,1.0);
\coordinate [label=$\bm{l}$] (h) at(0.3,0.0) ;

\draw[-latex] (0,0.5) -- (0.75,1.25);
\coordinate [label=$\bm{\nu}$] (h) at(0.75,1.25) ;
\draw[-latex] (0,0.5) -- (-0.75,1.25);
\coordinate [label=$\bm{\sigma}$] (h) at(-0.75,1.25) ;

\draw[-latex] (-1,-3) -- (3.5,1.5);
\coordinate [label=$\mathbf{\xi}$] (h) at(3.5,1.5) ;
\draw[-latex] (-1,-3) -- (-4,0);
\coordinate [label=$\mathbf{\eta}$] (h) at(-4,0) ;

\draw[-latex] (-1,-3) -- (3.5,-3);
\coordinate [label=$\mathbf{x}$] (h) at(3.5,-3) ;
\draw[-latex] (-1,-3) -- (-1,0);
\coordinate [label=$\mathbf{y}$] (h) at(-1,0) ;

\coordinate [label=$O$] (h) at(-1,-3.5) ;
\draw[<->] (-0.2,-3) arc (0:30:1.25);
\coordinate [label=left:$\theta$] (a) at (-0.2,-2.75);

\draw [thin,dash dot] (-0.8,-1.3) -- (-0.05,-2.05);
\coordinate [label=$\xi_1$] (h) at(-0.25,-2.25) ;
\draw [thin,dash dot] (1.5,1.0) -- (2.25,0.25);
\coordinate [label=$\xi_2$] (h) at(2.25,0.25) ;
\draw [thin,dash dot] (-0.8,-1.3) -- (-1.75,-2.25);
\coordinate [label=$\eta_0$] (h) at(-1.75,-2.25) ;
\end{tikzpicture}
}
\caption{Fractured media and the corresponding coordinates systems}
\label{fig:Fractures1D}
\end{figure}
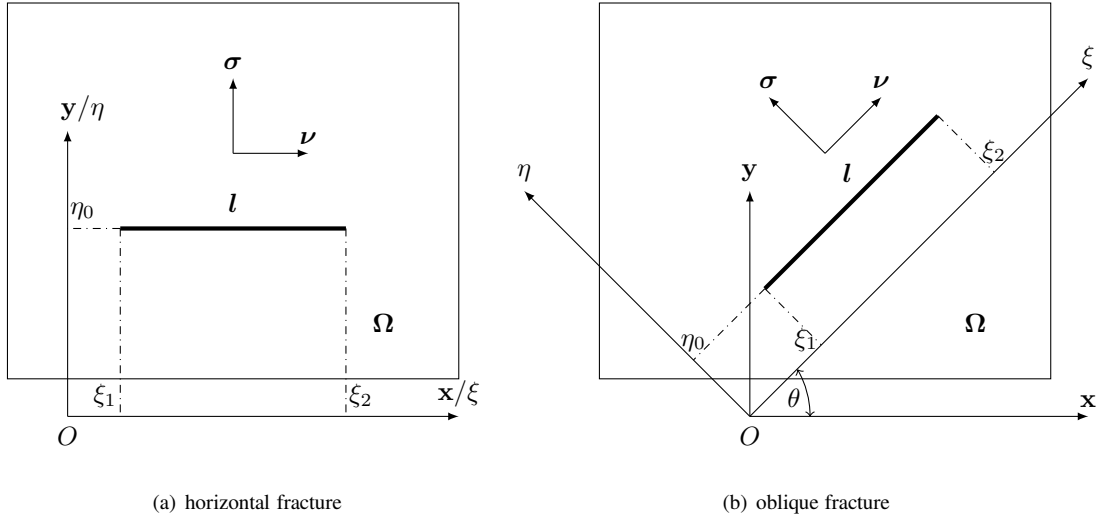

We follow the treatment in hybrid-dimensional discrete fracture model to reduce the fractures from two-dimensional narrow strips in Figure \ref{fig:Fractures2D} to one-dimensional line segments $\bm{l}$, see Figure \ref{fig:Fractures1D}. Moreover, we establish the global coordinates system $xOy$ and local coordinates system $\xi O\eta$  associated with the fracture. We denote by ${\bm{\nu}}$ and ${\bm{\sigma}}$ the tangential and normal unit vector of fracture $\bm{l}$, respectively. See Figure \ref{fig:Fractures1D} as an illustration.
For all $P\in \Omega$, the transformation between local coordinates $(\xi_P,\eta_P)$ and global coordinates $(x_P,y_P)$ is given as
\begin{equation}\label{transform}
\begin{bmatrix}
    \xi_P \\
    \eta_P
\end{bmatrix}=
\begin{bmatrix}
    \cos(\theta) & \sin(\theta) \\
    -\sin(\theta) & \cos(\theta)
\end{bmatrix}
\begin{bmatrix}
    x_P \\
    y_P
\end{bmatrix},
\end{equation}
where $\theta$ is the angle of the fracture.

According to the principle of superposition, the permeability tensor of the whole domain is composed of the permeability tensors of matrix and fracture, i.e. $\textbf{K}=\textbf{K}_m+\textbf{K}_f$. Due to the geometry of the fracture, the permeability tensor $\textbf{K}_f$ must be symmetric with its two characteristic directions being the tangential and normal directions of the fracture, namely $\bm{\nu}$ and $\bm{\sigma}$, respectively. Therefor, $\textbf{K}_f={K_{f,t}}\bm{\nu}\bm{\nu}^{T} + {K_{f,n}}\bm{\sigma}\bm{\sigma}^{T}$ by the spectral decomposition theorem.

Suppose the original thickness of the fracture strip in Figure \ref{fig:Fractures2D} is $\epsilon$ and the tangential permeability is $k_f$, and note that the gradient of pressure can be decomposed into $\displaystyle\nabla p=\frac{\partial p}{\partial {\nu}}\bm{\nu}+\frac{\partial p}{\partial {\sigma}}\bm{\sigma}$.

First, we consider the tangential component of the gradient of pressure $\displaystyle\frac{\partial p}{\partial {\nu}}\bm{\nu}$. The resulting flow in the original fracture strip is $\displaystyle\epsilon k_f \frac{\partial p}{\partial {\nu}}$ by Darcy's law, from which we can see $\epsilon k_f$ measures the conductivity of the fracture. In our model, the conductivity of the fracture is concentrated on this line segment and a line Dirac-$\delta$ function will be the choice to represent the permeability. Moreover, we should take the location of the fracture into account. Hence we have
\begin{equation}\label{localK}
K_{f,t}=\epsilon k_f\delta(\eta-\eta_0)\mathbbm{1}(\xi_1\leq \xi \leq \xi_2),
\end{equation}
where $\mathbbm{1}(\cdot)$ is the indicator function defined as $\mathbbm{1}(\text{expr})$ equals $1$ if expr is true while equals $0$ otherwise. By the coordinates transformation \eqref{transform}, the expression of $K_{f,t}$ under global coordinates $(x,y)$ is
\begin{equation}\label{globalK}
K_{f,t}=\epsilon k_f \delta(-\sin(\theta)x+\cos(\theta)y-\eta_0)\mathbbm{1}(\xi_1\leq \cos(\theta)x+\sin(\theta)y\leq \xi_2).
\end{equation}

Second, we consider the normal component of the gradient of pressure $\displaystyle\frac{\partial p}{\partial {\sigma}}\bm{\sigma}$. The effect of fracture on flow in this direction is negligible since its thickness $\epsilon$ is so tiny.

Therefore, instead of \eqref{EquiDimensionK}, we have the following hybrid-dimensional expression for permeability tensor in single fractured porous media
\begin{equation}\label{HybridDimensionK}
{\bf{K}} = \textbf{K}_m +  \epsilon k_f\delta(\cdot)\mathbbm{1}(\cdot)\bm{\nu\nu}^{T},
\end{equation}
where $\delta(\cdot)\mathbbm{1}(\cdot)$ is the shorthand of its full expression in \eqref{localK} and \eqref{globalK} under the local and global coordinates systems, respectively. Note that $\epsilon k_f$ measures the conductivity of the fracture, $\delta(\cdot)\mathbbm{1}(\cdot)$ contains the information of position of the fracture and  $\bm{\nu}$ indicates the direction of the fracture.

The expression above is only for fractured media with single fracture but can be extended to a fracture network as
\begin{equation}
\displaystyle{\bf{K}} = \textbf{K}_m + \sum^{L}_{i=1} \epsilon_i k_{fi}\delta_i(\cdot)\mathbbm{1}_i(\cdot)\bm{\nu}_i\bm{\nu}_i^{T},
\end{equation}
where $L$ is the number of fractures.

We call the expression $\textbf{K}$ in \eqref{HybridDimensionK} the hybrid-dimensional representation of permeability tensor because the matrix part $\textbf{K}_m$ in \eqref{HybridDimensionK} is of rank 2 and has a 2D support while the fracture part $\epsilon k_f\delta(\cdot)\mathbbm{1}(\cdot)\bm{\nu\nu}^{T}$ is of rank 1 and has a 1D support. Analogous to Section 2, we call the model problem \eqref{PoissonEq}, \eqref{HybridDimensionK} and \eqref{BVP} the hybrid-dimensional model for single-phase flow in fractured media.

\section{Finite element DFM scheme on non-conforming meshes}
This section explores the non-conforming DFM and analyze its relationship with traditional finite element DFM \cite{SuperposeStiffnessMat}. We first establish the variational form of the hybrid-dimensional model in Section \ref{Variational}. Then, in Section \ref{Numerical}, we construct the non-conforming DFM scheme by using linear Lagrange basis functions as finite element space in the variational form.
\subsection{Variational form of the hybrid-dimensional model}\label{Variational}
We define the variational space
$$H^1_D := \{v\in H^{1}(\Omega): v|_{\Gamma_D}=p_D\},$$ and
$$H^1_0 := \{v\in H^{1}(\Omega): v|_{\Gamma_D}=0\}.$$

Then the variational form of \eqref{PoissonEq}, \eqref{BVP} is to find a $p\in H^{1}_{D}$, such that the following variational equation holds for all $v\in H^{1}_{0}$,
\begin{equation}\label{variation1}
    \int_{\Omega}(\textbf{K}\nabla p)\cdot\nabla v ~dx dy = \int_{\Omega}f v~dx dy+\int_{\Gamma_N}q_N v~ds.
\end{equation}

Using the hybrid-dimensional representation of $\textbf{K}$ established in \eqref{HybridDimensionK}, we have $(\textbf{K}\nabla p)\cdot\nabla v = (\textbf{K}_m\nabla p)\cdot\nabla v + \epsilon k_f\delta(\cdot)\mathbbm{1}(\cdot)(\bm{\nu\nu}^{T}\nabla p)\cdot\nabla v$. The second term of right hand side can be further simplified as $\displaystyle\epsilon k_f\delta(\cdot)\mathbbm{1}(\cdot)(\bm{\nu\nu}^{T}\nabla p)\cdot\nabla v = \epsilon k_f\delta(\cdot)\mathbbm{1}(\cdot)\frac{\partial p}{\partial \nu}\bm{\nu}\cdot\nabla v = \epsilon k_f\delta(\cdot)\mathbbm{1}(\cdot)\frac{\partial p}{\partial\nu}\frac{\partial v}{\partial\nu}$.

Therefore, one can get the equivalent variational form as follows,
\begin{equation}\label{variation2}
    \int_{\Omega}(\textbf{K}_m\nabla p)\cdot\nabla v ~dx dy + \int_{\Omega}\epsilon k_f\delta(\cdot)\mathbbm{1}(\cdot)\frac{\partial p}{\partial\nu}\frac{\partial v}{\partial\nu} ~dx dy = \int_{\Omega}f v~dx dy+\int_{\Gamma_N}q_N v~ds,
\end{equation}

The following lemma helps us to rewrite the integration of the fracture term as a line integral on fracture.
\begin{lemma}\label{lemmadelta}
Let $\delta(\cdot)\mathbbm{1}(\cdot)$ be the shorthand of its full expression in \eqref{globalK}. For any continuous function $g$ on $\Omega$, we have
\begin{equation}
\int_{\Omega} \delta(\cdot)\mathbbm{1}(\cdot)g(x,y)~dx dy =\int_{l} g(x,y)~ ds,
\end{equation}
where the line segment ${\bm{l}}$ is the support of $\delta(\cdot)\mathbbm{1}(\cdot)$ as shown in Figure \ref{fig:Fractures1D}.
\end{lemma}

\begin{proof}
\begin{align*}
&\int_{\Omega} \delta(\cdot)\mathbbm{1}(\cdot)g(x,y)~dx dy\\
:=&\int_{\Omega} \delta(-\sin(\theta)x+\cos(\theta)y-\eta_0)\mathbbm{1}(\xi_1\leq \cos(\theta)x+\sin(\theta)y\leq \xi_2)g(x,y)~dx dy\\
=&\int_{x}\left(\int_{y} \delta(-\sin(\theta)x+\cos(\theta)y-\eta_0)\mathbbm{1}(\xi_1\leq \cos(\theta)x+\sin(\theta)y\leq \xi_2)g(x,y)~dy\right) dx\\
& ~ \text{Let}~ t=\cos(\theta)y:\\
=&\int_{x}\left(\sec(\theta)\int_{t} \delta(t-\sin(\theta)x-\eta_0)\mathbbm{1}(\xi_1\leq \cos(\theta)x+\tan(\theta)t\leq \xi_2)g(x,\sec(\theta)t)~dt\right) dx \\
=&\int_{x}\sec(\theta)\mathbbm{1}(\xi_1\leq \cos(\theta)x+\tan(\theta)(\sin(\theta)x+\eta_0)\leq \xi_2)g(x,\tan(\theta)x+\sec(\theta)\eta_0)~ dx\\
=&\int_{x}\sec(\theta)\mathbbm{1}(\cos(\theta)\xi_1-\sin(\theta)\eta_0\leq x\leq \cos(\theta)\xi_2-\sin(\theta)\eta_0)g(x,\tan(\theta)x+\sec(\theta)\eta_0)~ dx\\
=&\int_{\cos(\theta)\xi_1-\sin(\theta)\eta_0}^{\cos(\theta)\xi_2-\sin(\theta)\eta_0}\sec(\theta)g(x,\tan(\theta)x+\sec(\theta)\eta_0)~ dx\\
& ~ \text{Let}~\xi=\sec(\theta)\left(x+\sin(\theta)\eta_0\right):\\
=&\int_{\xi_1}^{\xi_2}g(\cos(\theta)\xi-\sin(\theta)\eta_0,\sin(\theta)\xi+\cos(\theta)\eta_0)~ d\xi\\
=&\int_{l}g(x,y)~ ds \qquad (\text{By using \eqref{transform} and the definition of line integral of first type })
\end{align*}
We assume $\theta\in(-\frac{\pi}{2},\frac{\pi}{2})$ in the above and the case $\theta=\frac{\pi}{2}$ can be proved by a limit process.
\end{proof}
Based on Lemma \ref{lemmadelta}, we can obtain the final version of the equivalent variational form of hybrid-dimensional model \eqref{PoissonEq}, \eqref{HybridDimensionK} and \eqref{BVP} of fractured porous media:
\begin{equation}\label{variation3}
    \int_{\Omega}(\textbf{K}_m\nabla p)\cdot\nabla v ~dx dy + \int_{\bm l}\epsilon k_f\frac{\partial p}{\partial\nu}\frac{\partial v}{\partial\nu} ~ds = \int_{\Omega}f v~dx dy+\int_{\Gamma_N}q_N v~ds.
\end{equation}
\begin{remark}
In \eqref{variation3}, the shape of fracture ${\bm l}$ actually can be a curve and the thickness $\epsilon$ and tangential permeability $k_f$ of the fracture can be a scalar function defined along ${\bm l}$. In the case of fracture network, we just add all fracture terms together in \eqref{variation3}, i.e.
\begin{equation}
\displaystyle\int_{\Omega}\left(\textbf{K}_m\nabla p\right)\cdot\nabla v ~dx dy + \sum^{L}_{i=1}\int_{{\bm l}_i}\epsilon_i k_{fi}\frac{\partial p}{\partial\nu_i}\frac{\partial v}{\partial\nu_i} ~ds = \int_{\Omega}f v~dx dy+\int_{\Gamma_N}q_N v~ds,
\end{equation}
where $L$ is the number of fractures.
Moreover, we would like to note that the left-hand side of \eqref{variation3} presents in a beautiful adjoint form. The second term of the left-hand side can be viewed as an oriented one-dimensional contribution to the bilinear form.
\end{remark}

\subsection{Numerical scheme of non-conforming DFM}\label{Numerical}

Without loss of generality, we demonstrate the discretization of the variational form \eqref{variation3} on unstructured triangular meshes using linear finite element spaces. One can extend the scheme to rectangular meshes without difficulty.

We adopt the linear finite element space $$V_h=\text{span}\{\Psi_1,\Psi_2,\ldots,\Psi_M, \Psi_{M+1}, \ldots, \Psi_N \},$$
where $\Psi_i$'s are the Lagrange linear basis with Lagrange property, i.e. $\Psi_i({\bm x}_j)=\begin{cases}
      1 ,& ~ i=j \\
      0 ,& ~ i\ne j
   \end{cases}$, in which $\bm{x}_j$'s are vertices of the triangulation, $M$ is the number of non-Dirichlet vertices (degrees of freedom) and $N$ is the total number of vertices in the triangulation.

The aim is to find $\bm{p}=(p_j)$ such that the linear system \eqref{variation4} holds.
\begin{equation}\label{variation4}
\sum^{N}_{j=1}\left\{\int_{\Omega}(\textbf{K}_m\nabla \Psi_i)\cdot\nabla \Psi_j ~dx dy + \int_{\bm l}\epsilon k_f\frac{\partial \Psi_i}{\partial\nu}\frac{\partial \Psi_j}{\partial\nu} ~ds\right\}p_j = \int_{\Omega}f \Psi_i~dx dy+\int_{\Gamma_N}q_N \Psi_i~ds,\quad i=1,2,\ldots,M.
\end{equation}
Note that the Dirichlet boundary condition gives us $p_j=p({\bm x}_j)$ for $\bm{x}_j\in\Gamma_D$, $j=M+1,\ldots,N$.

Denote by
$$a_{ij} = \int_{\Omega}(\textbf{K}_m\nabla \Psi_i)\cdot\nabla \Psi_j ~dx dy + \int_{\bm l}\epsilon k_f\frac{\partial \Psi_i}{\partial\nu}\frac{\partial \Psi_j}{\partial\nu} ~ds,\quad \bm{A}=(a_{ij}),$$ and $$b_i = \int_{\Omega}f \Psi_i~dx dy+\int_{\Gamma_N}q_N \Psi_i~ds,\quad \bm{b}=(b_i),$$
then the linear system can be written as $\bm{Ap}=\bm{b}$.

Like the traditional finite element methods, the global stiffness matrix $\bm{A}$ is obtained by assembling local stiffness matrices element-wisely. We illustrate this procedure in the following, and show the consistency of the non-conforming discrete fracture model (NDFM) with the traditional DFM \cite{SuperposeStiffnessMat} based on that. As for the assembly of $\bm{b}$, one can refer to \cite{Zienkiewicz}.

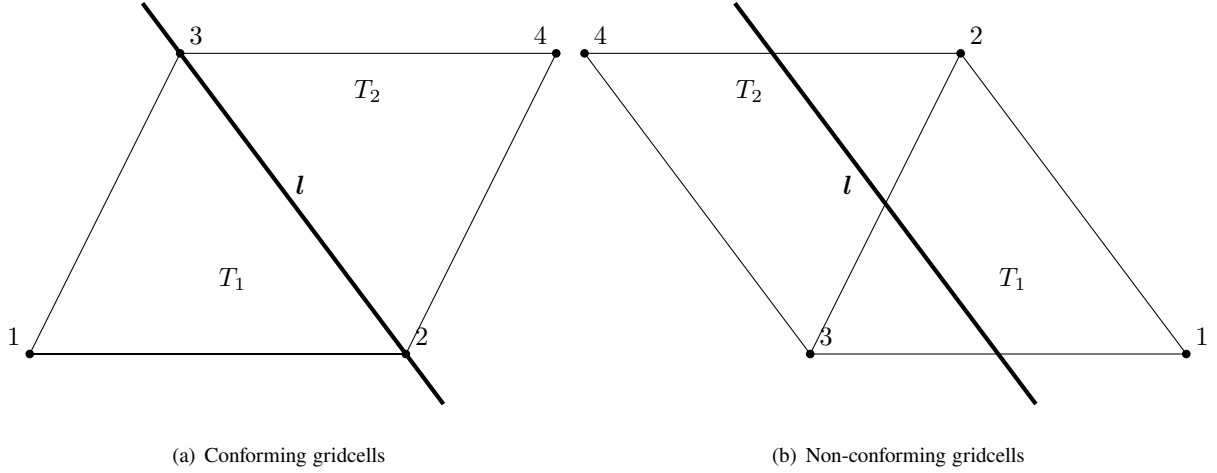
\begin{figure}[!htbp]
\centering
\subfigure[Conforming gridcells]
{
\begin{tikzpicture}[scale=1]
\draw (0,0) -- (5,0) -- (7,4)-- (2,4) -- (5,0) -- cycle;
\draw (0,0) -- (2,4);

\foreach \Point/\PointLabel in {(0,0)/1, (7,4)/4}
\draw[fill=black] \Point circle (0.05) node[above left] {$\PointLabel$};

\foreach \Point/\PointLabel in { (5,0)/2, (2,4)/3}
\draw[fill=black] \Point circle (0.05) node[above right] {$\PointLabel$};
\draw[black, ultra thick] (5.5,-0.66667) -- (1.5,4.66667);
\coordinate [label=$\bm{l}$] (h) at(3.6,2) ;
\coordinate [label=left:$T_1$] (a) at (3,1);
\coordinate [label=left:$T_2$] (b) at (4.8,3.5);
\end{tikzpicture}\label{fig:localconforming}
}
\subfigure[Non-conforming gridcells]
{
\begin{tikzpicture}
\draw (0,0) -- (5,0) -- (2,4) -- (-3,4) -- (0,0) -- cycle;
\draw (0,0) -- (2,4);
\foreach \Point/\PointLabel in {(5,0)/1,(2,4)/2,(0,0)/3, (-3,4)/4}
\draw[fill=black] \Point circle (0.05) node[above right] {$\PointLabel$};
\draw[black, ultra thick] (3,-0.66667) -- (-1,4.66667);
\coordinate [label=$\bm{l}$] (h) at(0.5,2) ;
\coordinate [label=left:$T_1$] (a) at (3,1);
\coordinate [label=left:$T_2$] (b) at (-0.5,3.5);
\end{tikzpicture}\label{fig:localnonconforming}
}
\caption{\label{fig:TwoCases} Locality of gridcells}
\end{figure}

Consider two possible cases in a triangulation of fractured media. A fracture may be conforming with the gridcells locally, see Figure \ref{fig:localconforming}. The other possibility is the case of non-conforming gridcells, as is shown in Figure \ref{fig:localnonconforming}. The global stiffness matrix $\bm{A}$ in this zone is assembled by three local stiffness matrices parts in both cases: elements parts $T_1$, $T_2$, and the fracture part $\bm{l}$.

The local stiffness matrix contributed by triangular elements $T_1$ and $T_2$ are given as follows,
\begin{equation*}
    T_1:\begin{bmatrix}
    a^{T_1}_{11} & a^{T_1}_{12} & a^{T_1}_{13} & 0 \\
    a^{T_1}_{21} & a^{T_1}_{22} & a^{T_1}_{23} & 0 \\
    a^{T_1}_{31} & a^{T_1}_{32} & a^{T_1}_{33} & 0 \\
    0 & 0 & 0 & 0
\end{bmatrix},\qquad
T_2:\begin{bmatrix}
    0 & 0 & 0 & 0 \\
    0 & a^{T_2}_{22} & a^{T_2}_{23} & a^{T_2}_{24} \\
    0 & a^{T_2}_{32} & a^{T_2}_{33} & a^{T_2}_{34} \\
    0 & a^{T_2}_{42} & a^{T_2}_{43} & a^{T_2}_{44},
\end{bmatrix}
\end{equation*}
where
$$a^{T_n}_{i,j}=\int_{T_n}(\textbf{K}_m\nabla\Psi_i)\cdot\nabla\Psi_{j}~dxdy,\quad n=1,2;~i,j=1,2,3,4$$
The local stiffness matrices contributed by the fracture $\bm{l}$ in \textit{Case (a)} and \textit{Case (b)} are shown as follows,
\[\bm{l}\left(\textit{Case a}\right):\begin{bmatrix}
    0 & 0 & 0 & 0 \\
    0 & a^{l}_{22} & a^{l}_{23} & 0 \\
    0 & a^{l}_{32} & a^{l}_{33} & 0 \\
    0 & 0 & 0 & 0
\end{bmatrix},\qquad
\bm{l}\left(\textit{Case b}\right):\begin{bmatrix}
    a^{l}_{11} & a^{l}_{12} & a^{l}_{13} & a^{l}_{14} \\
    a^{l}_{21} & a^{l}_{22} & a^{l}_{23} & a^{l}_{24} \\
    a^{l}_{31} & a^{l}_{32} & a^{l}_{33} & a^{l}_{34} \\
    a^{l}_{41} & a^{l}_{42} & a^{l}_{43} & a^{l}_{44}
\end{bmatrix}\]
where
$$a^{l}_{i,j}=\int_{\bm{l}}\epsilon k_f\frac{\partial\Psi_i}{\partial \nu}\frac{\partial\Psi_j}{\partial \nu}ds,\quad i,j=1,2,3,4.$$

Note that in \textit{Case (a)}, where the grid is conforming with fracture, the only non-zero entries in local stiffness matrix of fracture are $a^{l}_{22},a^{l}_{23},a^{l}_{32},a^{l}_{33}$, because $\Psi_1=\Psi_4=0$ along $\bm{l}$. On the other hand, \textit{Case (b)} shows the complete form of local stiffness matrix resulted from the fracture $\bm{l}$. Finally, the matrix $\bm{A}$ in this local zone is assembled in the following manners:

$\textit{Case (a)}$:
\begin{equation}\label{ConformingA}
\begin{bmatrix}
    a^{T_1}_{11} & a^{T_1}_{12} & a^{T_1}_{13} & 0 \\
    a^{T_1}_{21} & a^{T_1}_{22}+a^{T_2}_{22}+a^{l}_{22} & a^{T_1}_{23}+a^{T_2}_{23}+a^{l}_{23} & a^{T_2}_{24} \\
    a^{T_1}_{31} & a^{T_1}_{32}+a^{T_2}_{32}+a^{l}_{32} & a^{T_1}_{33}+a^{T_2}_{33}+a^{l}_{33} & a^{T_2}_{34} \\
    0 & a^{T_2}_{42} & a^{T_2}_{43} & a^{T_2}_{44}
\end{bmatrix}
\end{equation}

$\textit{Case (b)}$:
\begin{equation}\label{NonconformingA}
\begin{bmatrix}
    a^{T_1}_{11}+a^{l}_{11} & a^{T_1}_{12}+a^{l}_{12} & a^{T_1}_{13}+a^{l}_{13} & a^{l}_{14} \\
    a^{T_1}_{21}+a^{l}_{21} & a^{T_1}_{22}+a^{T_2}_{22}+a^{l}_{22} & a^{T_1}_{23}+a^{T_2}_{23}+a^{l}_{23} & a^{T_2}_{24}+a^{l}_{24} \\
    a^{T_1}_{31}+a^{l}_{31} & a^{T_1}_{32}+a^{T_2}_{32}+a^{l}_{32} & a^{T_1}_{33}+a^{T_2}_{33}+a^{l}_{33} & a^{T_2}_{34}+a^{l}_{34} \\
    a^{l}_{41} & a^{T_2}_{42}+a^{l}_{42} & a^{T_2}_{43}+a^{l}_{43} & a^{T_2}_{44}+a^{l}_{44}
\end{bmatrix}
\end{equation}
We obtain the global stiffness matrix $\bm{A}$ after the loop of all triangles and fractures.
\begin{remark}
Just like the conforming mesh can be viewed as a special case of non-conforming mesh, the stiffness matrix \eqref{ConformingA} is also a special case of the general form \eqref{NonconformingA} when cell interfaces matching the fractures. Moreover, we note that \eqref{ConformingA} is exactly the same local stiffness matrix built in the traditional Galerkin finite element DFM \cite{FEMDFM}, in which the aforementioned coupling procedure for stiffness matrix was explained as superposition principle based on conforming mesh. Therefore we can conclude that, theoretically, when the fracture is conforming with gridcells, our scheme degenerates to classical finite element DFM. The numerical verification of this conclusion will be demonstrated in Section 5.
\end{remark}
\begin{remark}
We demonstrate how to implement the quadrature for the integration on fractures, especially for curved fractures.
We describe the fracture $\displaystyle\bm{l}$ as a parametric equation $\displaystyle x=x(t), y=y(t), \alpha\leq t\leq \beta$, where the equation is non-degenerate, i.e. $\displaystyle\dot{x}^2+\dot{y}^2>0, \forall t\in[\alpha,\beta]$. The line integral $\displaystyle\int_{\bm l}\epsilon k_f\frac{\partial \Psi_i}{\partial\nu}\frac{\partial \Psi_j}{\partial\nu} ds$ can be written as a definite integral $\displaystyle \sum_{K\in \mathscr{P}}\int^{t^{K}_{b}}_{t^{K}_{a}} \epsilon k_f \frac{[\nabla \Psi_i\cdot (\dot{x},\dot{y})][\nabla \Psi_j\cdot (\dot{x},\dot{y})]}{\sqrt{\dot{x}^2+\dot{y}^2}}dt$, where $\mathscr{P}$ is the collection of elements on the path of the fracture and $t^{K}_{a}, t^{K}_{b}$ are the parameters at which the fracture intersects with the element $K$.
The intersections between the fractures and the elements have to be found to compute the parameters $t^{K}_{a}, t^{K}_{b}$ on each element. To evaluate this definite integration in practice, we can choose Gauss or Gauss-Lobatto quadrature rule to generate the quadrature points between $t^{K}_{a}, t^{K}_{b}$. Note that if the fracture is a straight line and $\epsilon, k_f$ are constant, the midpoint rule is enough.
\end{remark}

\section{Numerical tests}
In this section, we provide eight numerical tests, roughly in an increasing order of the geometrical
complexity, to show the performance of the NDFM.
Part of the numerical examples are either chosen from common fracture settings \cite{pEDFM,CrossShaped} or well known benchmarks \cite{Benchmark,benchmark2,benchmark3,ThesisXFEM,Hydrocoin,LMFEM2D} so that one can refer to these articles for more details about reliable reference solutions.
Some others are given to shown either the consistency with traditional finite element DFM on conforming meshes, or the rate of convergence by comparing with an analytical solution. The last two experiments are for curved fractures and 3D cases. All source codes are available at \url{https://github.com/ziyaoxu/Nonconforming-DFM.git}.
Special thanks go to the authors of \cite{Benchmark} for sharing the data of grids and results \cite{Website} computed by a number of DFM algorithms. These works enable us to compare and evaluate our model with the existing ones. We declare that all the data and figures used in Example \ref{ex4}, \ref{ex5} and \ref{ex6} for the reference, comparison, and evaluation of our numerical results come from them.

\begin{figure}[htbp]
\centering
\includegraphics[width=3in]{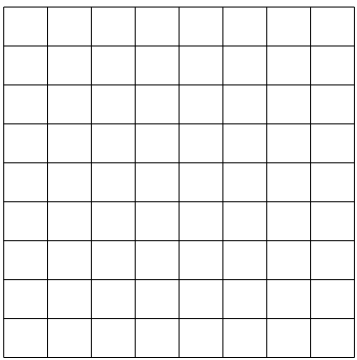}
\caption{A uniform rectangular mesh ($N_x=8, N_y=9$) } \label{fig:UniformRect}
\end{figure}

For the sake of simplicity, the flow in these tests is driven by boundary conditions instead of source term, i.e. $f=0$. In all the examples, we use linear Lagrange shape functions to construct the finite element space. Both unstructured triangular meshes and uniform rectangular meshes are employed to discretize the computational domain $\Omega_h$. A typical uniform rectangular mesh is shown in Figure \ref{fig:UniformRect}. The uniform rectangular meshes won't be exhibited later since they are all similar with this one, but the unstructured triangular discretizations will be shown in each individual examples if necessary.

\begin{figure}[!htbp]
\subfigure[Domain and boundary conditions]{\includegraphics[width = 3in]{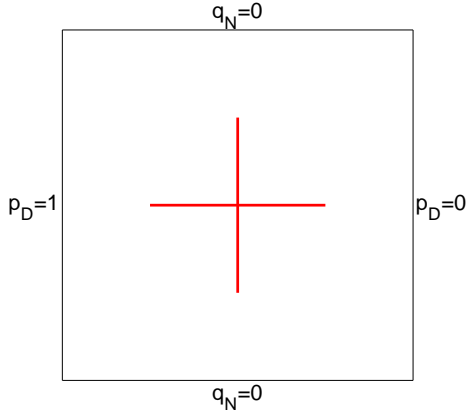}} \subfigure[Reference solution on fully resolved mesh]{\includegraphics[width = 3in]{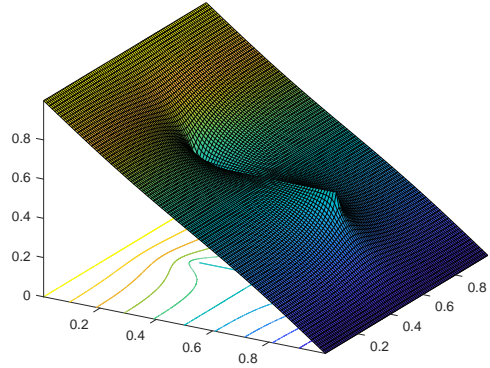}} \caption{Fracture setting and reference solution of Example \ref{ex1}}\label{fig:DomainandReference}
\end{figure}

\begin{ex} \label{ex1}
\textbf{Cross-shaped fractures}

In this example, we test a fractured media with simple cross-shaped fractures and compare the numerical results on different uniform rectangular meshes with a reference solution on fully resolved mesh.
The computational domain is $[0,1]\times[0,1]$ with cross-shaped fractures $[0.25,0.75]\times[0.4995,0.5005]\bigcup[0.4995,0.5005]\times[0.25,0.75]$ laying on it. The permeability of porous matrix and fractures regions are $1$ and $10^{8}$, respectively. Moreover, the Dirichlet boundary conditions $p_D=1$ and $p_D=0$ are imposed on the left and right boundaries respectively, and the top and bottom boundaries are set to be impervious, i.e. $q_N=0$. See Figure \ref{fig:DomainandReference}(a) for an illustration of the domain and boundary conditions. This test case is the same as the one given in \cite{pEDFM} so one can refer to their result for more detailed comparison.
\end{ex}

\begin{figure}[!htbp]
\subfigure[Solution on $10\times10$ mesh]{\includegraphics[width = 3in]{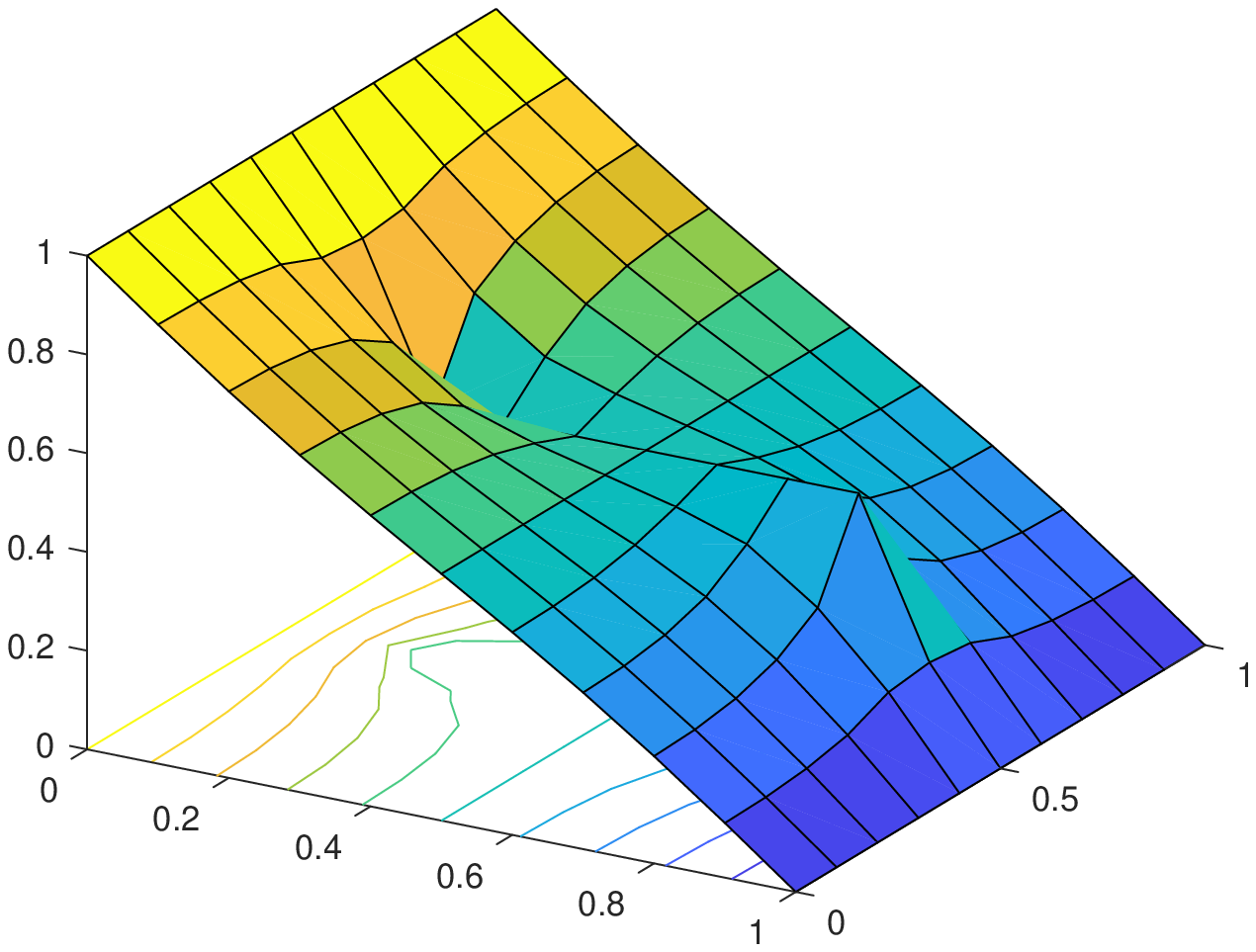}} \subfigure[Solution on $11\times11$ mesh]{\includegraphics[width = 3in]{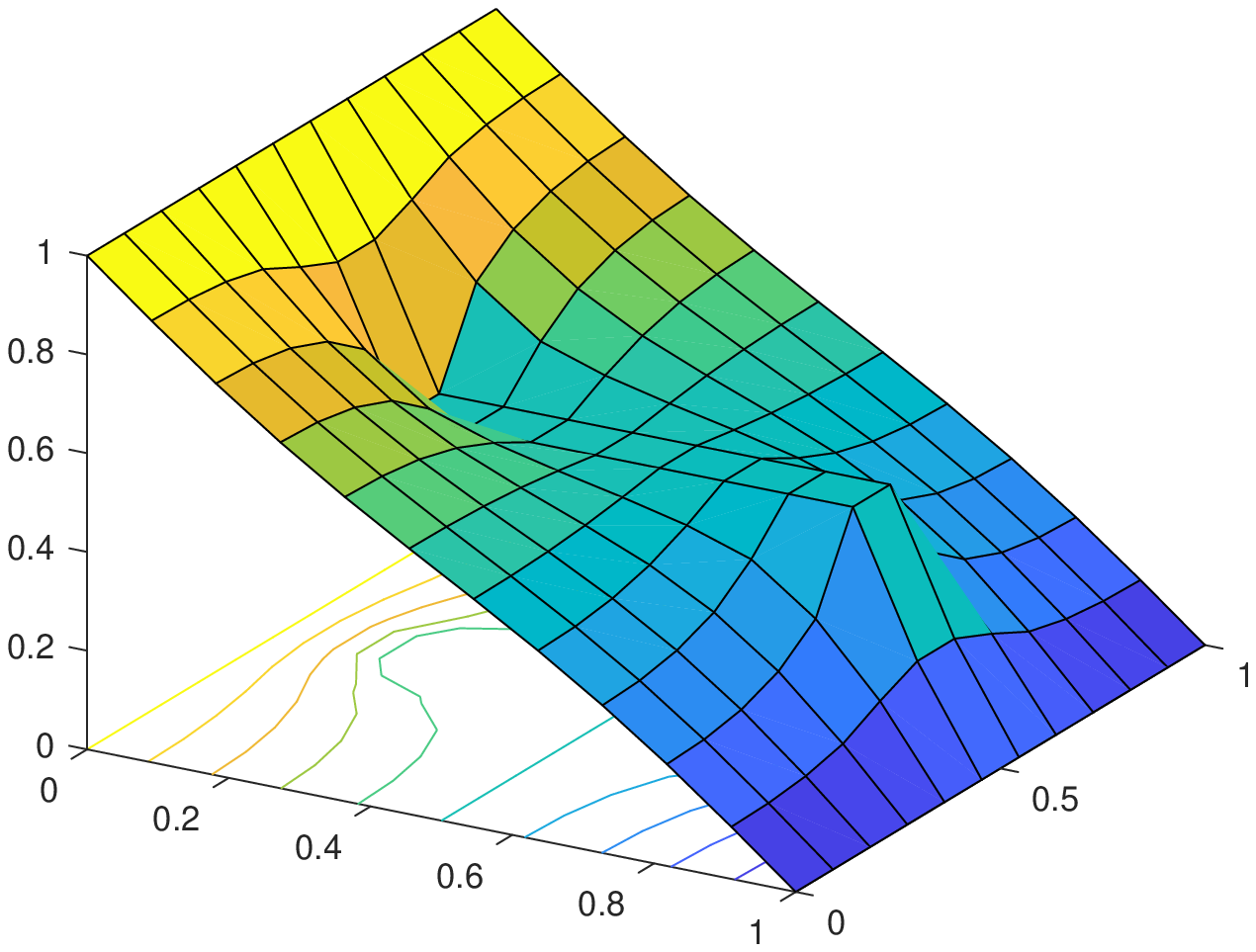}}\\
\subfigure[Solution on $20\times20$ mesh]{\includegraphics[width = 3in]{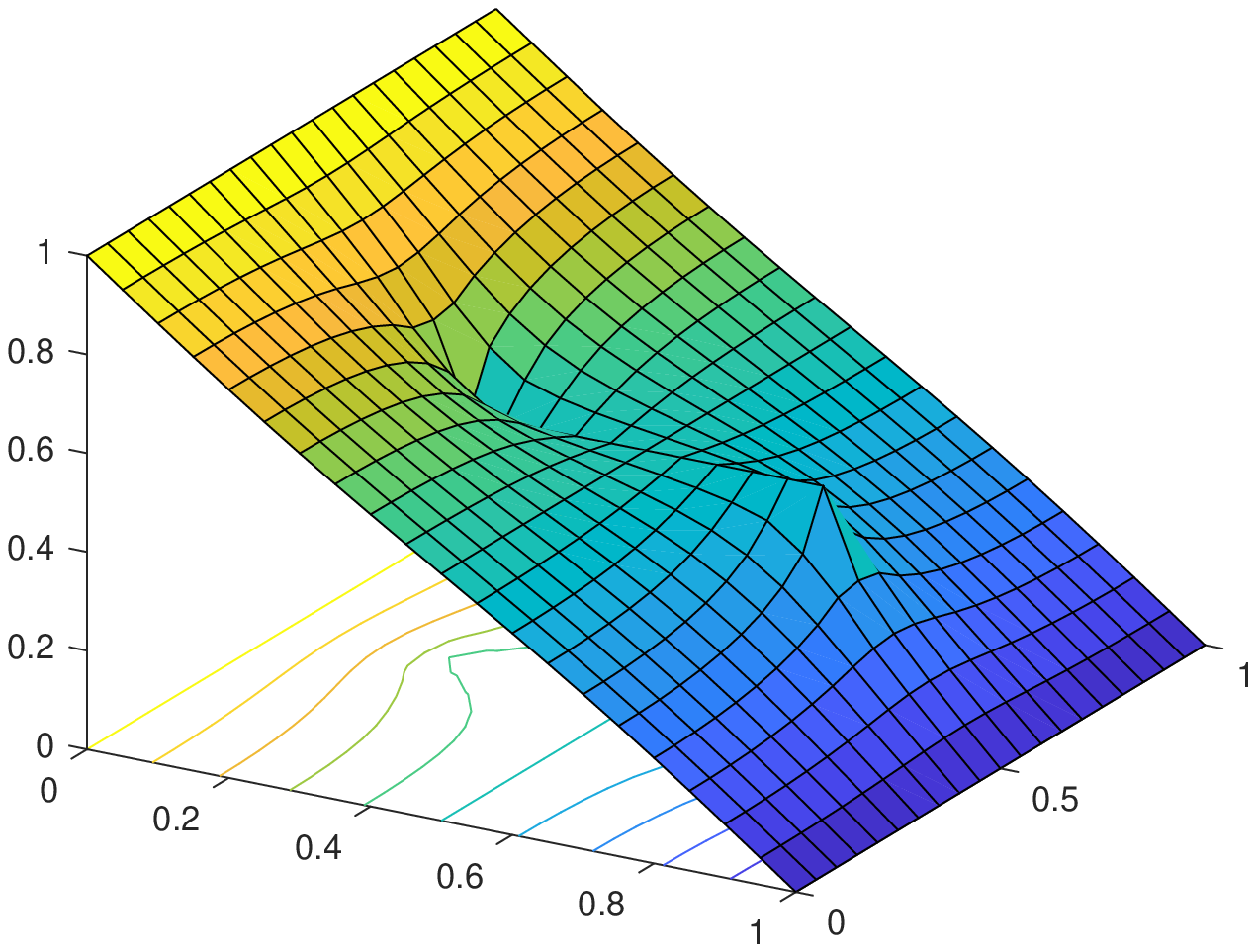}} \subfigure[Solution on $21\times21$ mesh]{\includegraphics[width = 3in]{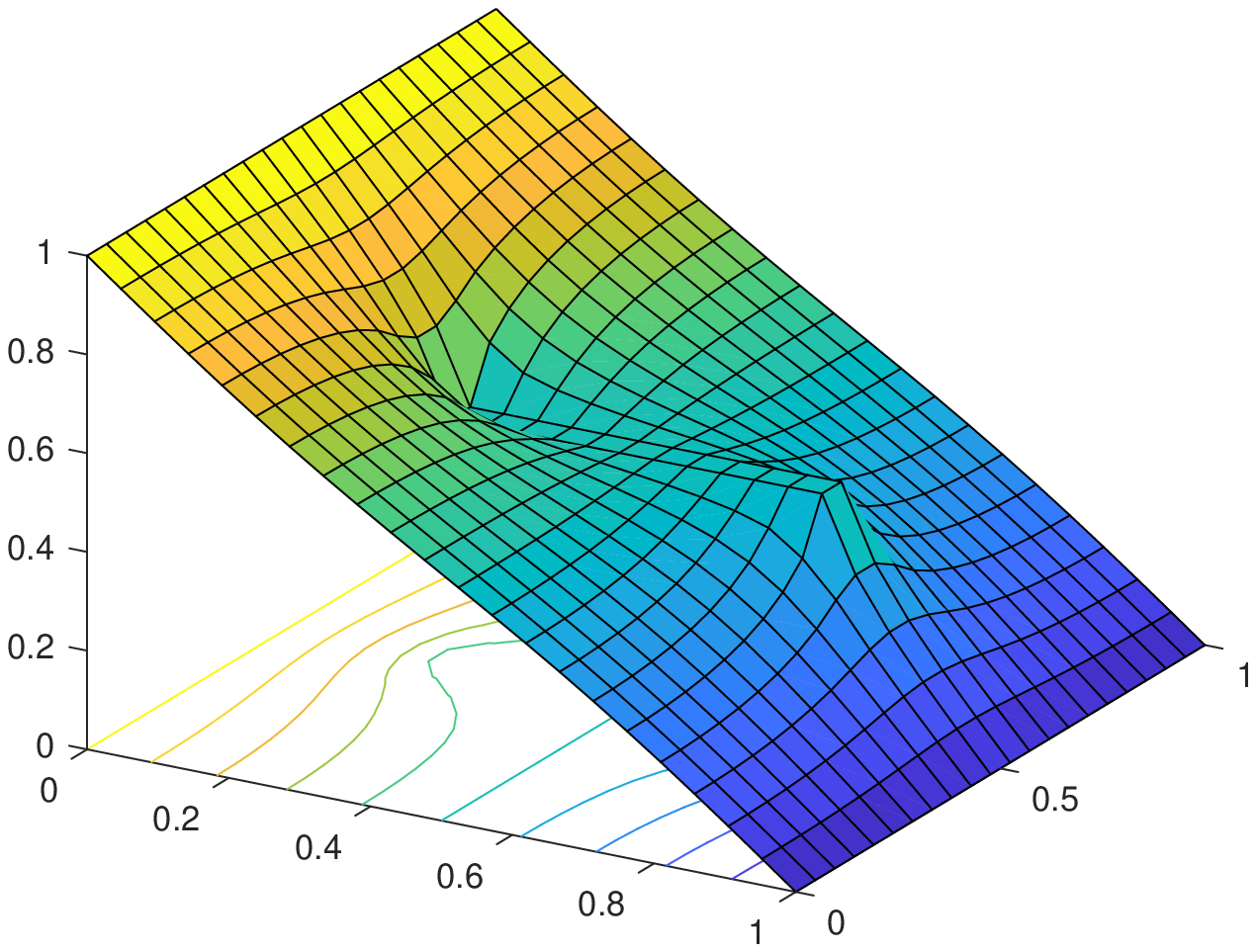}}\\
\subfigure[Solution on $50\times50$ mesh]{\includegraphics[width = 3in]{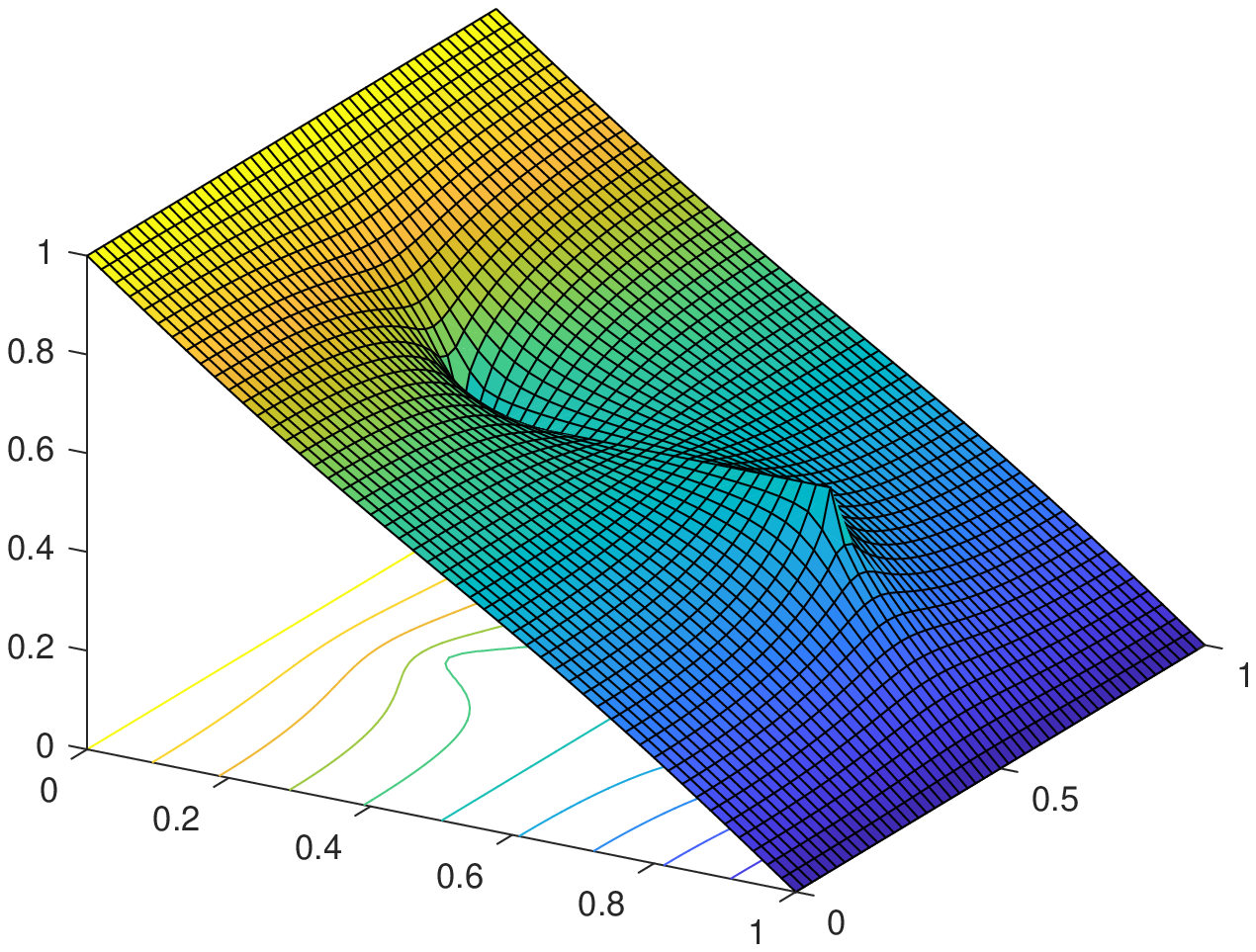}} \subfigure[Solution on $51\times51$ mesh]{\includegraphics[width = 3in]{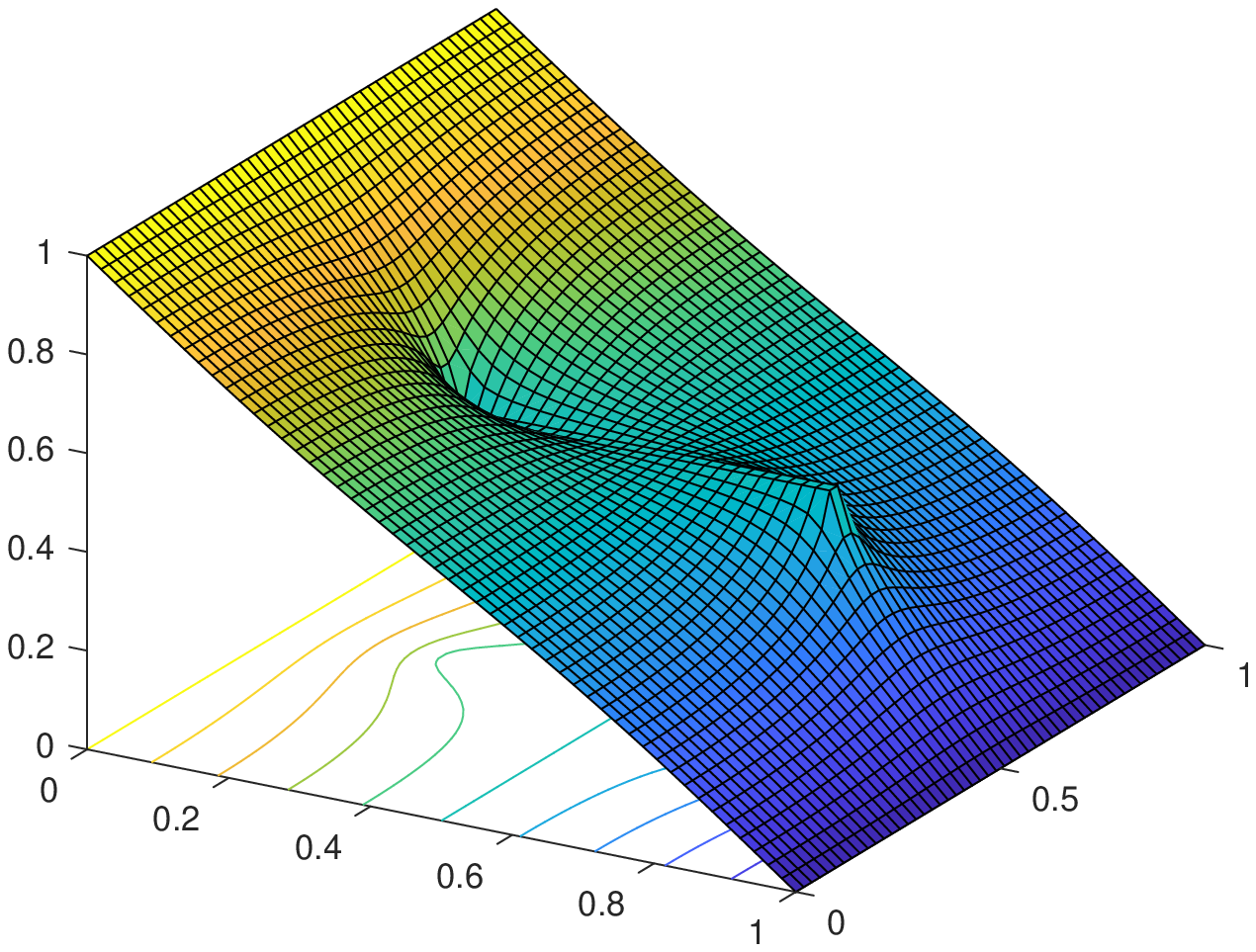}}
\caption{Solutions on different meshes of Example \ref{ex1}}\label{fig:ex1Solutions}
\end{figure}

\begin{figure}[!htbp]
\subfigure[Pressure along $y=0.45$]{\includegraphics[width = 3in]{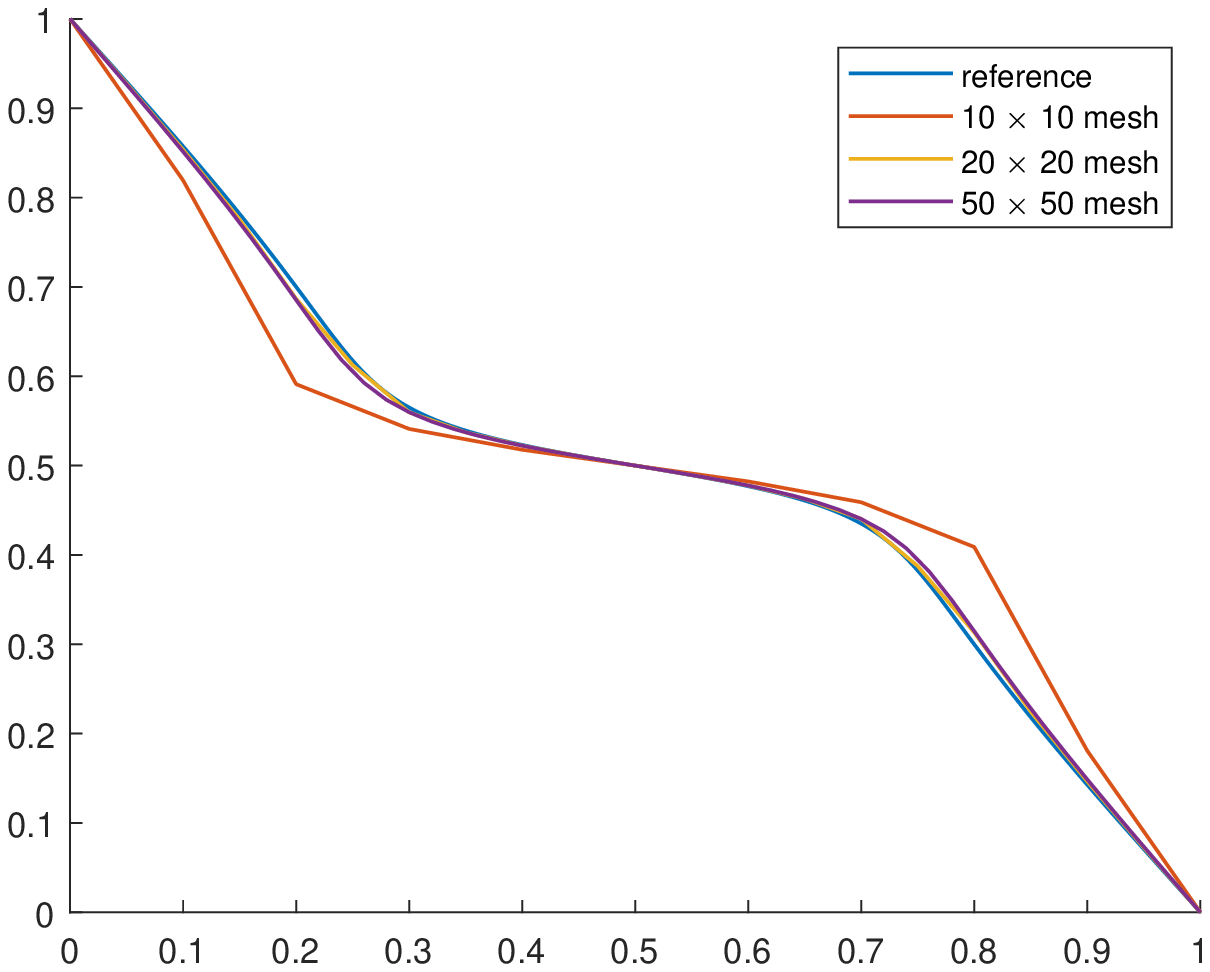}}
\subfigure[Pressure along $y=0.45$]{\includegraphics[width = 3in]{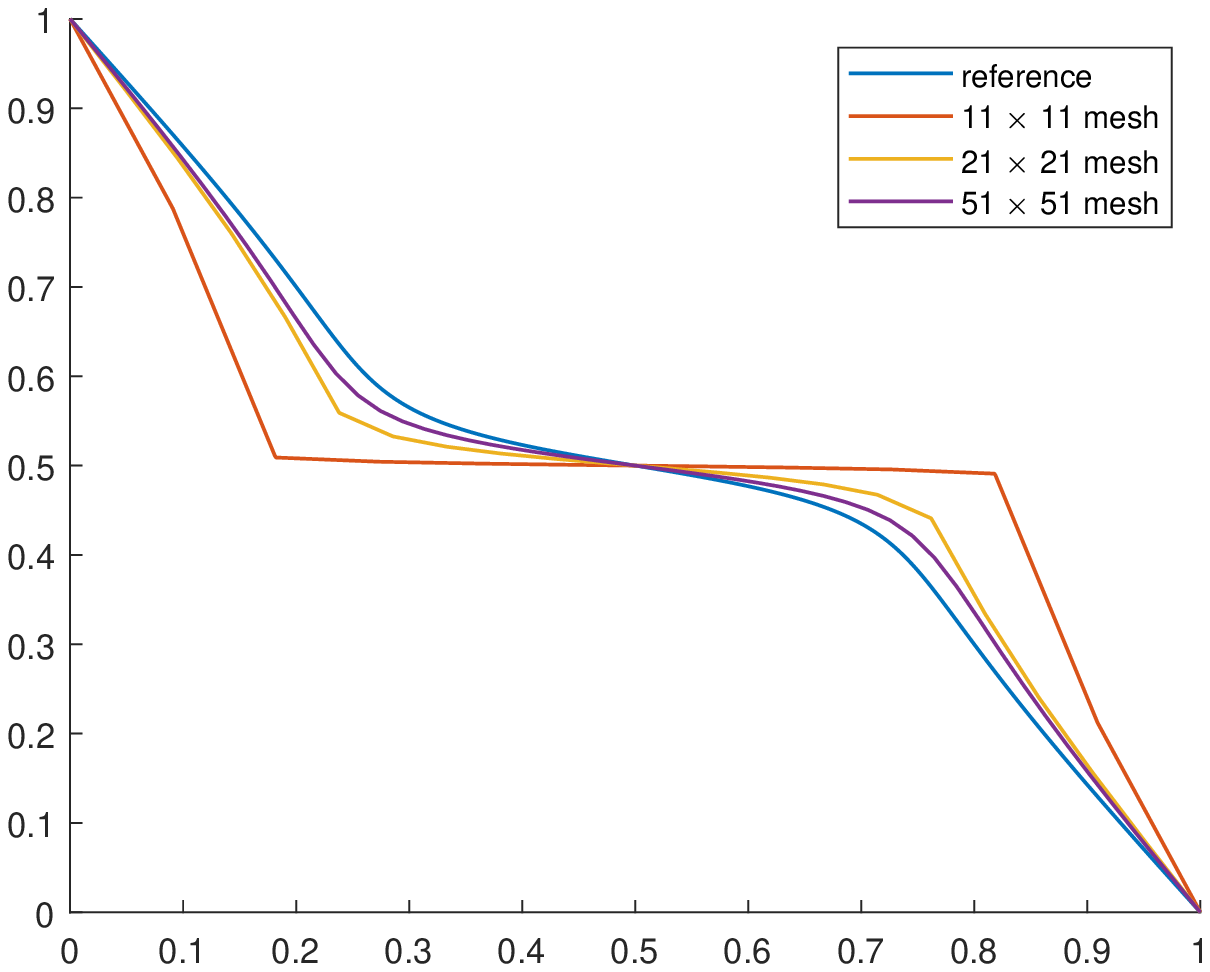}}\\
\subfigure[Pressure along $x=0.3$]{\includegraphics[width = 3in]{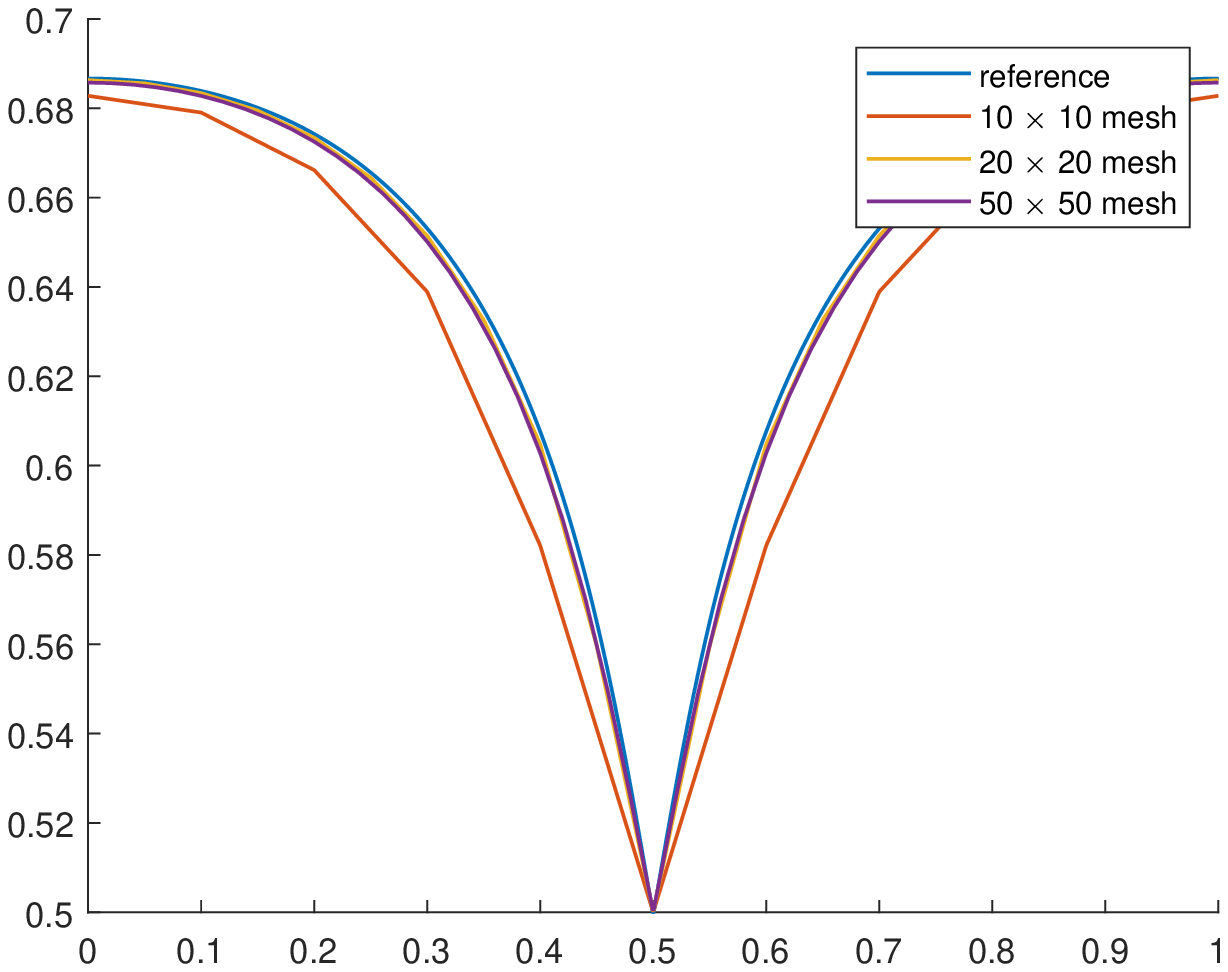}} \subfigure[Pressure along $x=0.3$]{\includegraphics[width = 3in]{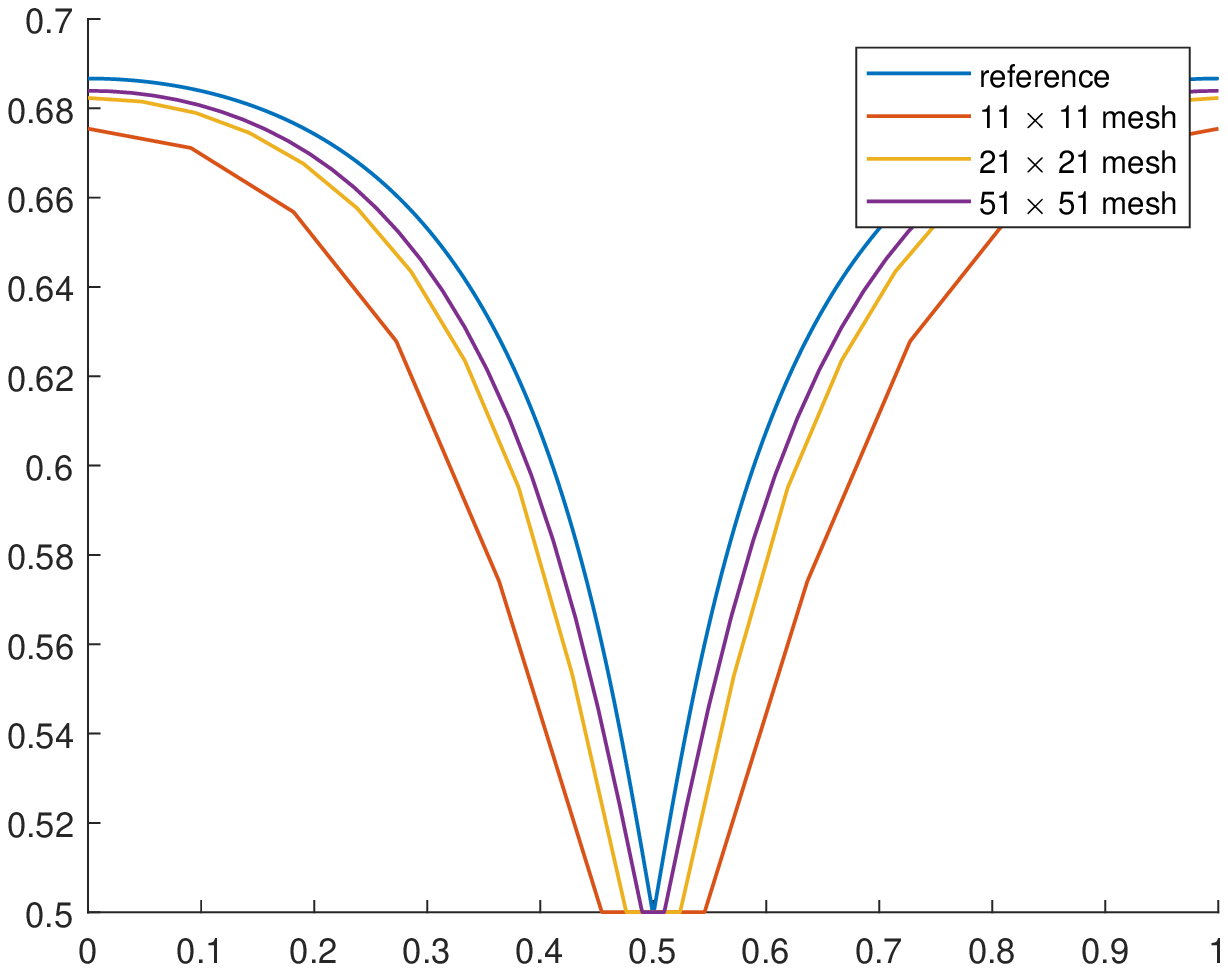}}\\
\subfigure[Pressure along $x=0.4$]{\includegraphics[width = 3in]{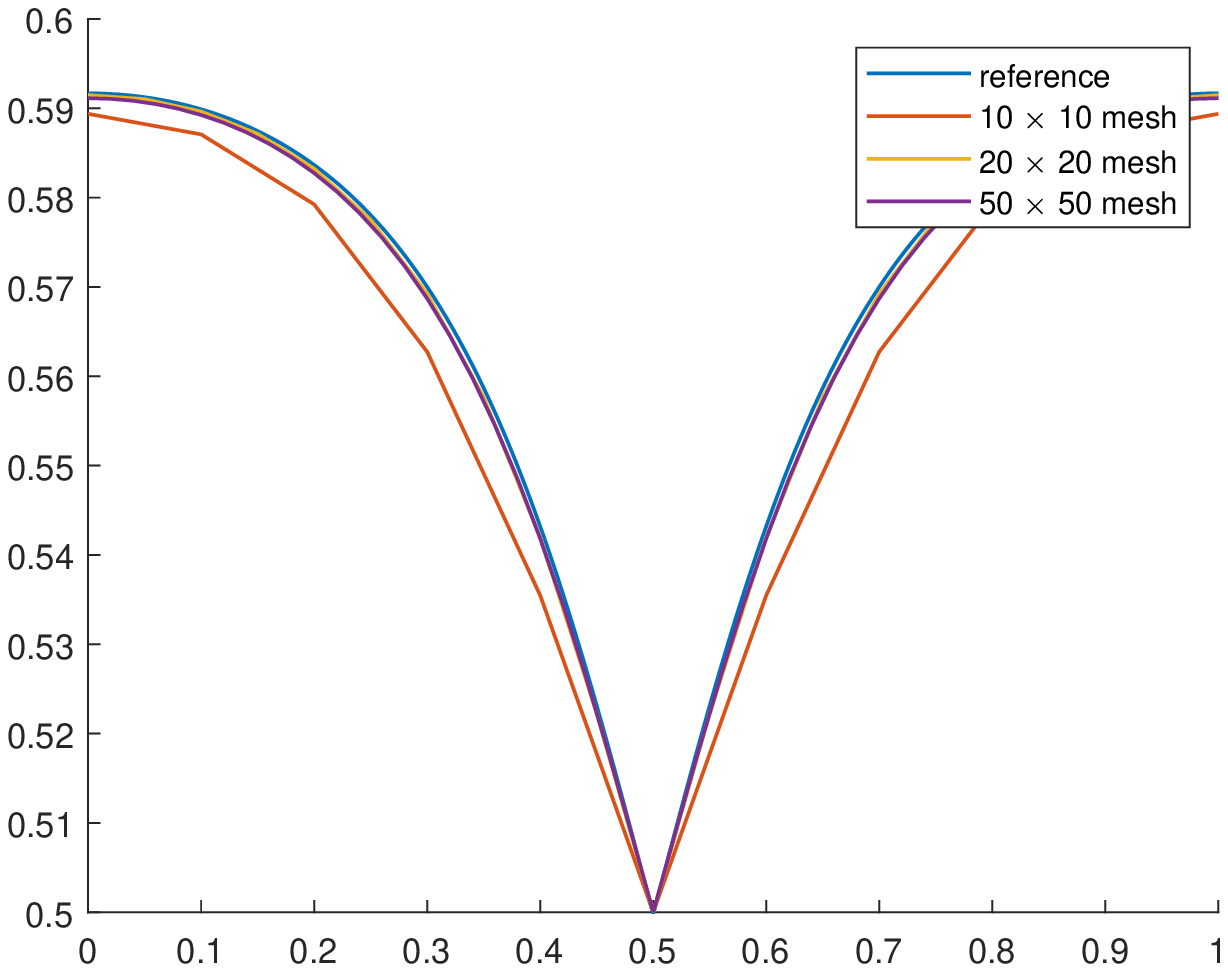}} \subfigure[Pressure along $x=0.4$]{\includegraphics[width = 3in]{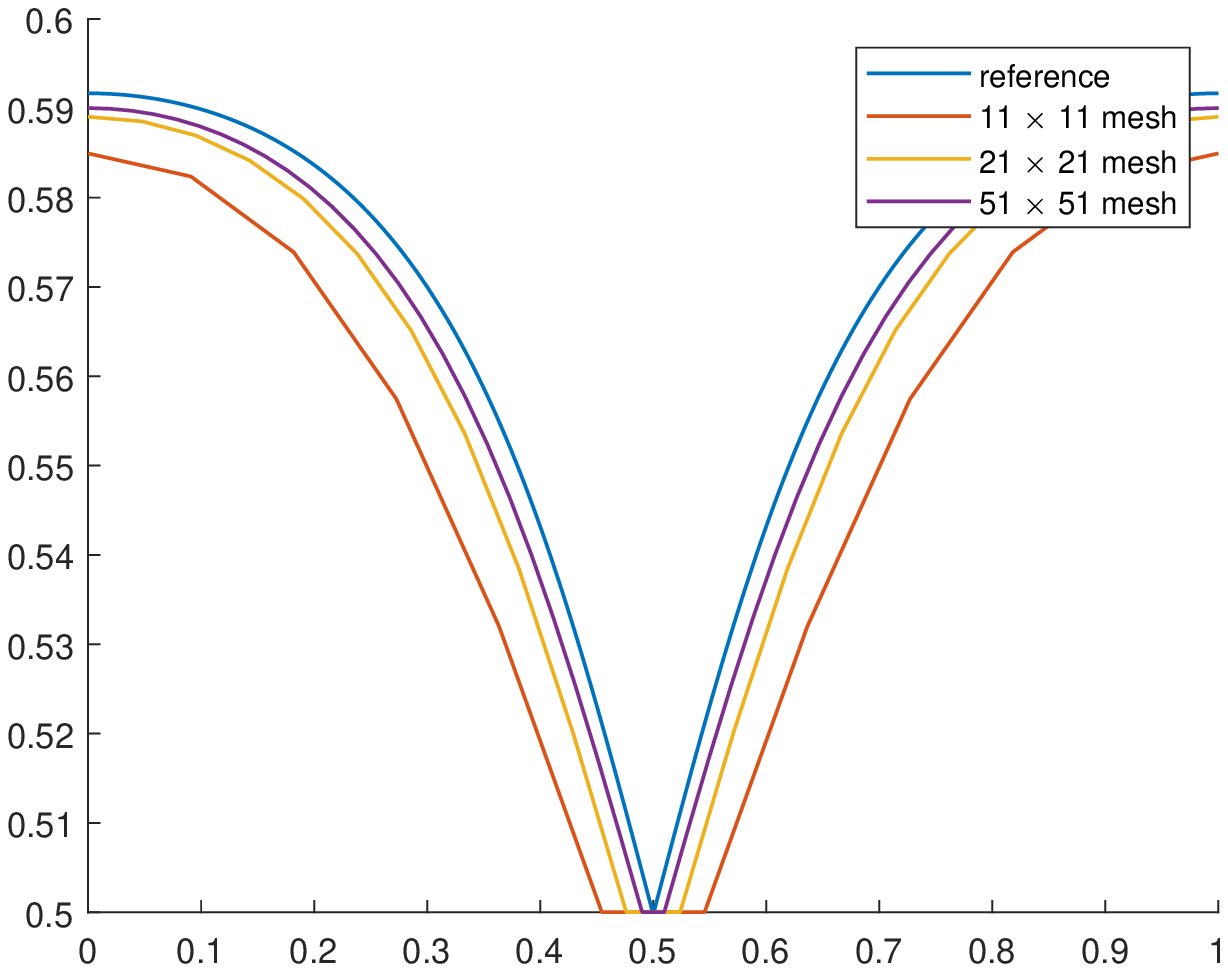}}
\caption{Pressure profiles along different lines of Example \ref{ex1}}\label{fig:ex1Slices}
\end{figure}

We employ the standard Galerkin finite element method on a fully resolved rectangular mesh ($N_x=N_y=1001$) to give the reference solution, with the fractured media treated as a heterogeneous media, i.e. the equi-dimensional model presented in Section 2. The surface and contour of the reference pressure are given in the Figure \ref{fig:DomainandReference}(b), and the results of the non-conforming DFM algorithm on different rectangular meshes are shown in Figure \ref{fig:ex1Solutions}.
Moreover, a comparison of pressure profiles sliced along $y=0.5, x=0.3$ and $x=0.4$ are shown in Figure \ref{fig:ex1Slices}. By comparing these results with the reference solution, we can see that they match well, especially when the number of grids increase.

One may also notice that when the number of grids are odd (the cases shown in the right column), the pressure is flat along some gridcells passed through by fractures, which seems not as good as that with even number grids (the cases shown in the left column). This phenomenon is reasonable since when fractures go through the center of gridcells under such a symmetric domain and boundary condition setting, the flat pressure is the best that an linear approximation can do.

\begin{figure}[!htbp]
\subfigure[Conforming triangulation\cite{DistMesh} and fracture networks]{\includegraphics[width = 3in]{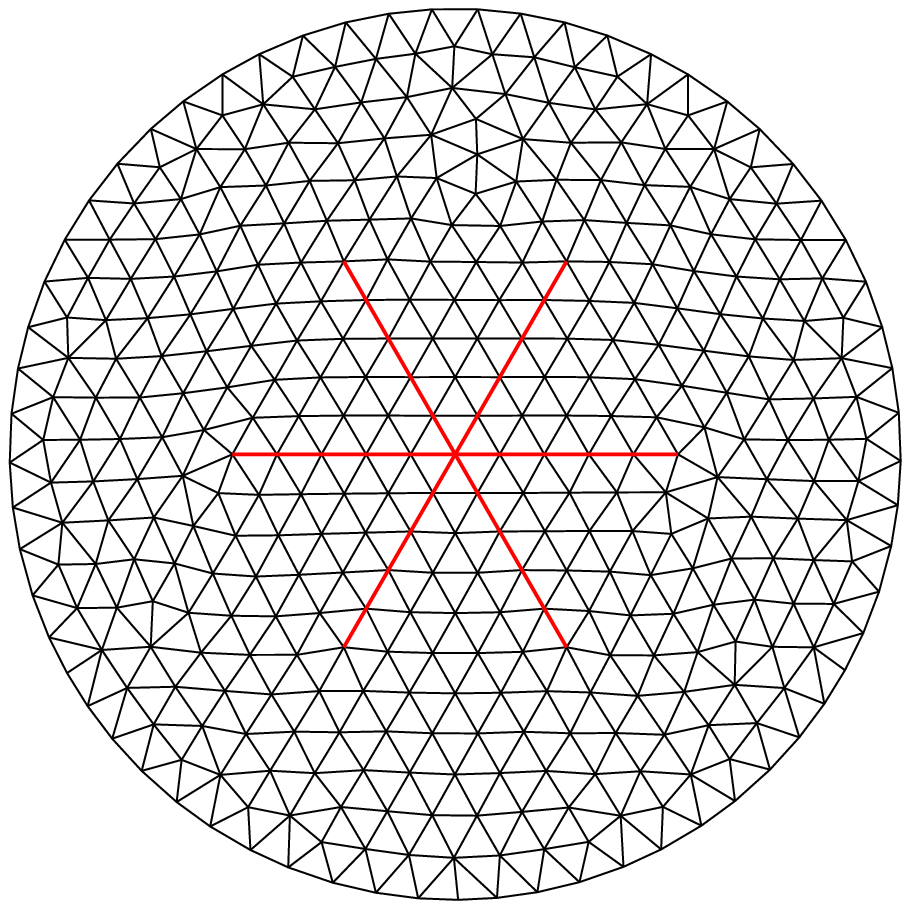}} \subfigure[Numerical result of DFM]{\includegraphics[width = 3in]{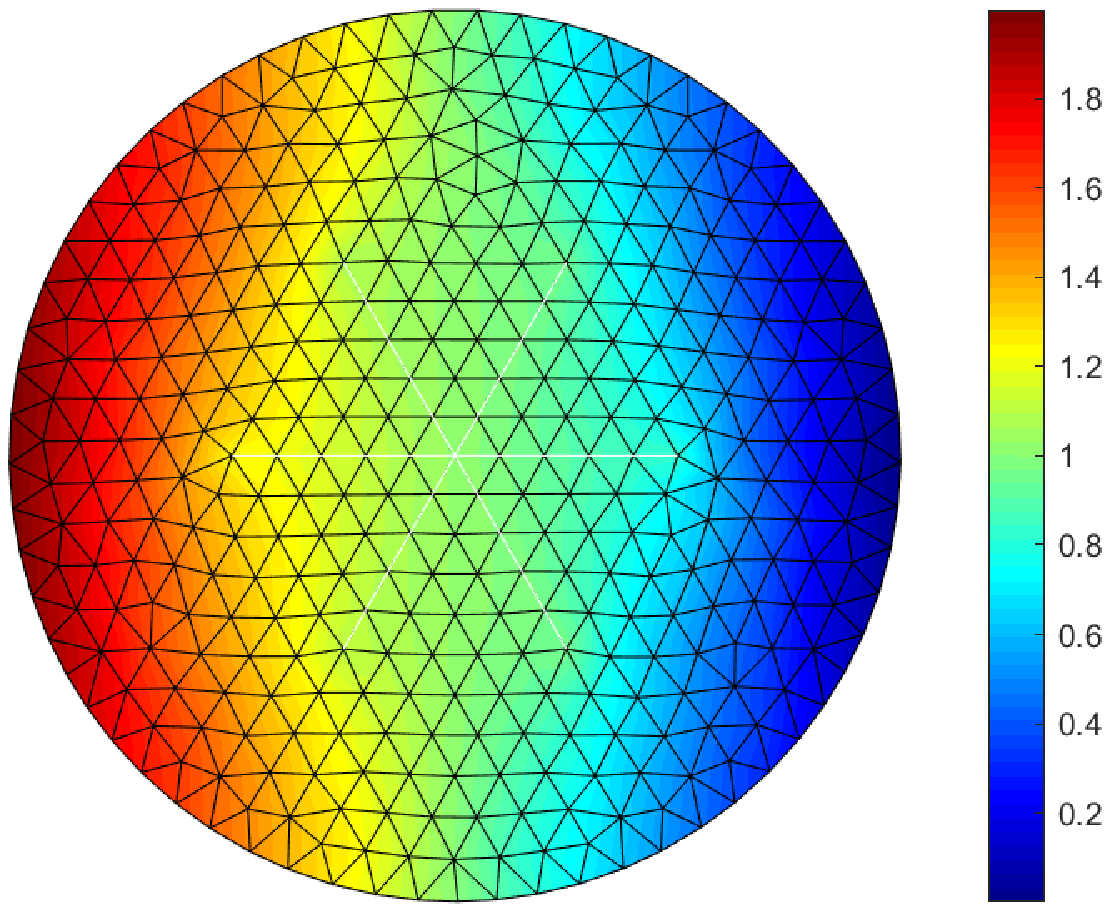}} \caption{Triangulation, fracture networks and numerical result of DFM in \ref{ex2}}\label{fig:ex2DFM}
\end{figure}

\begin{ex} \label{ex2}
\textbf{Consistency test}

In this example, we demonstrate the consistency between the NDFM and traditional DFM under conforming meshes. The most straightforward way to shown this property is to construct a sequence of non-conforming meshes that converges to a conforming mesh. Then we can numerically prove the consistency if the the numerical results of NDFM converge to that of DFM during this process. However, to avoid the interpolation used in the comparison of solutions among different meshes, we choose the following equivalent way. we first construct a fractured porous media with a conforming triangulation on it. Then, instead of changing the mesh, we keep the same triangulation but set a sequence of different fracture networks which is non-conforming on this mesh to converge to the aforementioned fracture network. Finally, we compute the numerical solutions on those fractured media. If the solutions of non-conforming fracture networks converge to the solution of the conforming one, the consistency is verified.

We set the porous matrix to be the unite circle $\Omega=\{(x,y)|x^2+y^2\leq 1\}$ with permeability $1$. The conforming fracture network consist of three fractures centered at the origin with length $1$, thickness $10^{-4}$, angles $\theta=0, \pi/3, -\pi/3$ and permeability $k_f=10^{4},2\times10^{4},3\times10^{4}$, respectively. See Figure \ref{fig:ex2DFM}(a) for an illustration. The sequence of non-conforming fracture networks are generated by adding smaller and smaller rotation $\Delta\theta$ and shift $(\Delta x,\Delta y)$ on the original fracture network. The boundary conditions of the fractured media is Dirichlet boundary with $p_D=1-x$.

\end{ex}

\begin{figure}[!htbp]
\subfigure[Non-conforming triangulation and fracture network 1]{\includegraphics[width = 3in]{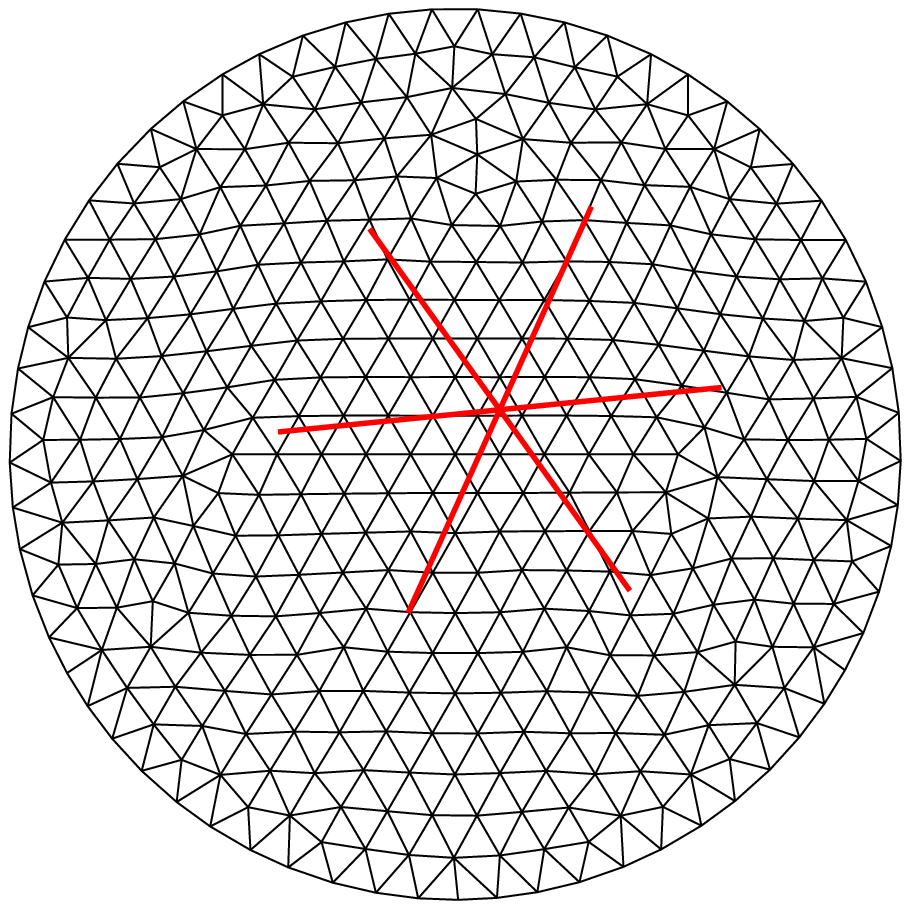}}
\subfigure[Numerical result of corresponding NDFM]{\includegraphics[width = 3in]{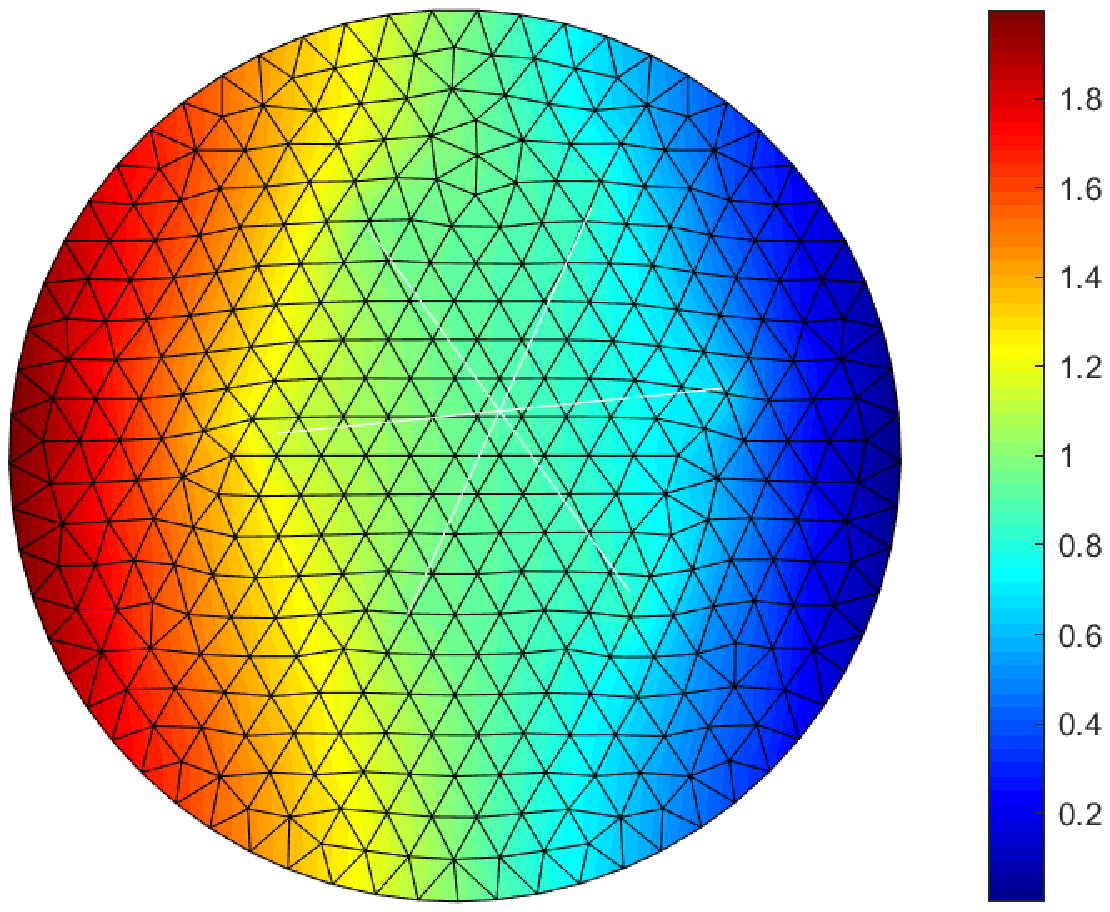}} \\
\subfigure[Non-conforming triangulation and fracture network 2]{\includegraphics[width = 3in]{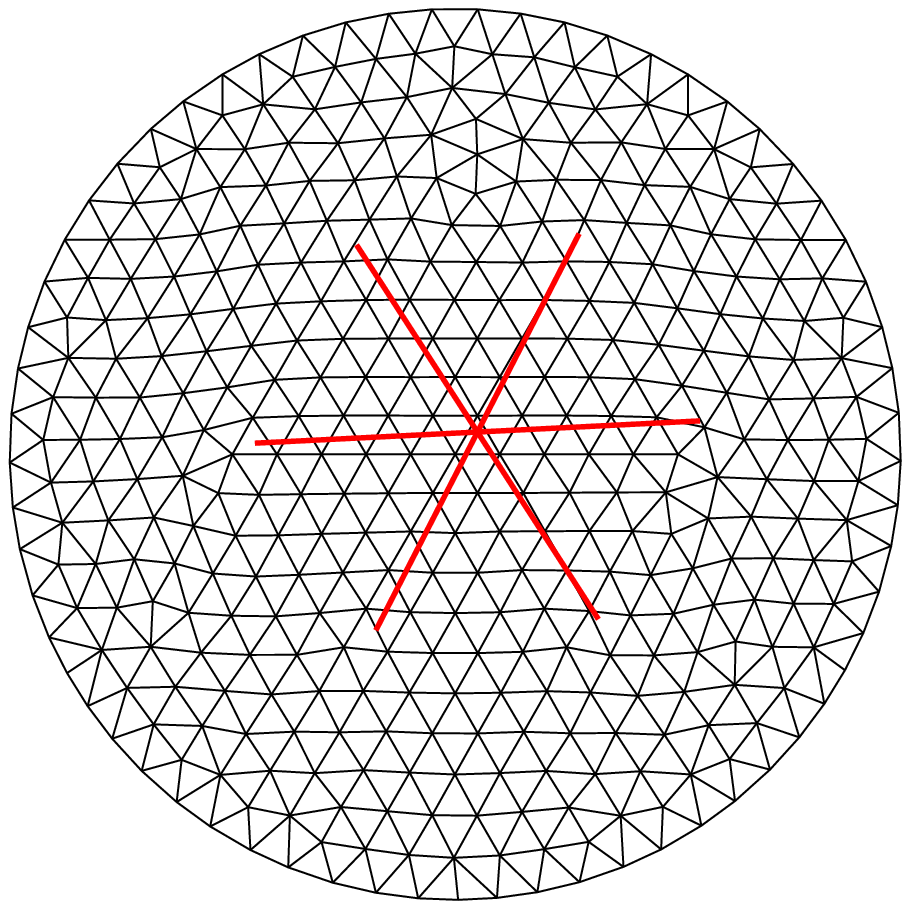}}
\subfigure[Numerical result of corresponding NDFM]{\includegraphics[width = 3in]{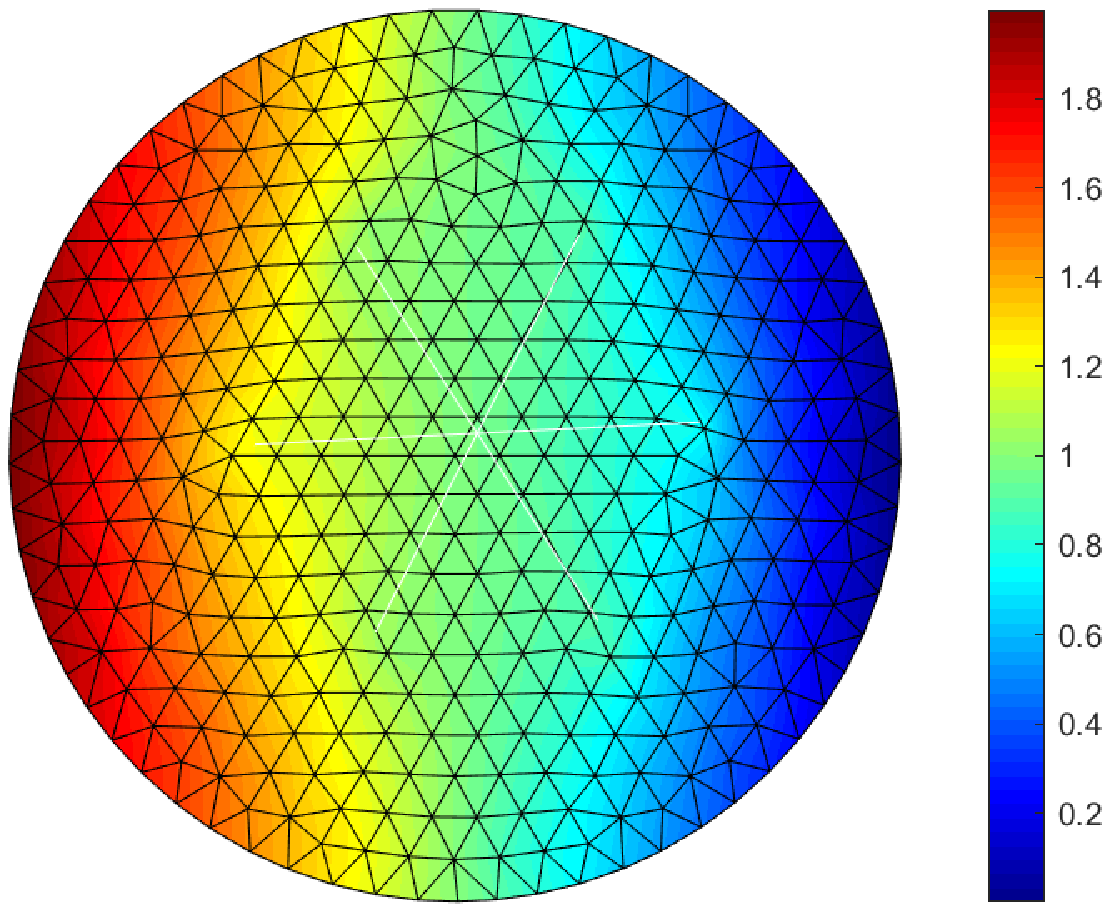}} \\
\subfigure[Non-conforming triangulation and fracture network 3]{\includegraphics[width = 3in]{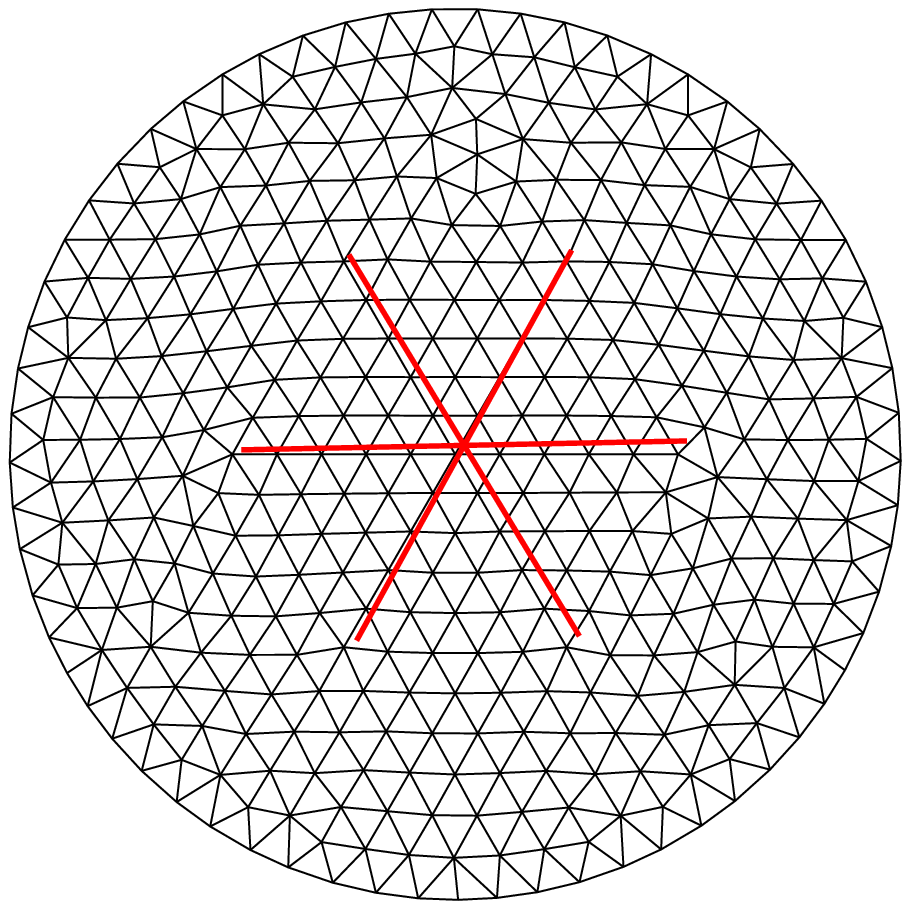}}
\subfigure[Numerical result of corresponding NDFM]{\includegraphics[width = 3in]{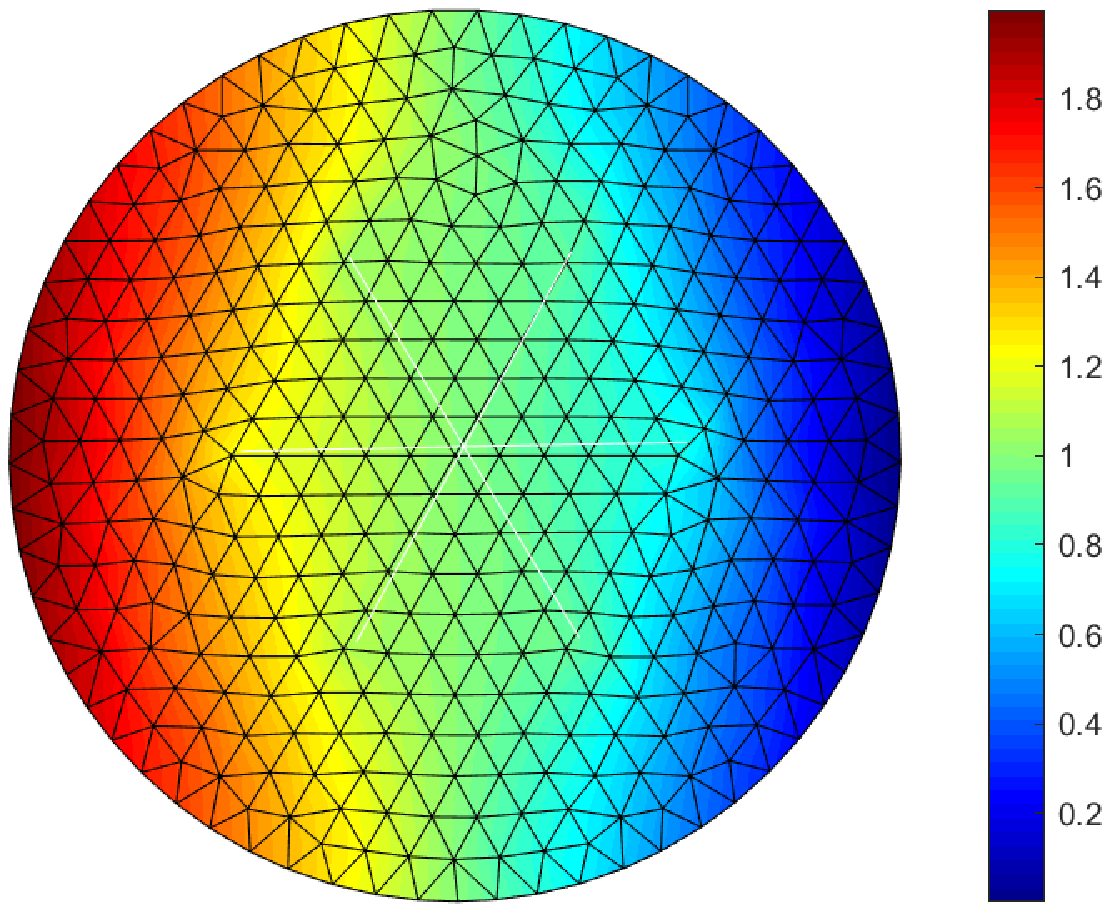}} \\
\caption{Triangulation, fracture networks and corresponding numerical results of NDFM in \ref{ex2}}\label{fig:ex2NDFM}
\end{figure}

The numerical solution of the fractured media in which the fractures and triangulation are conforming is shown in Figure \ref{fig:ex2DFM}(b). We only draw the the solutions of first three non-conforming fracture networks in Figure \ref{fig:ex2NDFM} because the others are not visually distinguishable with the conforming one. The maximum norm of the differences between the solutions of non-conforming and conforming networks as well as the scales of their shift and rotation are summarized in Table \ref{tab:ex2differences}. From the table we can see the solutions of NDFM converge to that of DFM as the fracture networks converge to the conforming one, which proves the consistency of NDFM.

\begin{table}[!htb]
\centering
\caption{\label{tab:ex2differences} Differences between the solutions of NDFM and DFM in Example \ref{ex2}}
\begin{tabular}{|c|c c c|c |}
  \hline
Number &$\Delta x$ & $\Delta y$ & $\Delta \theta$ & $||p_{\text{DFM}}-p_\text{{NDFM}}||_{\infty}$ \\\hline
1&1E-01&1E-01&1E-01&1.46E-01 \\\hline
2&5E-02&5E-02&5E-02&1.03E-01  \\\hline
3&2E-02&2E-02&2E-02&8.13E-02\\\hline
4&1E-02&1E-02&1E-02&5.53E-02 \\ \hline
5&1E-03&1E-03&1E-03&8.73E-03 \\ \hline
6&1E-04&1E-04&1E-04&9.79E-04 \\ \hline
7&1E-06&1E-06&1E-06& 9.92E-06 \\ \hline
8&1E-08&1E-08&1E-08& 9.92E-08 \\ \hline
\end{tabular}
\end{table}

\begin{figure}[!htbp]
\subfigure[$\theta=0$]{\includegraphics[width = 3in]{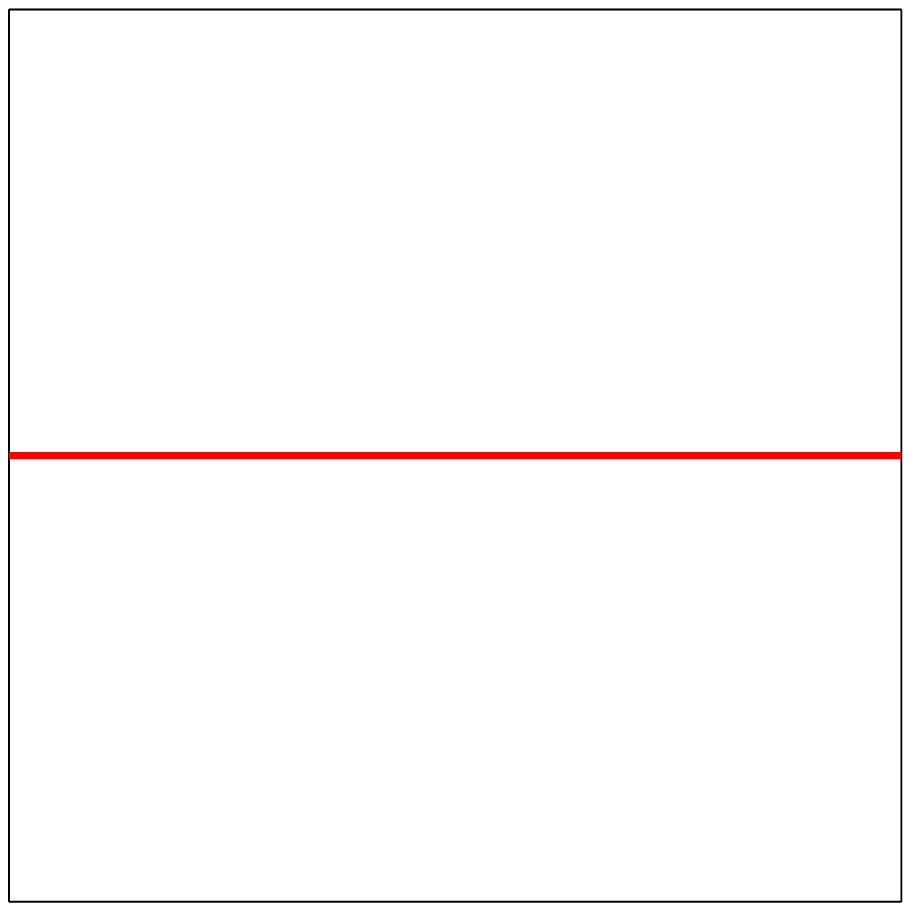}} \subfigure[$\theta=5.3$]{\includegraphics[width = 3in]{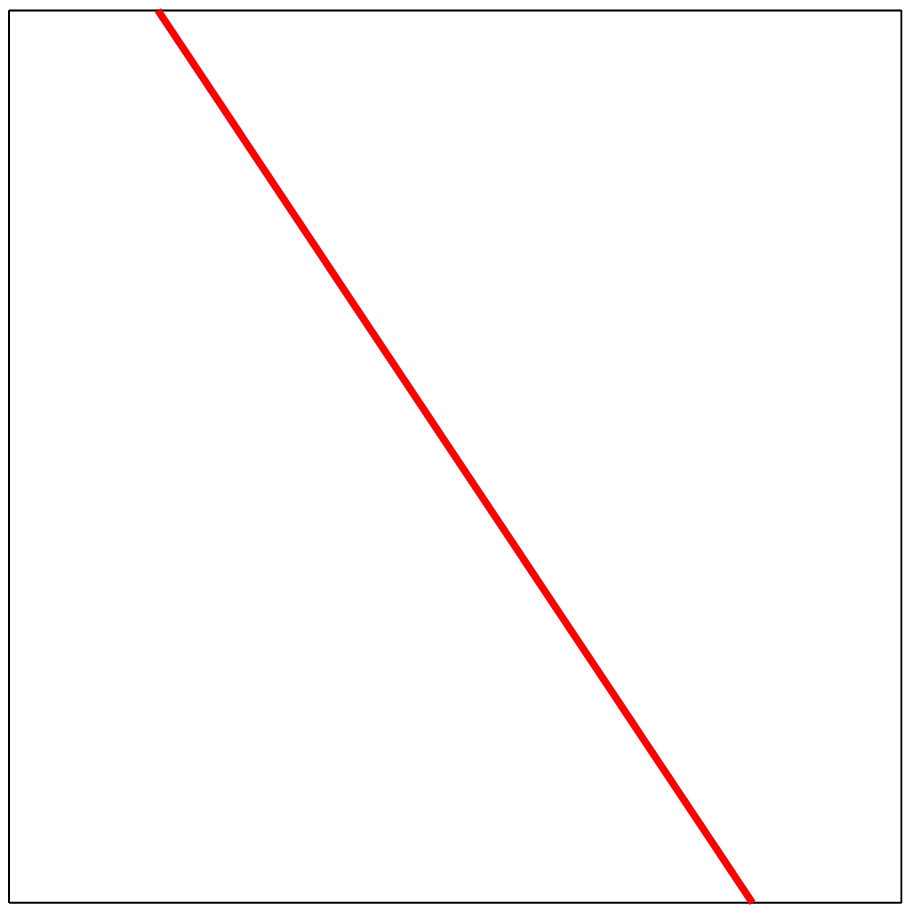}} \caption{Fracture settings with different $\theta$ of Example \ref{ex3}}\label{fig:ex3FractureAngles}
\end{figure}

\begin{ex} \label{ex3}
\textbf{Convergence test}

In this example, we test the convergence of our algorithm by the following problem:
$$  -\nabla\cdot({\bf{K}}\nabla p)=0,\quad x\in\Omega,$$
where ${\bf{K}} = {\bf{I}} + 2\delta(-\sin(\theta)x+\cos(\theta)y)
\begin{bmatrix}
    \cos^2(\theta) & \sin(\theta)\cos(\theta) \\
    \sin(\theta)\cos(\theta) & \sin^2(\theta)
\end{bmatrix}$, $\theta$ is an arbitrary fixed number and $\Omega=[-\pi,\pi]\times[-\pi,\pi]$.

One can verify that $$p(x,y)=\sin(\cos(\theta)x+\sin(\theta)y)e^{|-\sin(\theta)x+\cos(\theta)y|}$$ is the analytic solution of the above equation under corresponding Dirichlet boundary conditions for any fixed $\theta$. This problem can be interpreted as a single fracture whose angle is $\theta$ with the product of thickness and permeability $\epsilon k_f=2$ going through the square domain with permeability $k_m=1$ of porous matrix. By choosing different $\theta$'s and meshes, one can test the convergence of the algorithm comprehensively.
\end{ex}
We implemented three tests with different settings of $\theta$'s and rectangular meshes, in which an $\theta=5.3$ is randomly chosen to make the fracture intersecting with gridcells obliquely. See Figure \ref{fig:ex3FractureAngles} for details of settings of fracture. The numerical results are gathered in the Table \ref{tab:ex3}. From the results we can conclude that the NDFM algorithm is convergent, and it is confirmed by previous and subsequent numerical tests. Moreover, we find the rate of convergence is optimal under the conforming meshes and suboptimal under non-conforming meshes. These results suggest we choose conforming meshes if triangulation allows, but the choice of non-conforming meshes are always safe. What's more, based on the consistency of the algorithm shown in the Example \ref{ex2}, it's  reasonable to use a pseudo-conforming mesh, i.e. a non-conforming mesh which is relatively close to conforming meshes, to improve the accuracy of simulation, in the case that a conforming mesh is really hard to generate.

\begin{table}[!htb]
\centering
\caption{\label{tab:ex3} Convergence test for Example \ref{ex3}.}
\begin{tabular}{|c|c|cc|cc|cc|}
  \hline
& $N_x\times N_y$ &$||err||_{L^1(\Omega)}$ & order &$||err||_{L^2(\Omega)}$ & order&$||err||_{L^\infty(\Omega)}$ & order\\\hline
	&$20\times20$&  1.48E-00 &  --       &3.15E-01 & -- &2.32E-01 & -- \\
	&$40\times40$&  3.70E-01 & 2.00 & 7.88E-02 & 2.00 &5.82E-02 & 1.99 \\
$\theta=0$&$80\times80$& 9.24E-02 & 2.00 & 1.97E-02 & 2.00&1.45E-02 & 2.00 \\
    &$160\times160$& 2.31E-02 & 2.00 & 4.93E-03 & 2.00 &3.63E-03 & 2.00 \\
    &$320\times320$& 5.77E-03 &  2.00 & 1.23E-03 & 2.00&9.07E-04 & 2.00 \\
	&$640\times640$&1.44E-03  & 2.00  & 3.08E-04 & 2.00&2.27E-04 & 2.00 \\ \hline
	&$21\times21$& 1.92E-00 &  --       &3.72E-01 & -- &2.12E-01 & -- \\
	&$41\times41$&  6.43E-01 & 1.63 & 1.24E-01 & 1.64 &5.55E-02 & 2.00 \\
$\theta=0$&$81\times81$& 2.35E-01 & 1.48 & 4.67E-02 & 1.43&2.19E-02 & 1.37 \\
    &$161\times161$& 9.46E-02 & 1.32 & 2.00E-02 & 1.24&1.04E-02 & 1.09 \\
    &$321\times321$& 4.16E-02 & 1.19 & 9.24E-03 & 1.12&5.05E-03 & 1.05 \\
	&$641\times641$&1.94E-02  & 1.10  & 4.45E-03 & 1.06& 2.29E-03 & 1.02 \\ \hline
	&$20\times20$& 1.67E-00 &  --  &4.02E-01 & -- &4.18E-01 & -- \\
	&$40\times40$& 5.15E-01 & 1.70 & 1.11E-01 & 1.86 & 1.04E-01 & 2.01 \\
$\theta=5.3$&$80\times80$& 1.99E-01 &  1.37 & 4.10E-02 & 1.43&2.80E-02 & 1.89 \\
    &$160\times160$& 9.08E-02 & 1.13 & 2.02E-02 & 1.02 &1.63E-02 & 0.78 \\
    &$320\times320$& 4.81E-02 & 0.92 & 1.13E-02 & 0.83&9.15E-03 & 0.83 \\
	&$640\times640$&2.35E-02  & 1.04  & 5.72E-03 & 0.99&5.16E-03 & 0.83 \\ \hline
\end{tabular}
\end{table}

\begin{ex} \label{ex4}
\textbf{Hydrocoin}

This example is originally a benchmark for heterogeneous groundwater flow presented in the international Hydrocoin project \cite{Hydrocoin}.
The governing equation is
$$-\nabla\cdot({\bf{K}}\nabla H)=0,\quad x\in\Omega,$$
where H is the target variable hydraulic head, and ${\bf{K}}$ is the hydraulic conductivity.

A slight modification for geometrical parameters are made in \cite{Benchmark} and we will follow their settings and compare the results from NDFM with their benchmark reference
solution. One can refer to \cite{Hydrocoin,ThesisXFEM} for the numerical results based on the original geometry. The domain and boundary conditions used in this test are shown in Figure \ref{fig:ex4domain} with detailed coordinates listed in Table \ref{tab:ex4coord}. There are two fractures crossing through porous matrix, with central axis 2-7, 4-8 and thickness $\epsilon_1=5\sqrt{2},\epsilon_2=33\sqrt{5}/5$, respectively. The hydraulic conductivity ${\bf{K}}$ is $10^{-8}$ \si[per-mode=symbol]{\meter\per\second} in porous matrix and $10^{-6}$ \si[per-mode=symbol]{\meter\per\second} in fractured region. Moreover, as shown in Figure \ref{fig:ex4domain}, the left, right and bottom boundaries are impervious, i.e. $q_N=0$, while the top one is a Dirichlet boundary with $H=$ height, i.e. the $z$ coordinate.
\begin{figure}[htbp]
\centering
\includegraphics[width=4.5in]{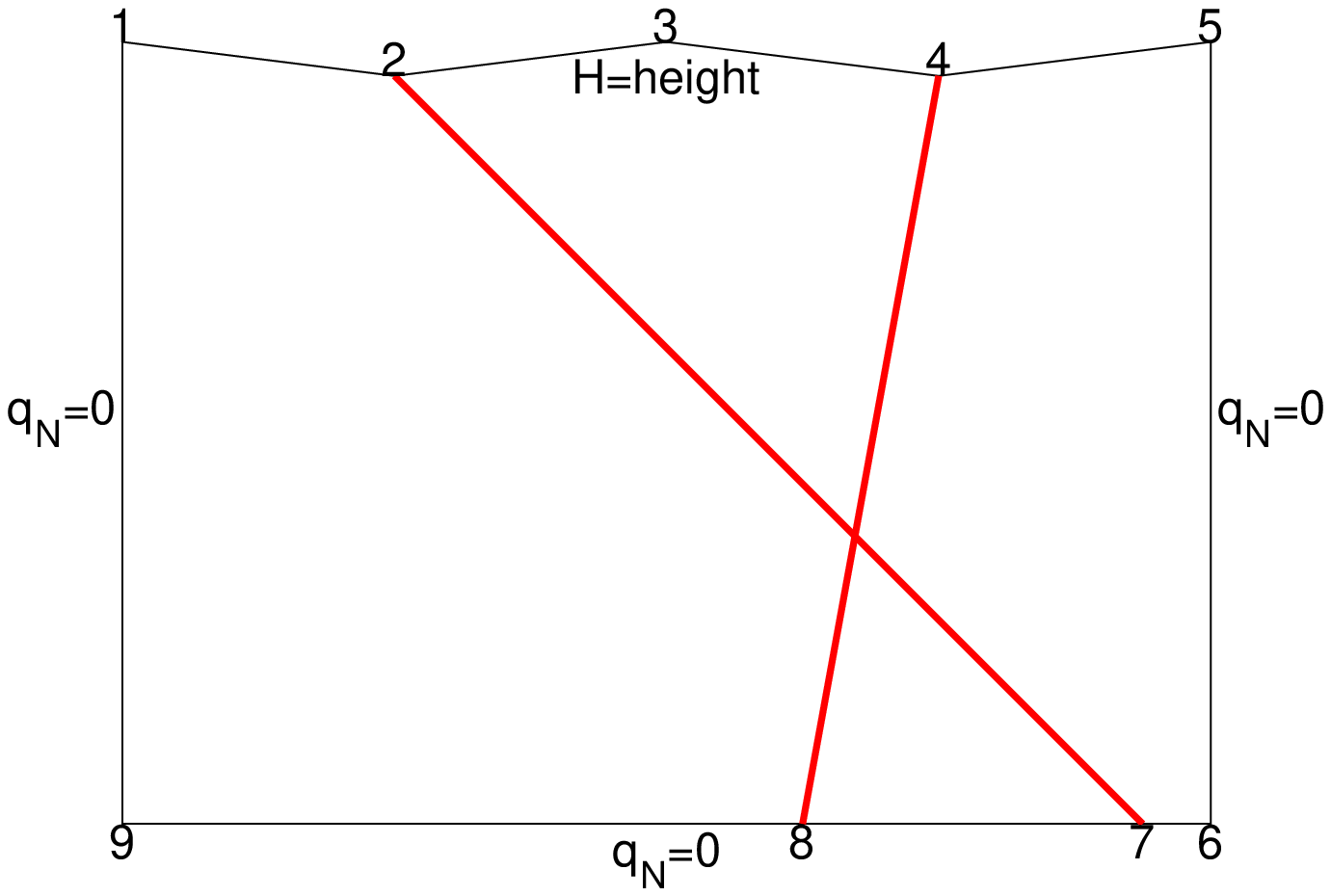}
\caption{Domain and boundary conditions of Example \ref{ex4}} \label{fig:ex4domain}
\end{figure}

\begin{table}[!htb]
\centering
\caption{\label{tab:ex4coord} Coordinates of the numbers labeled in Figure \ref{fig:ex4domain}}
\begin{tabular}{c c c|c c c}
  \hline
point & $x$ (\si{m}) & $z$ (\si{m}) & point & $x$ (\si{m}) & $z$ (\si{m})\\\hline
1&0& 150 &6&1600&-1000 \\\hline
2&400&100&7&1500&-1000  \\\hline
3&800&150&8&1000&-1000\\\hline
4&1200&100&9&0&-1000 \\ \hline
5&1600&150& - & - & -\\ \hline
\end{tabular}
\end{table}
\end{ex}

\begin{figure}[!htbp]
\subfigure[Triangulation of the domain (coarse mesh)]{\includegraphics[width = 3in]{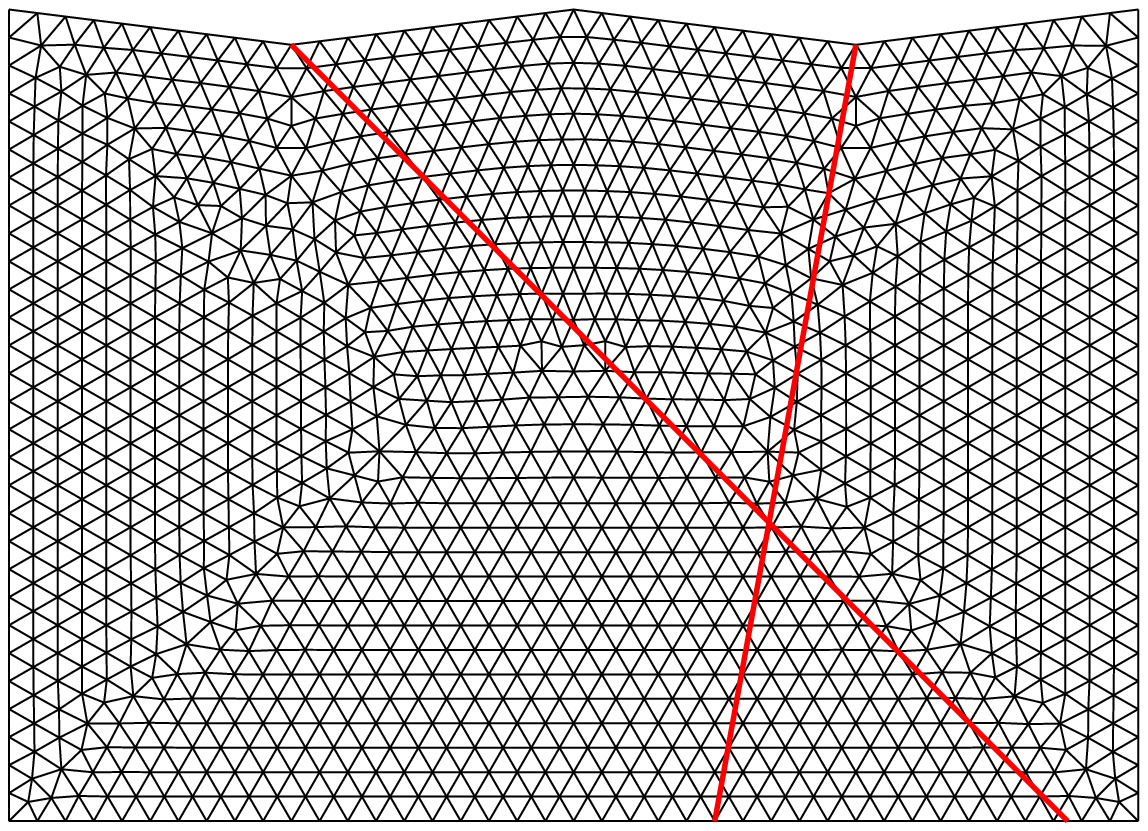}} \subfigure[Numerical result (coarse mesh)]{\includegraphics[width = 3in]{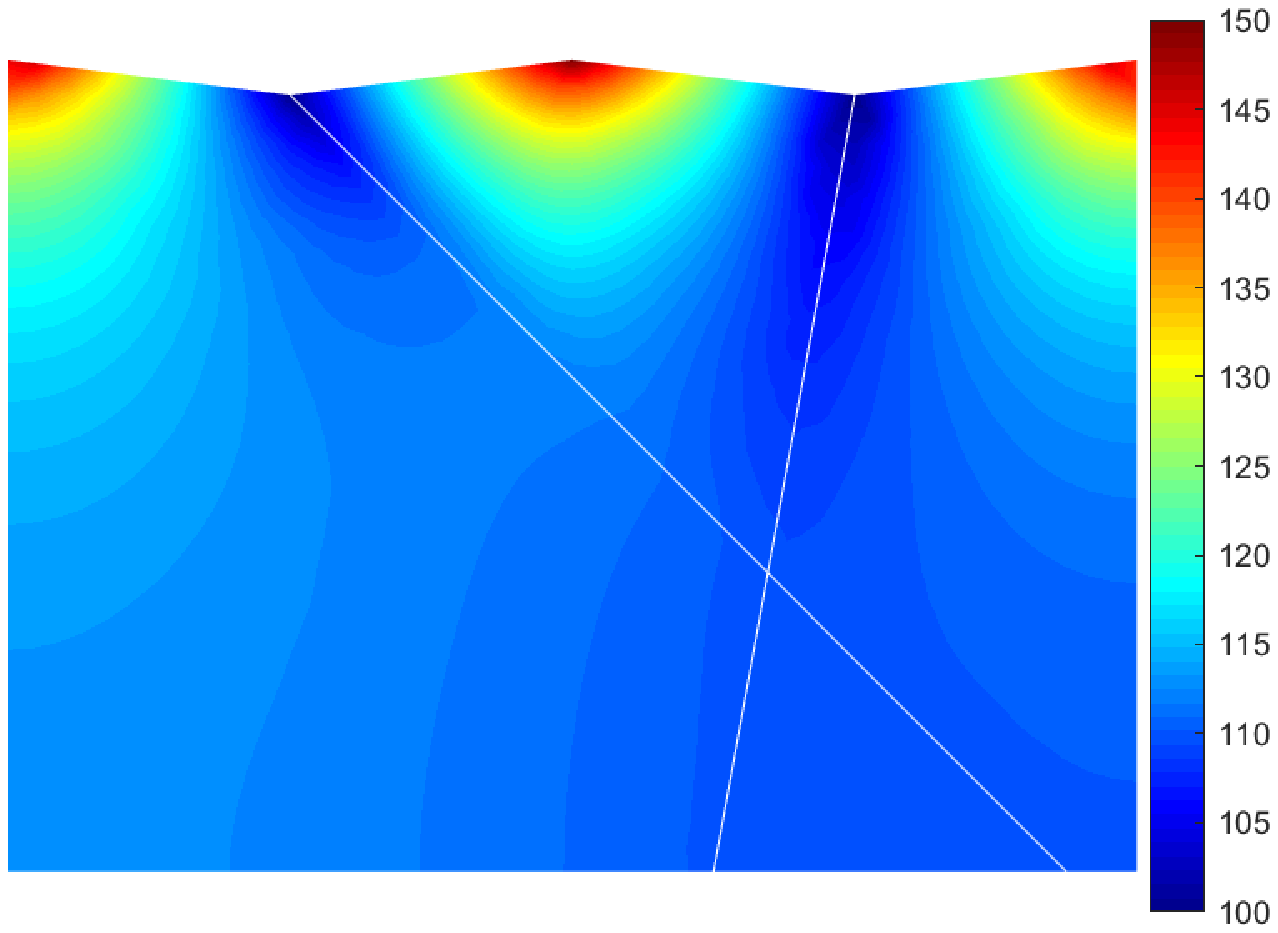}}\\
\subfigure[Triangulation of the domain (fine mesh)]{\includegraphics[width = 3in]{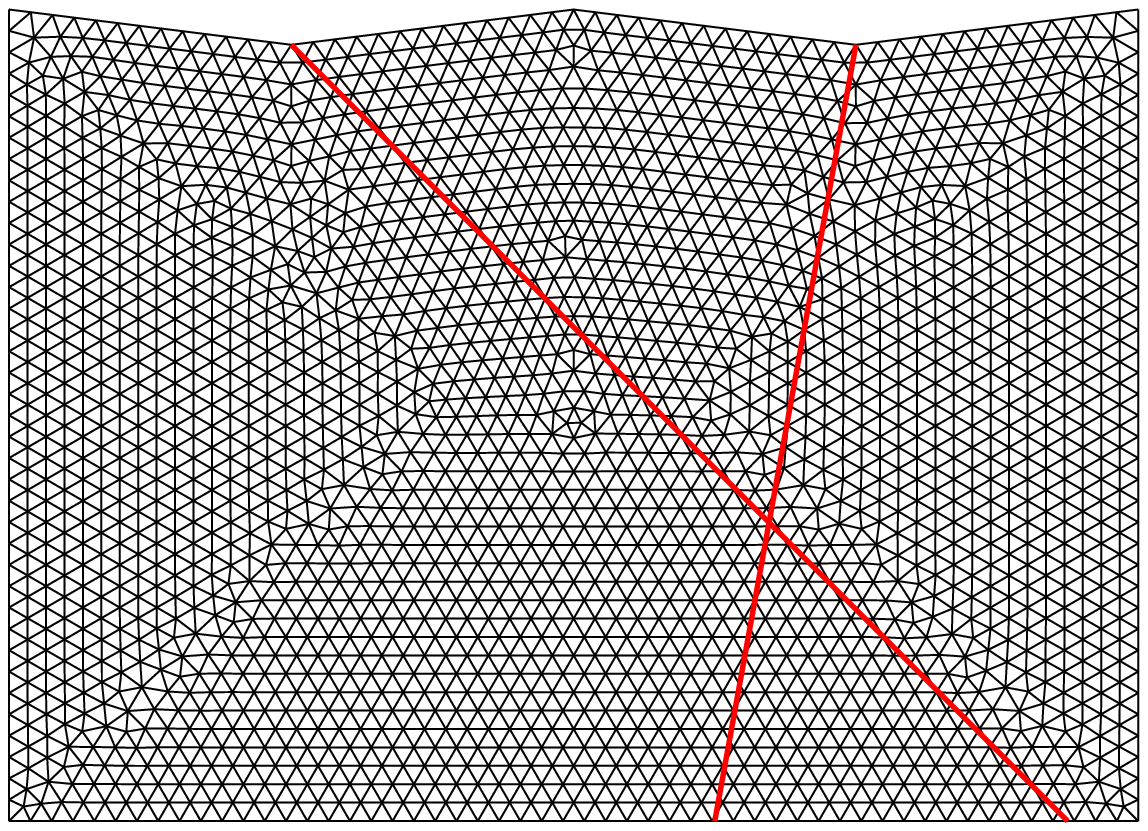}} \subfigure[Numerical result (fine mesh)]{\includegraphics[width = 3in]{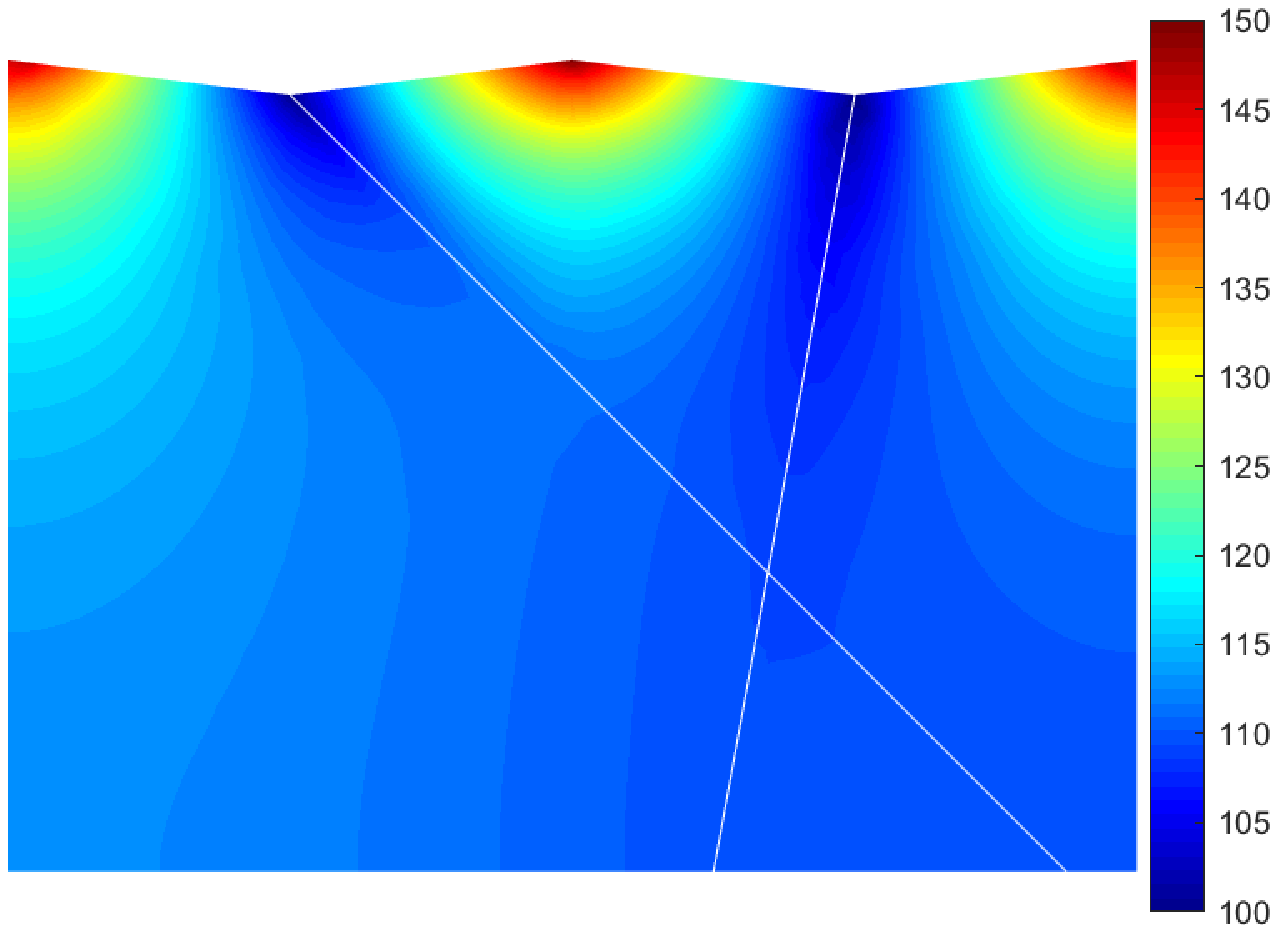}}
\caption{Triangulation and numerical results of Example \ref{ex4}\label{fig:ex4MeshSolution}}
\end{figure}

To demonstrate the performance of the NDFM more comprehensively, we employ this method on two triangular meshes with different grid sizes. The coarse mesh contains 2612 triangular elements while the fine mesh contains 4511 triangular elements. The meshes and contour plots of numerical results under different triangulation are presented in Figure \ref{fig:ex4MeshSolution}.

To evaluate our method, we need to compare its results with the reference solution, as well as solutions from other methods. A description of the participants in this test is given as follows.
The reference solution is provided by mimetic finite difference (MFD) method for equi-dimensional model problem \eqref{PoissonEq}, \eqref{EquiDimensionK}, \eqref{BVP} on a very fine mesh containing 424921 matrix elements and 19287 fracture elements, with total degrees of freedom (d.o.f) 889233. Other models and methods participate in this example for comparison and evaluation are vertex-centered control volume discrete fracture model (Box-DFM), cell-centered two point flux approximation control volume discrete fracture model (CC-DFM), embedded discrete fracture model (EDFM), mortar-flux discrete fracture model (Mortar-DFM), primal extended finite element method (P-XFEM), and dual extended finite element method (D-XFEM). Some of them are restricted on conforming meshes while the others can be employed on non-conforming meshes. For more introductions of above methods, see \cite{Benchmark}.

\begin{figure}[!htbp]
\subfigure[Slices of hydraulic head (coarse mesh)]{\includegraphics[width = 3in]{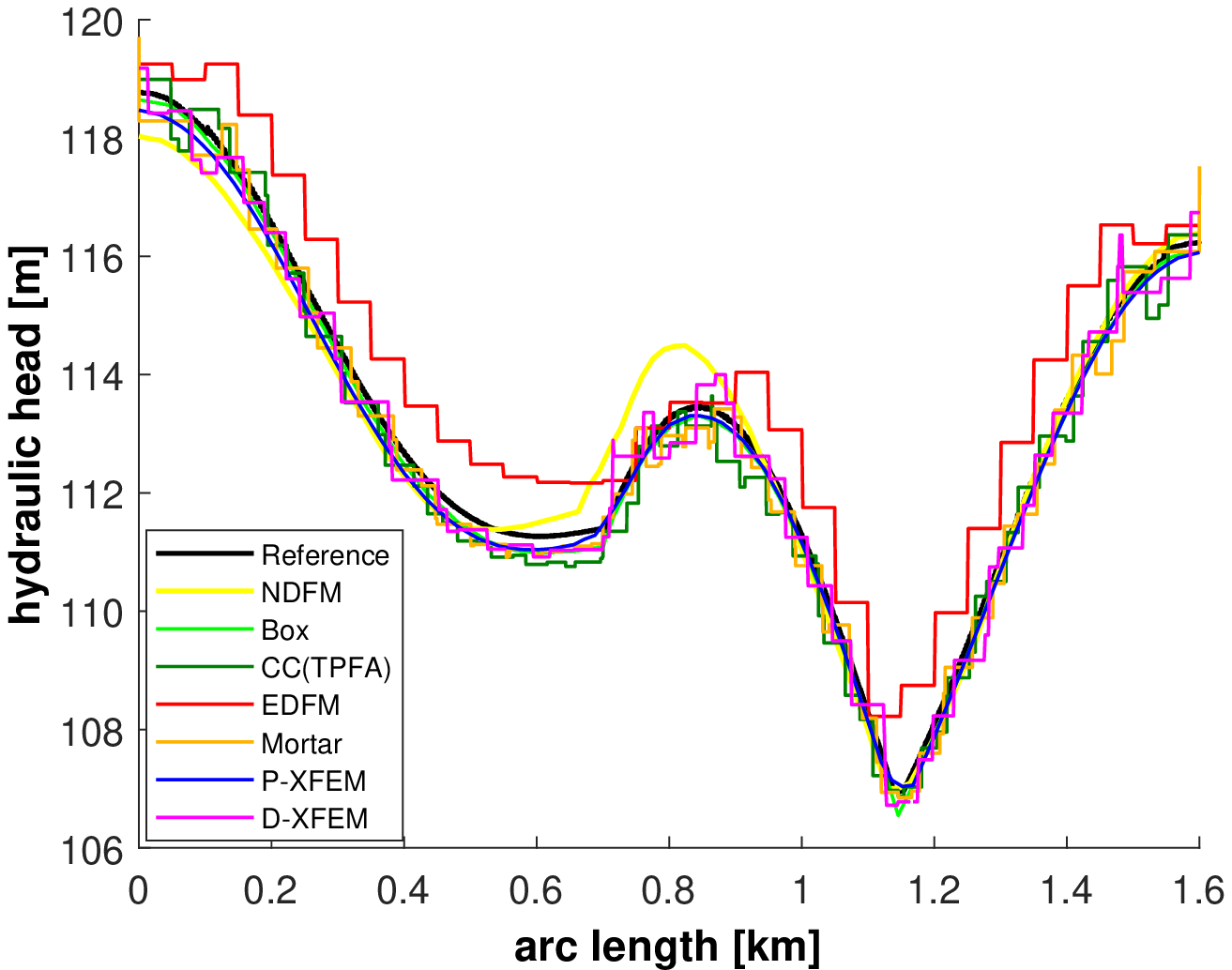}}
\subfigure[Slices of hydraulic head (fine mesh)]{\includegraphics[width = 3in]{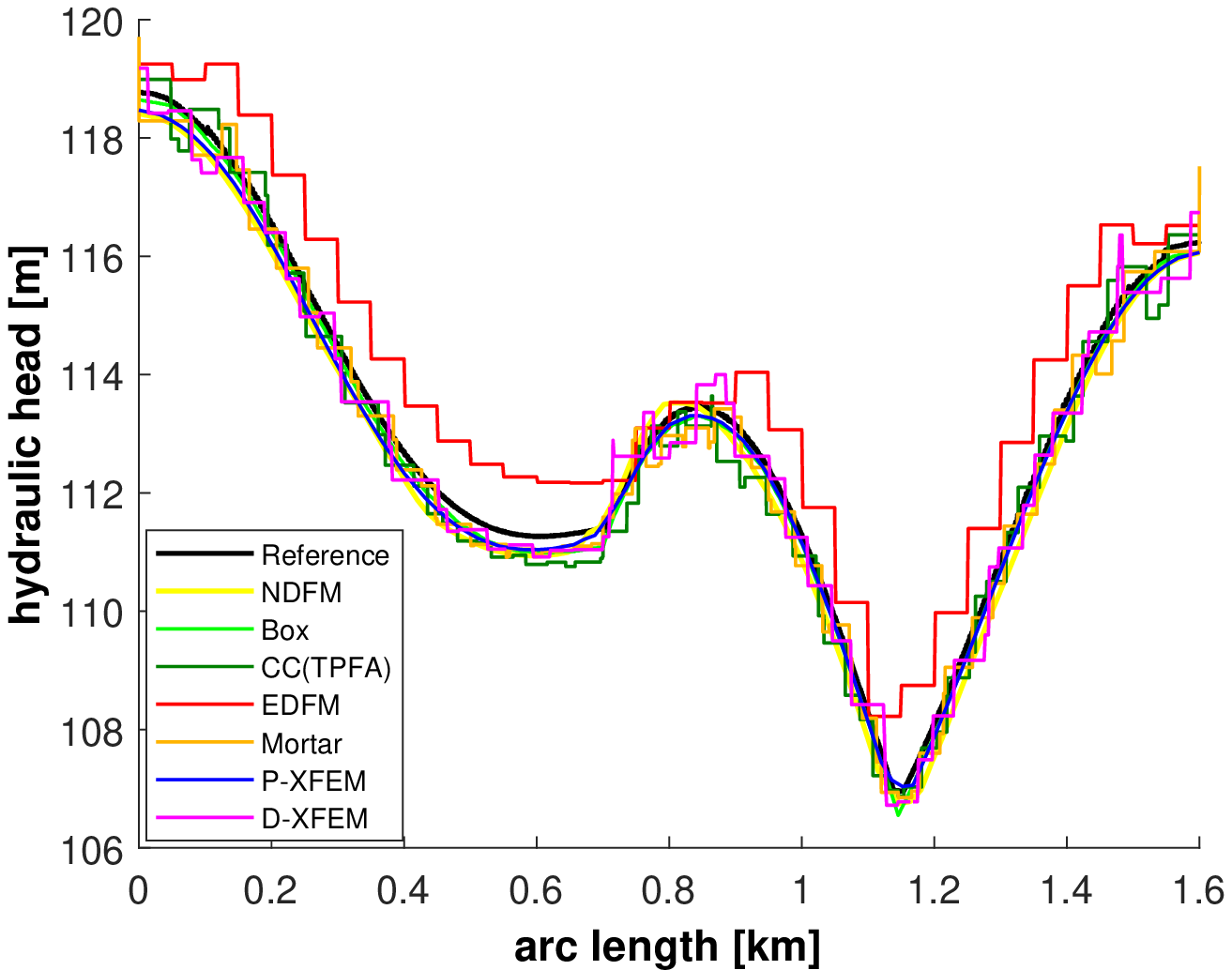}}
\caption{Comparison of slices of hydraulic head along a horizontal line $z=-200$\si{m} of Example \ref{ex4}}\label{fig:ex4Slices}
\end{figure}

We slice profiles of hydraulic head along the horizontal line $z=-200\si{\m}$ on the coarse and fine meshes and plot them in Figure \ref{fig:ex4Slices} (a) and (b), respectively, together with the slices of reference solution and solutions of other methods.
From the plots we can see the profile on coarse mesh has a little derivation from the reference solution but the profile from fine mesh match it well.

In addition to the comparison of the slices along specific line, another way to evaluate the accuracy of a method is to compute the relative error of numerical solutions based on the reference solution of MFD on fine mesh. Following \cite{Benchmark}, the relative $L^2$ errors of solutions on matrix and fractures are computed using the formulas given below.
\begin{equation}\label{ErrmFormula}
    err^{2}_\text{m} = \frac{1}{|\Omega|(\Delta p_\text{ref})^{2}} \sum_{i,j} |T^{i}_\text{ref} \cap T^{j}_\text{m}|\left(p_\text{ref}|_{T^{i}_\text{ref}}-p_\text{m}  |_{T^{j}_\text{m}}\right)^{2},
\end{equation}
\begin{equation}\label{ErrfFormula}
        err^{2}_\text{f} = \frac{1}{|\Gamma|(\Delta p_\text{ref})^{2}} \sum_{i,l} |T^{i}_\text{ref} \cap T^{l}_\text{f}|\left(p_\text{ref}|_{T^{i}_\text{ref}}-p_\text{f}  |_{T^{l}_\text{f}}\right)^{2},
\end{equation}
where $\Omega$ and $\Gamma$ are the matrix and fracture region, respectively, $|\cdot|$ means 2-D measure in \eqref{ErrmFormula} and 1-D measure in $\eqref{ErrfFormula}$, $\Delta p_\text{ref}=\max{p_\text{ref}}-\min{p_\text{ref}}$ is the range of reference solution, $T^{i}_\text{ref} , i=1,2,\ldots,I,$ are the fine elements used in the reference solution, $T^{j}_\text{m}, j=1,2,\ldots,J,$ and
$T^{l}_\text{f}, l=1,2,\ldots,L,$ are the matrix elements and fracture elements in the methods to be evaluated, respectively.

A summary of errors of different methods on matrix and fractures are listed in Table \ref{tab:ex4}, together with some other important aspects of the methods, such as the requirement for meshes, degrees of freedom (d.o.f), sparsity and conditional number ($||\cdot||_2$-cond) of the resulting system of linear equations.

\begin{table}[!htb]
\centering
\caption{\label{tab:ex4} Evaluation data of different algorithms in Example \ref{ex4}.}
\begin{tabular}{|c|c|c|c|c|c|c|}
  \hline
method & d.o.f & mesh & err$_m$ &err$_f$&sparsity&$||\cdot||_2$-cond\\\hline
Box-DFM&$1496$& conforming & 9.3E-03       &3.3E-03 & 4.5\textperthousand &5.4E03 \\\hline
CC-DFM&$1459$& conforming & 1.1E-02 & 1.1E-02 & 2.7\textperthousand &3.5E04  \\\hline
EDFM&$1044$& non-conforming & 1.5E-02 & 8.3E-03 & 4.7\textperthousand &3.9E04\\\hline
Mortar-DFM&$3647$& conforming & 1.0E-02  & 7.2E-03 & 1.5\textperthousand & 9.0E12 \\ \hline
P-XFEM &$1667$& non-conforming & 1.2E-02     & 3.2E-03 & 6.5\textperthousand &2.7E09\\ \hline
D-XFEM&3514& non-conforming & 1.2E-02  & 6.9E-03 & 1.7\textperthousand & 6.2E12\\ \hline
NDFM& 1337 &non-conforming &1.1E-02 & 1.2E-02 & 5.1\textperthousand &7.0E04\\ \hline
NDFM& 2296 &non-conforming &8.8E-03 & 6.1E-03 & 3.0\textperthousand &1.3E05\\ \hline
\end{tabular}
\end{table}

\begin{figure}[!htbp]
\centering
\includegraphics[width=4.5in]{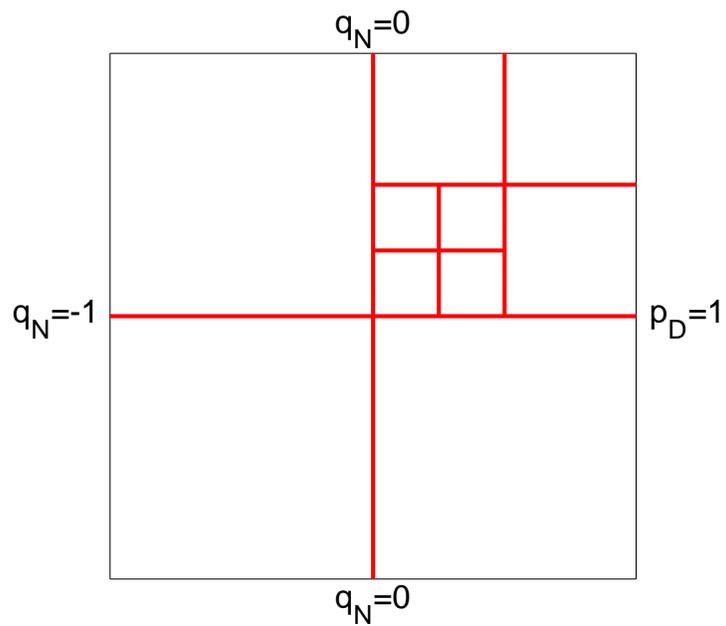}
\caption{Domain and boundary conditions of Example \ref{ex5}} \label{fig:ex5DomainSolution}
\end{figure}

\begin{ex} \label{ex5}
\textbf{Regular Fracture Network}

This test case is originally from \cite{benchmark3} and modified by \cite{Benchmark}, which simulates a regular fracture network in a square porous media. The domain $\Omega=[0,1]\times[0,1]$. The central axis of fractures are $y=0.5, x=0.5, x=0.75, y=0.75, x=0.625, y=0.625$, respectively, and all fractures have a uniform thickness $\epsilon = 10^{-4}$. See Figure \ref{fig:ex5DomainSolution} as an illustration. The permeability is $1$ in porous matrix and $10^{4}$ in fractured region. Moreover, the left boundary is an inflow boundary with $q_N=-1$, the right boundary is a Dirichlet boundary with pressure $p_D=1$, and the top and bottom boundaries are impervious, i.e. $q_N=0$.
\end{ex}

\begin{figure}[!htbp]
\subfigure[Numerical result on rectangular mesh $N_x=N_y=25$]{\includegraphics[width = 3in]{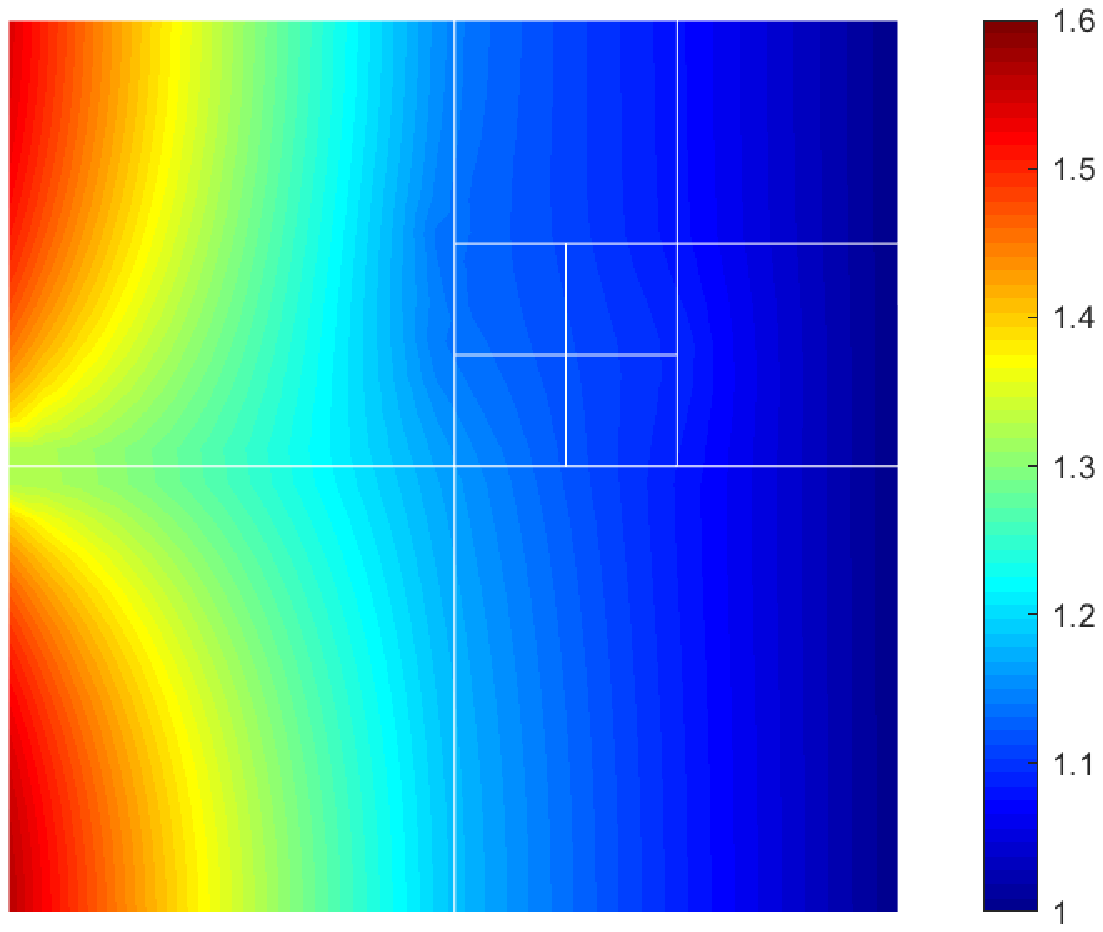}}
\subfigure[Numerical result on rectangular mesh $N_x=N_y=35$]{\includegraphics[width = 3in]{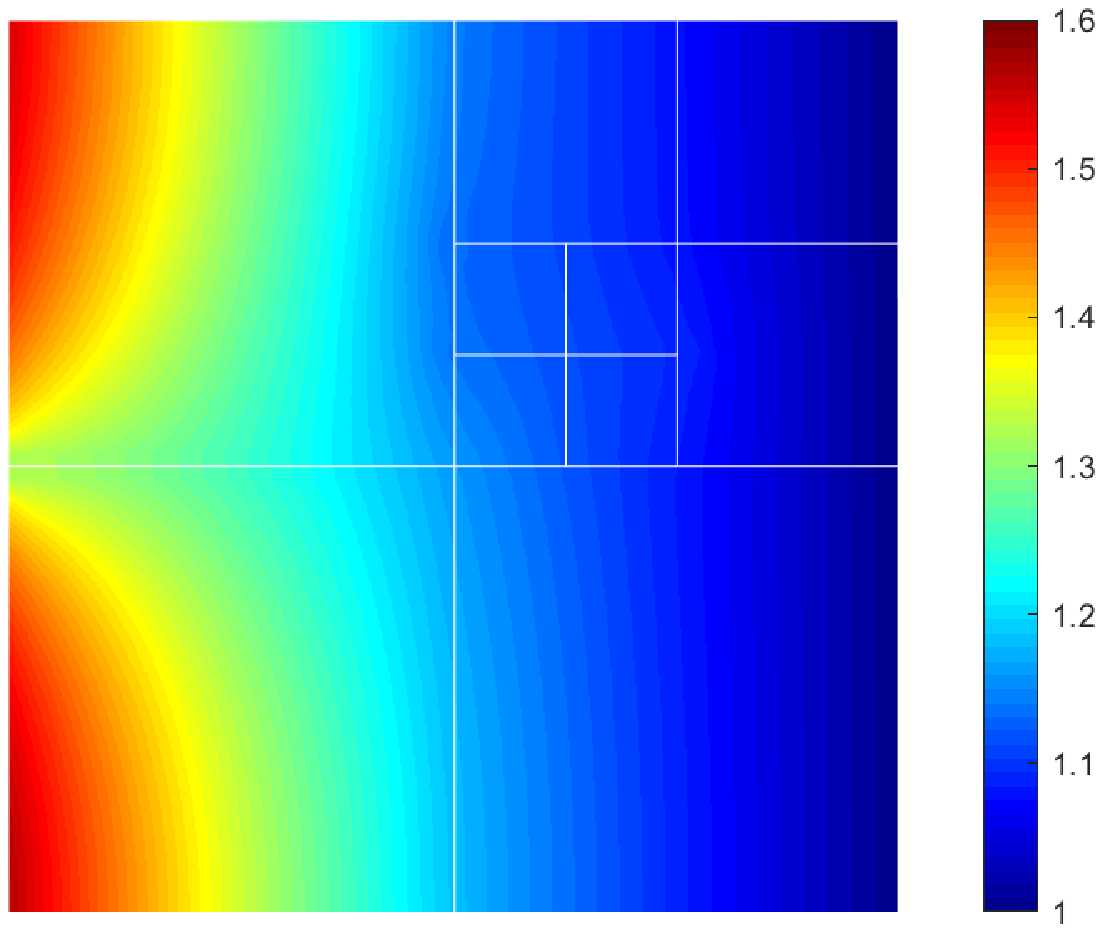}}
\caption{Numerical results of Example \ref{ex5} on different meshes}\label{fig:ex5Solutions}
\end{figure}

\begin{figure}[!htbp]
\subfigure[Slice of pressure along $y=0.7$ (coarse mesh)]{\includegraphics[width = 3in]{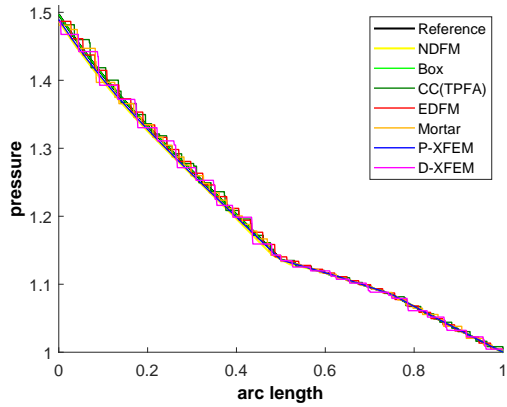}}
\subfigure[Slice of pressure along $y=0.7$ (fine mesh)]{\includegraphics[width = 3in]{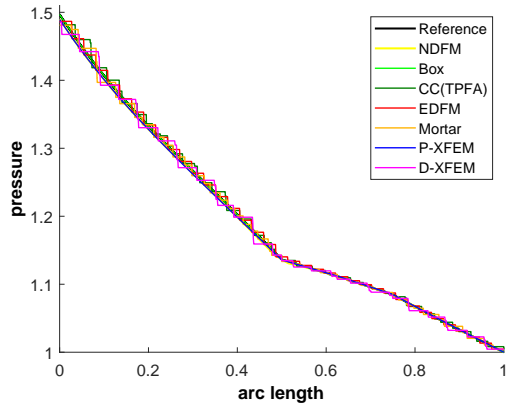}}\\
\subfigure[Slice of pressure along $x=0.5$ (coarse mesh)]{\includegraphics[width = 3in]{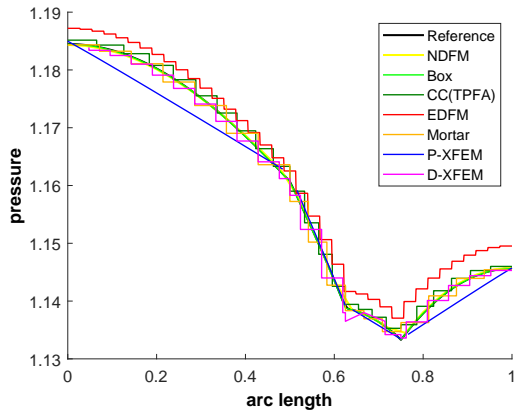}}
\subfigure[Slice of pressure along $x=0.5$ (fine mesh)]{\includegraphics[width = 3in]{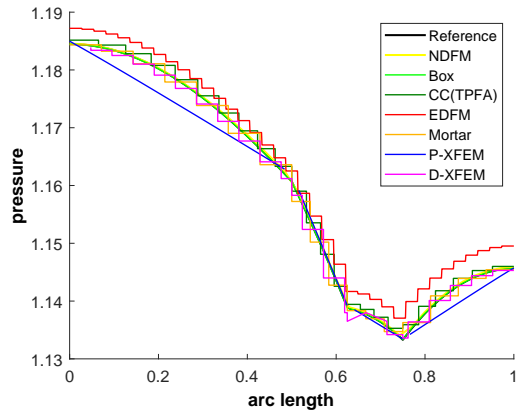}}
\caption{Comparison of slices of pressure along lines $y=0.7$ and $x=0.5$ of Example \ref{ex5} on different meshes}\label{fig:ex5Slices}
\end{figure}

We employ the NDFM on two rectangular meshes of different grid sizes. The coarse mesh contains $25\times25$ subrectangles while the fine mesh has $35\times35$ ones. The contour plots of the numerical results are presented in Figure \ref{fig:ex5Solutions}.
What's more, we slice the profiles of pressure along the horizontal line $y=0.7$  and vertical line $x=0.5$ on the coarse and fine meshes, respectively, and plot them in Figure \ref{fig:ex5Slices} together with the slices of reference solution and solutions of other methods. The reference solution is provided by MFD method for equi-dimensional model on a very fine nonuniform grid containing 1136456 matrix elements and 38600 fracture elements, with total d.o.f 2352280. As we can see from the figures, the profiles match the reference solution perfectly.
In addition to the slices, we summarize the relative errors on matrix and fractures, as well as other important aspects of different methods in Table \ref{tab:ex5}. As we can see, the NDFM methods successfully gain very accurate solution though the meshes used are coarse.
\begin{table}[!htb]
\centering
\caption{\label{tab:ex5} Evaluation data of different algorithms in Example \ref{ex5}.}
\begin{tabular}{|c|c|c|c|c|c|c|}
  \hline
method & d.o.f & mesh & err$_m$ &err$_f$&sparsity&$||\cdot||_2$-cond\\\hline
Box-DFM&$1422$& conforming & 6.7E-03&1.1E-03 & 4.7\textperthousand &7.9E03 \\\hline
CC-DFM&$1481$& conforming & 1.1E-02 & 5.0E-03 & 2.7\textperthousand &5.6E04  \\\hline
EDFM&$1501$& non-conforming & 6.5E-03 & 4.0E-03 & 3.3\textperthousand &5.6E04\\\hline
Mortar-DFM&$3366$& conforming & 1.0E-02  & 7.4E-03 & 1.8\textperthousand & 2.4E06 \\ \hline
P-XFEM &$1632$& non-conforming & 1.7E-02     & 6.0E-03 & 7.8\textperthousand &6.8E09\\ \hline
D-XFEM&$4474$& non-conforming & 9.6E-03  & 8.9E-03 & 1.3\textperthousand & 1.2E06\\ \hline
NDFM&$650$& non-conforming & 1.3E-02  & 8.9E-03 & 13.1\textperthousand & 1.8E04\\ \hline
NDFM&$1260$& non-conforming & 8.8E-03  & 6.4E-03 & 6.9\textperthousand & 6.4E04\\ \hline
\end{tabular}
\end{table}

\begin{figure}[!htbp]
\centering
\includegraphics[width=4.5in]{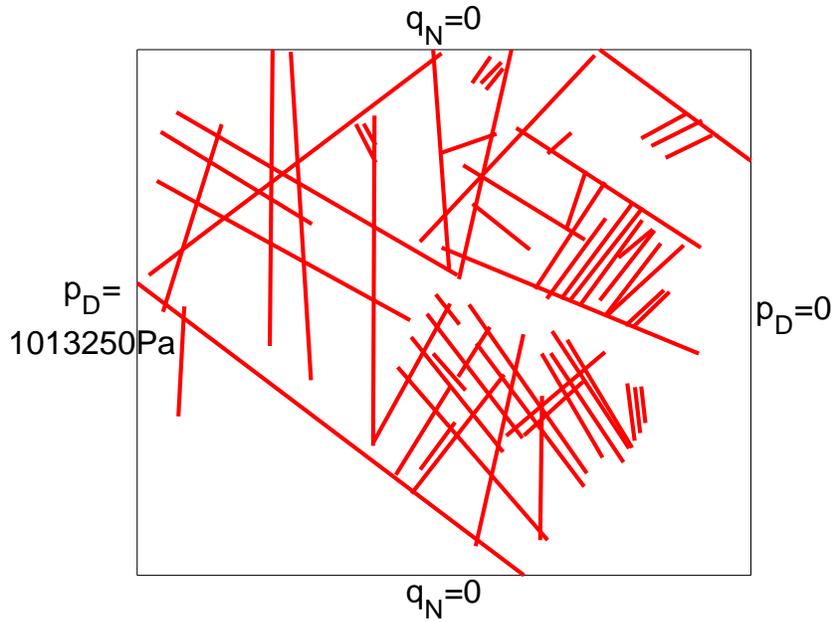}
\caption{Domain and boundary conditions of Example \ref{ex6}} \label{fig:ex6DomainSolution}
\end{figure}

\begin{ex} \label{ex6}
\textbf{a Realistic Case}

This example is a benchmark problem in \cite{Benchmark} modified from a real set of fractures from an interpreted outcrop in
the Sotra island. One can also refer to the simulation results in \cite{benchmark2} for reference. The test case is a complex fracture network containing 63 fractures with different lengths and connectivity. The domain is set to be $\Omega=[0,700\si{m}]\times[0,600\si{m}]$ with permeability $10^{-14}$\si{\m^2}. The fractures on it are shown in Figure \ref{fig:ex6DomainSolution} (a), with uniform permeability $10^{-8}$\si{\m^2} and thickness $\epsilon=10^{-2}$\si{\m}. The detailed geometric data of fractures is attached in appendix. The top and bottom boundary are impervious, i.e. $q_N=0$ while the left and right boundary are Dirichlet boundary with pressure $p_D=1013250$ \si{Pa} and $p_D=0$, respectively.
\end{ex}

\begin{figure}[!htbp]
\subfigure[Numerical result on rectangular mesh $N_x=105,N_y=90$ ]{\includegraphics[width = 3in]{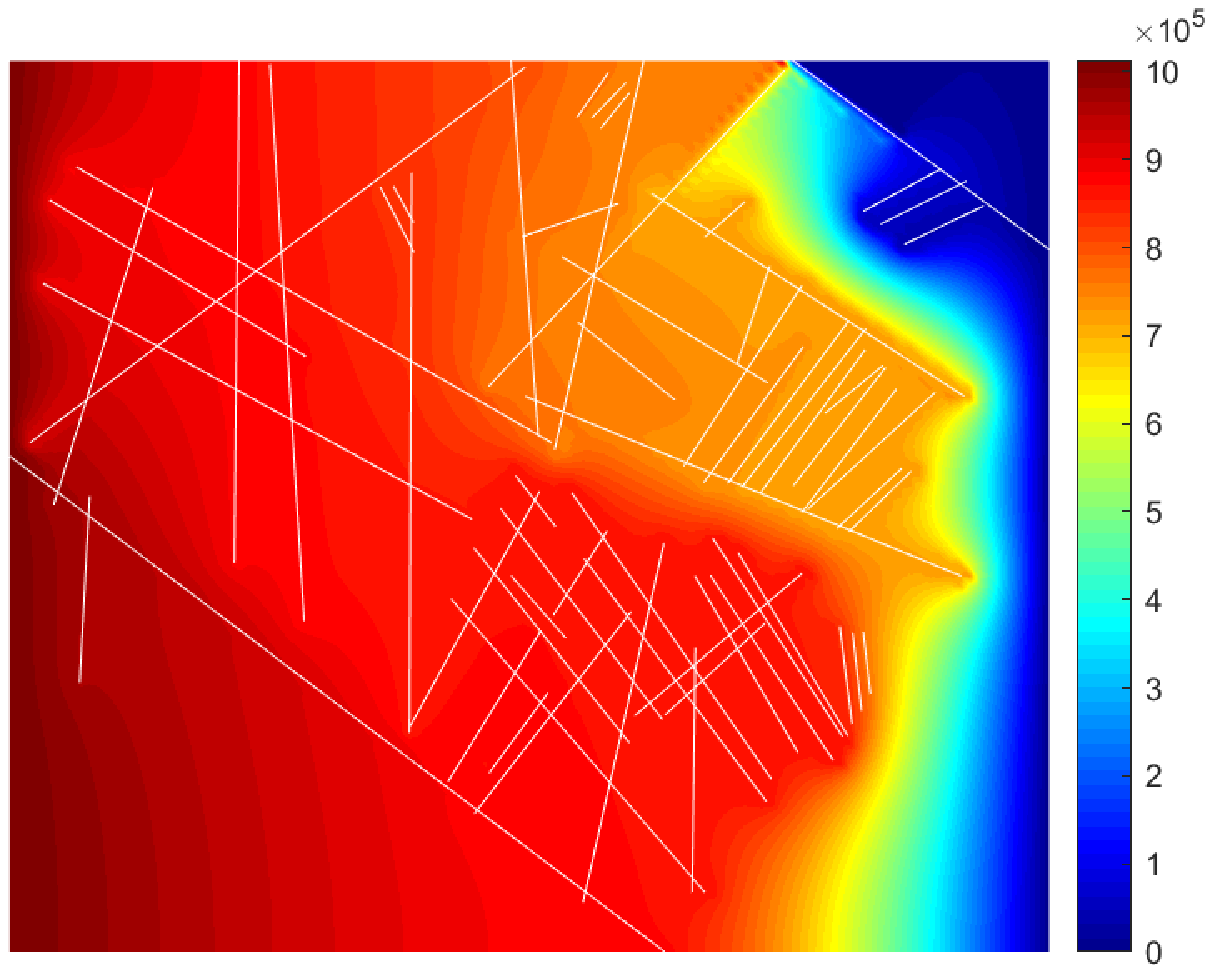}}
\subfigure[Numerical result on rectangular mesh $N_x=175,N_y=150$]{\includegraphics[width = 3in]{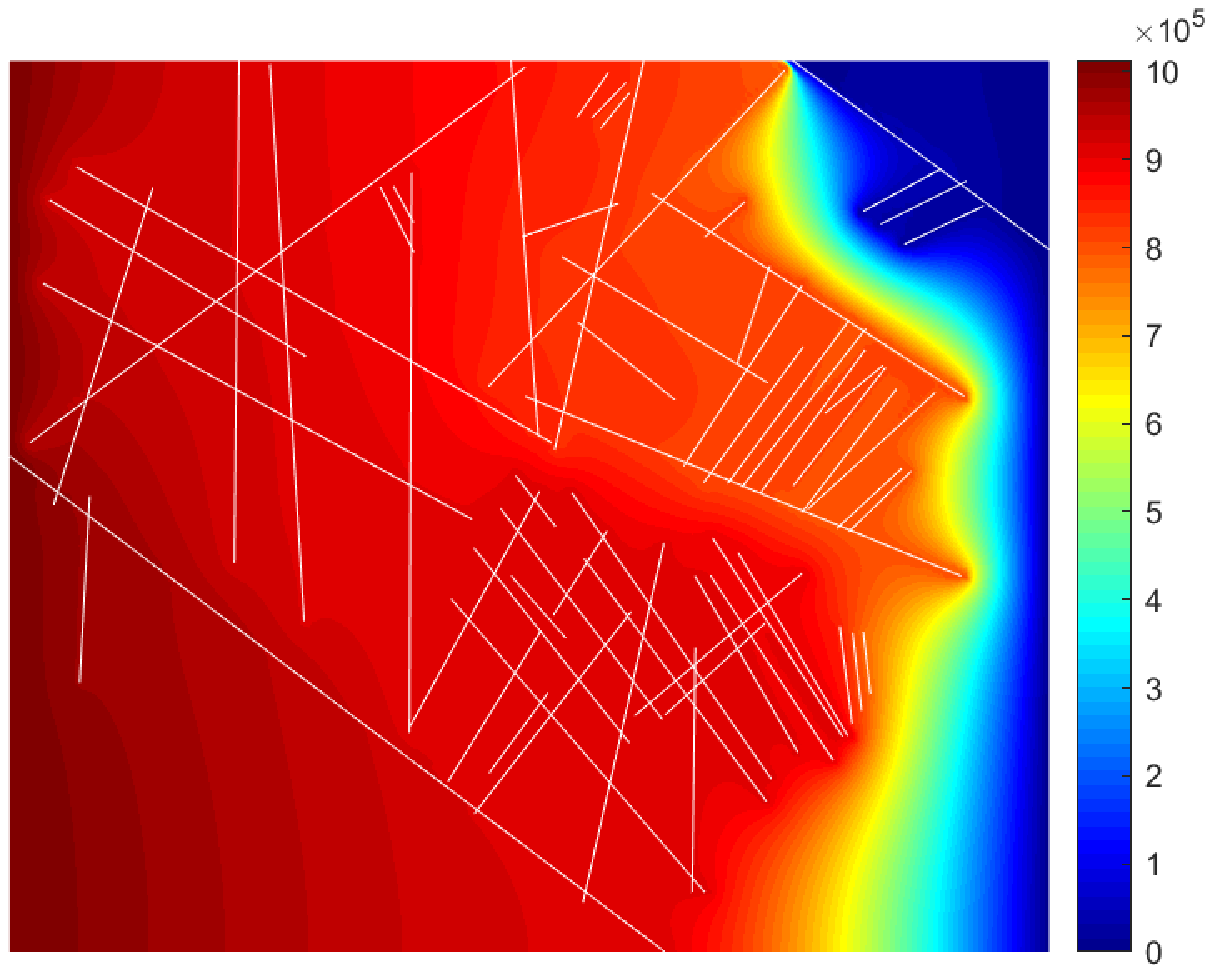}}
\caption{Numerical results of Example \ref{ex6} on different meshes}\label{fig:ex6Solutions}
\end{figure}

We implement the NDFM on a coarse mesh containing $105\times90$ rectangular elements and a fine mesh containing $175\times150$ rectangular elements. The contour plots of numerical results are presented in Figure \ref{fig:ex6Solutions}.

\begin{figure}[!htbp]
\subfigure[Slice of pressure along $y=500$\si{\m} (d.o.f=9464)]{\includegraphics[width = 3in]{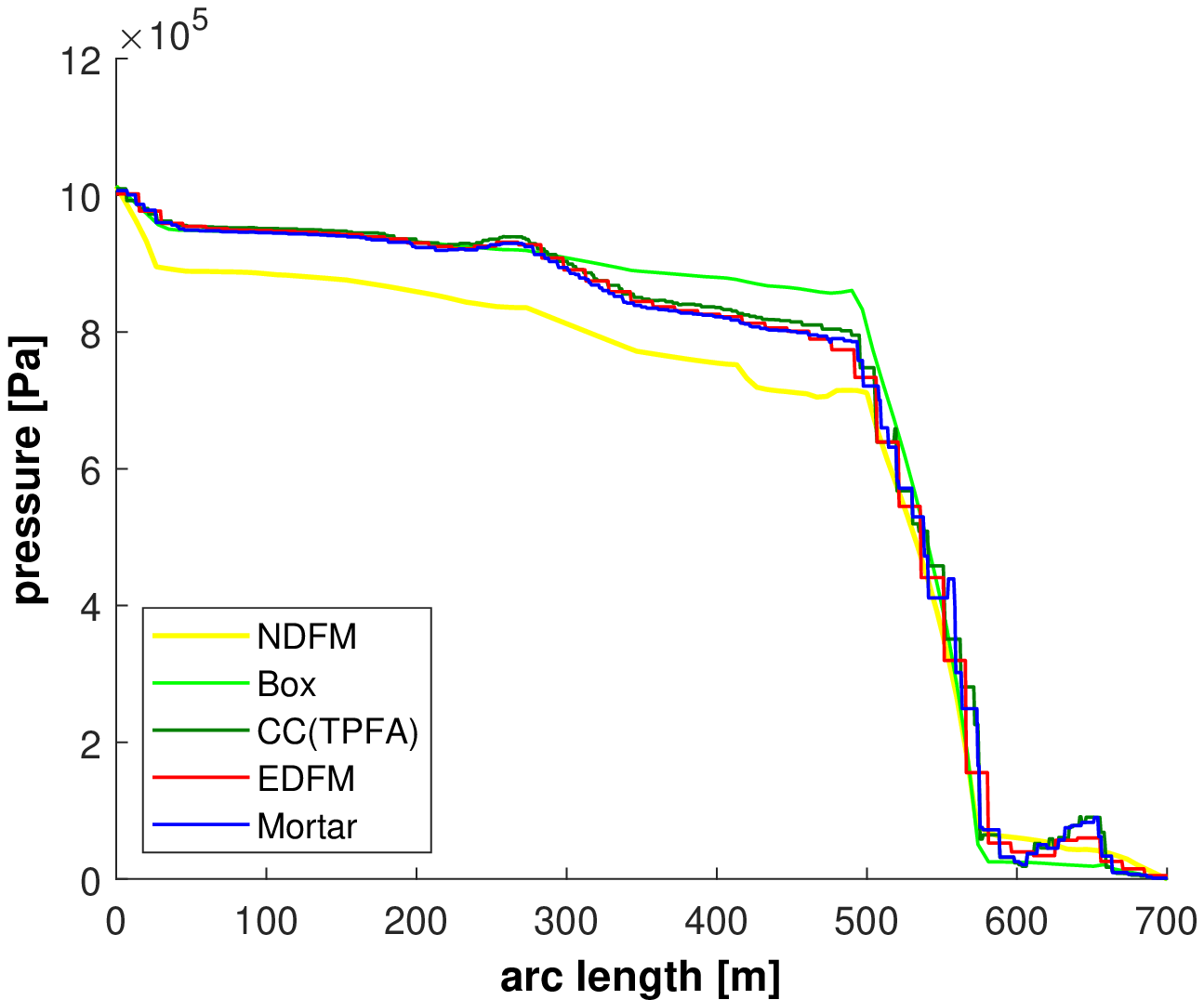}}
\subfigure[Slice of pressure along $y=500$\si{m} (d.o.f=26274)]{\includegraphics[width = 3in]{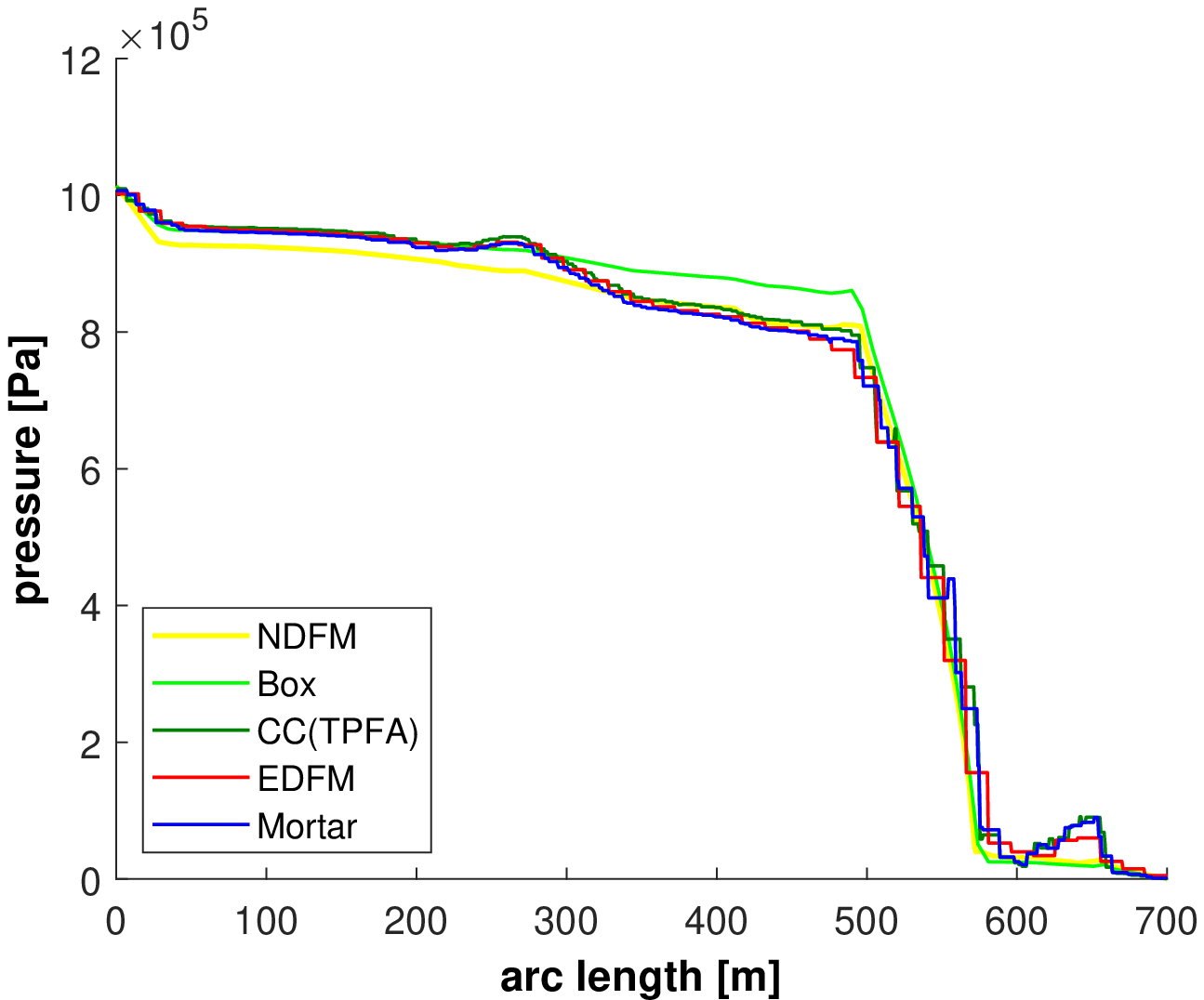}}\\
\subfigure[Slice of pressure along $x=625$\si{\m} (d.o.f=9464)]{\includegraphics[width = 3in]{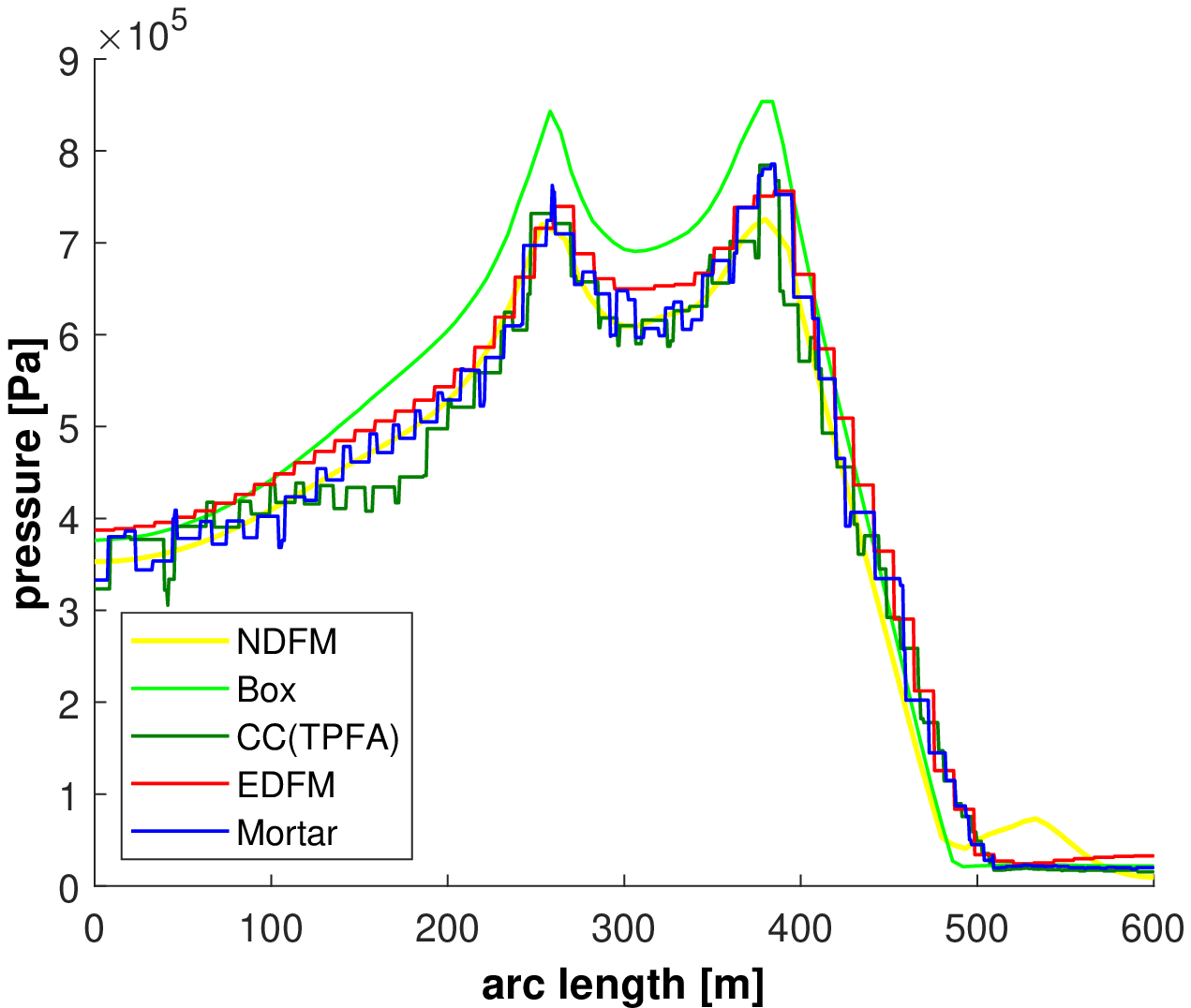}}
\subfigure[Slice of pressure along $x=625$\si{\m} (d.o.f=26274)]{\includegraphics[width = 3in]{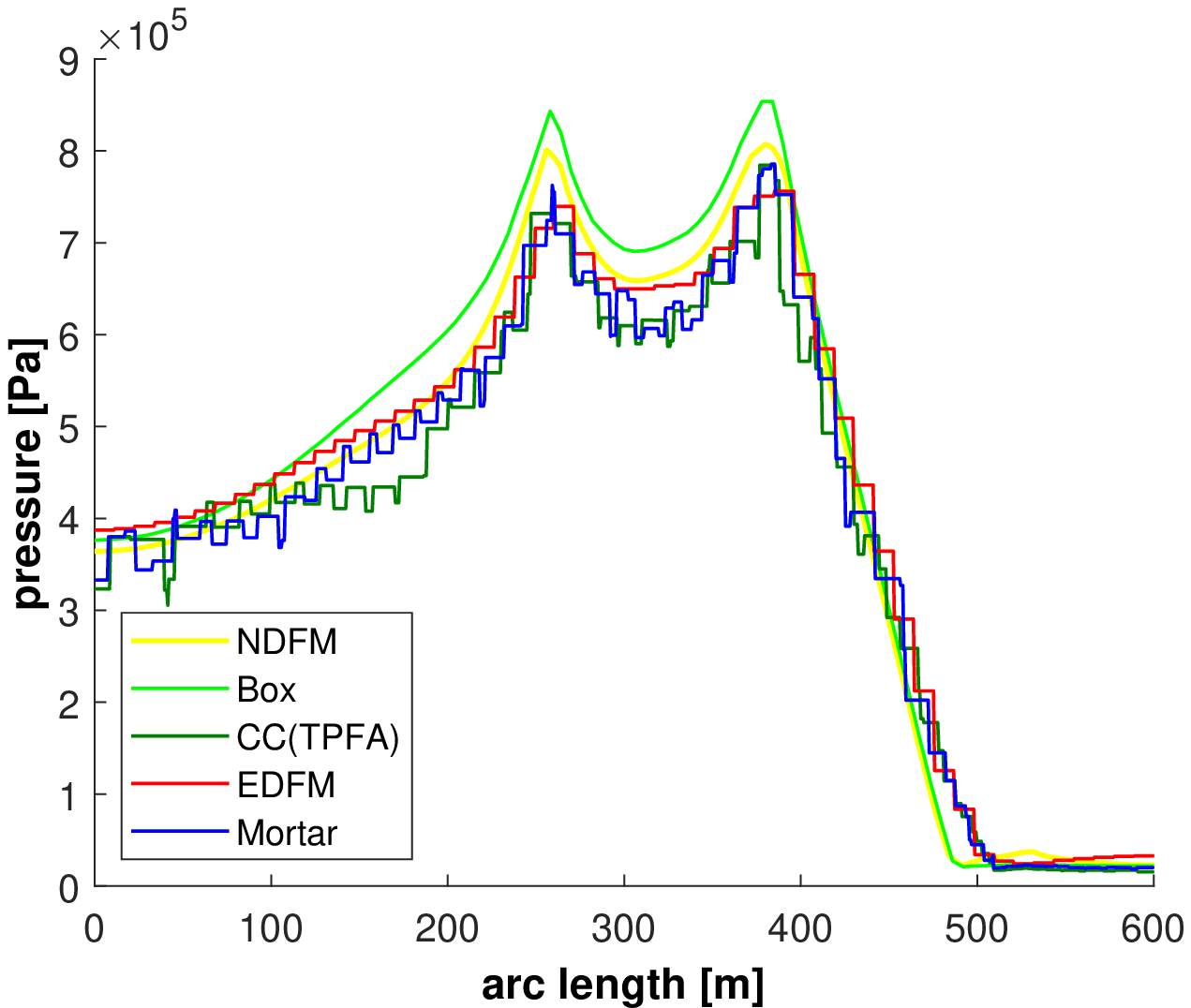}}
\caption{Comparison of slices of pressure along lines $y=500$\si{\m} and $x=625$\si{m} of Example \ref{ex6}\label{fig:ex6Slices}}
\end{figure}

We give a comparison of profiles of pressure of different methods as we did previously, but unfortunately the reference solution is not available in \cite{Website} because of the demand of tremendous gridcells in this test. For a similar reason, the XFEM-class method, although can be employed in principle, is too inconvenient to be practically implemented here. Therefore, there are only four participants in the comparison. We slice the solution of NDFM along the horizontal line $y = 500\si{\m}$ and vertical line $x = 625\si{m}$ on the coarse and fine meshes, respectively, and plot them with the profiles of Box-DFM, CC-DFM, EDFM, and Mortar-DFM in Figure \ref{fig:ex6Slices}. The degrees of freedom, sparsity and conditional number of linear systems in different methods are gathered in Table \ref{ex6}.

\begin{table}[!htb]
\centering
\caption{\label{tab:ex6} Evaluation data of different algorithms in Example \ref{ex6}.}
\begin{tabular}{|c|c|c|c|c|}
  \hline
method & d.o.f & mesh &sparsity&$||\cdot||_2$-cond\\\hline
Box-DFM&$5563$& conforming  & 1.2\textperthousand &9.3E05 \\\hline
CC-DFM&$8481$& conforming & 0.5\textperthousand &5.3E06  \\\hline
EDFM&$3599$& non-conforming & 1.4\textperthousand &4.7E06\\\hline
Mortar-DFM&$25258$& conforming & 0.2\textperthousand & 2.2E17 \\ \hline
NDFM&$9464$& non-conforming & 0.9\textperthousand & 3.3E06\\ \hline
NDFM&$26274$& non-conforming & 0.3\textperthousand & 1.5E07\\ \hline
\end{tabular}
\end{table}

\begin{figure}[!htbp]
\subfigure[Domain and fracture setting 1]{\includegraphics[width = 3in]{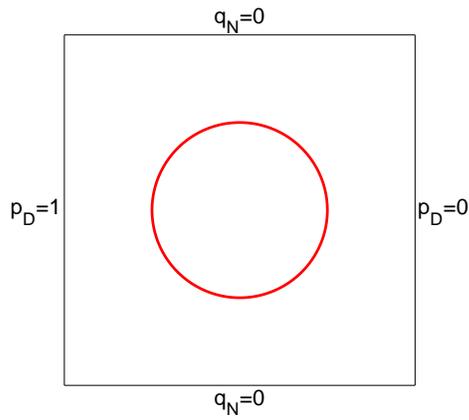}}
\subfigure[Domain and fracture setting 2]{\includegraphics[width = 3in]{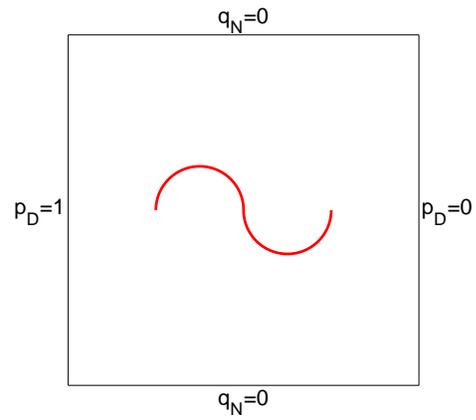}}\\
\subfigure[Reference solution 1]{\includegraphics[width = 3in]{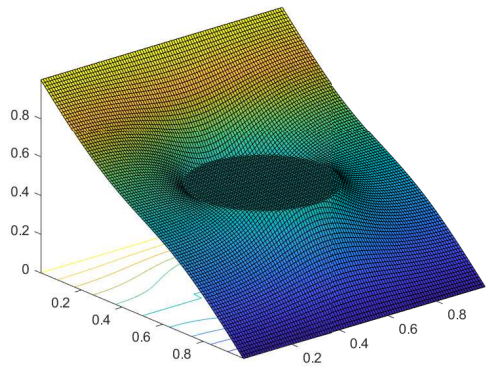}}
\subfigure[Reference solution 2]{\includegraphics[width = 3in]{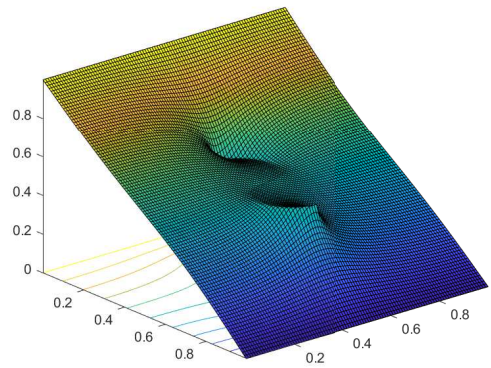}}\\
\caption{Fracture settings and reference solutions of Example \ref{ex7}\label{fig:ex7references}}
\end{figure}

\begin{figure}[!htbp]
\subfigure[Solution on $30\times30$ mesh, fracture 1]{\includegraphics[width = 3in]{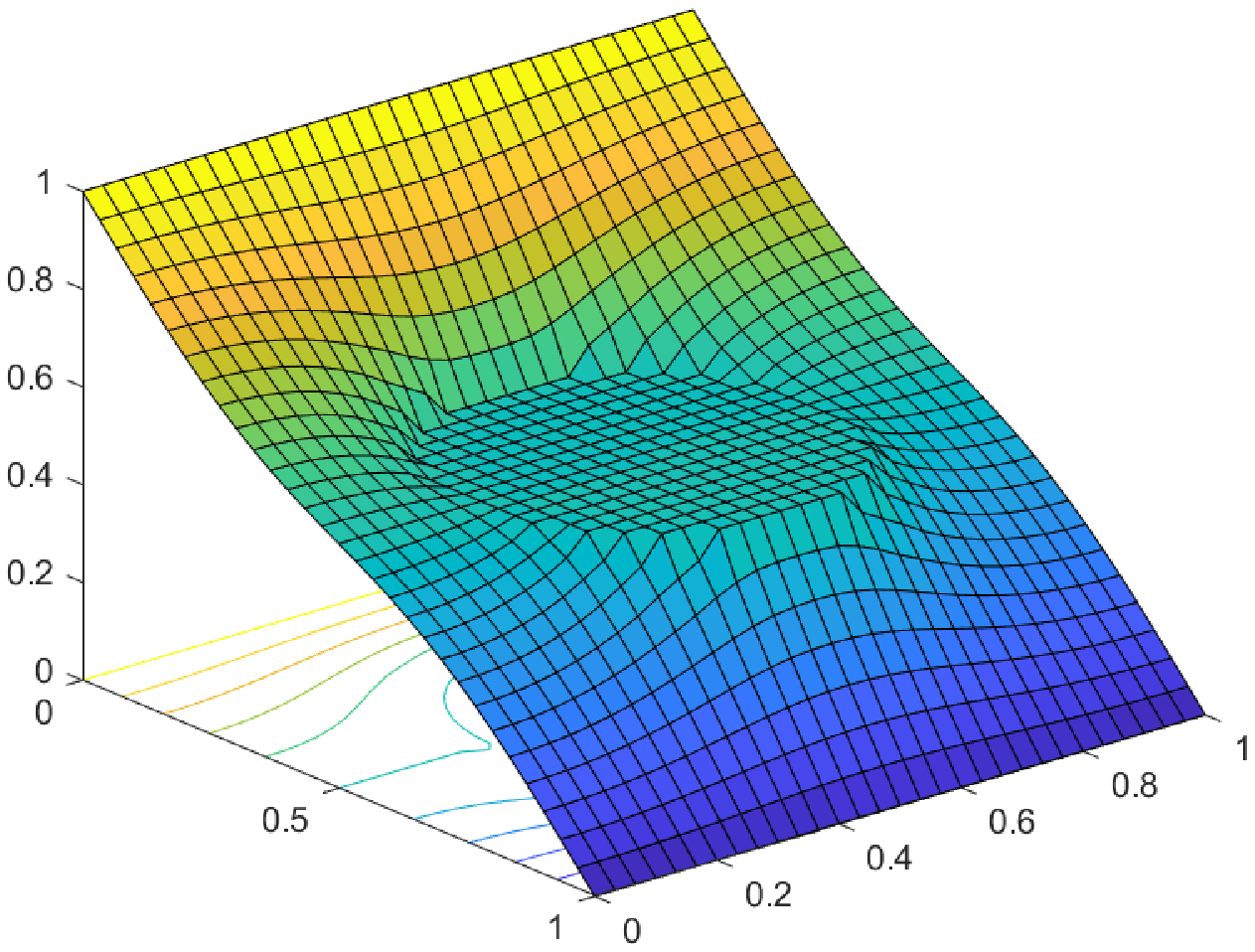}}
\subfigure[Solution on $30\times30$ mesh, fracture 2]{\includegraphics[width = 3in]{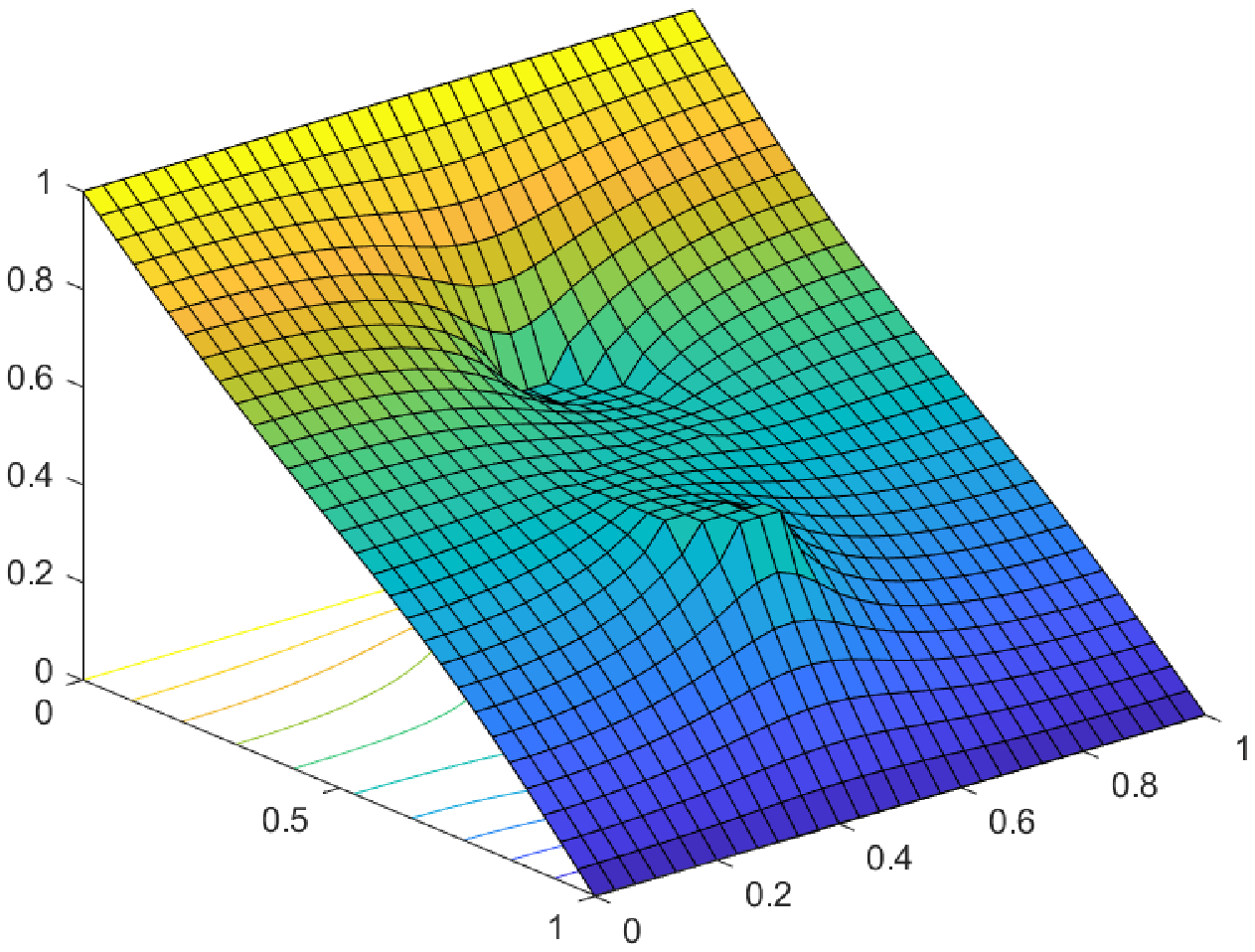}}\\
\subfigure[Solution on $60\times60$ mesh, fracture 1]{\includegraphics[width = 3in]{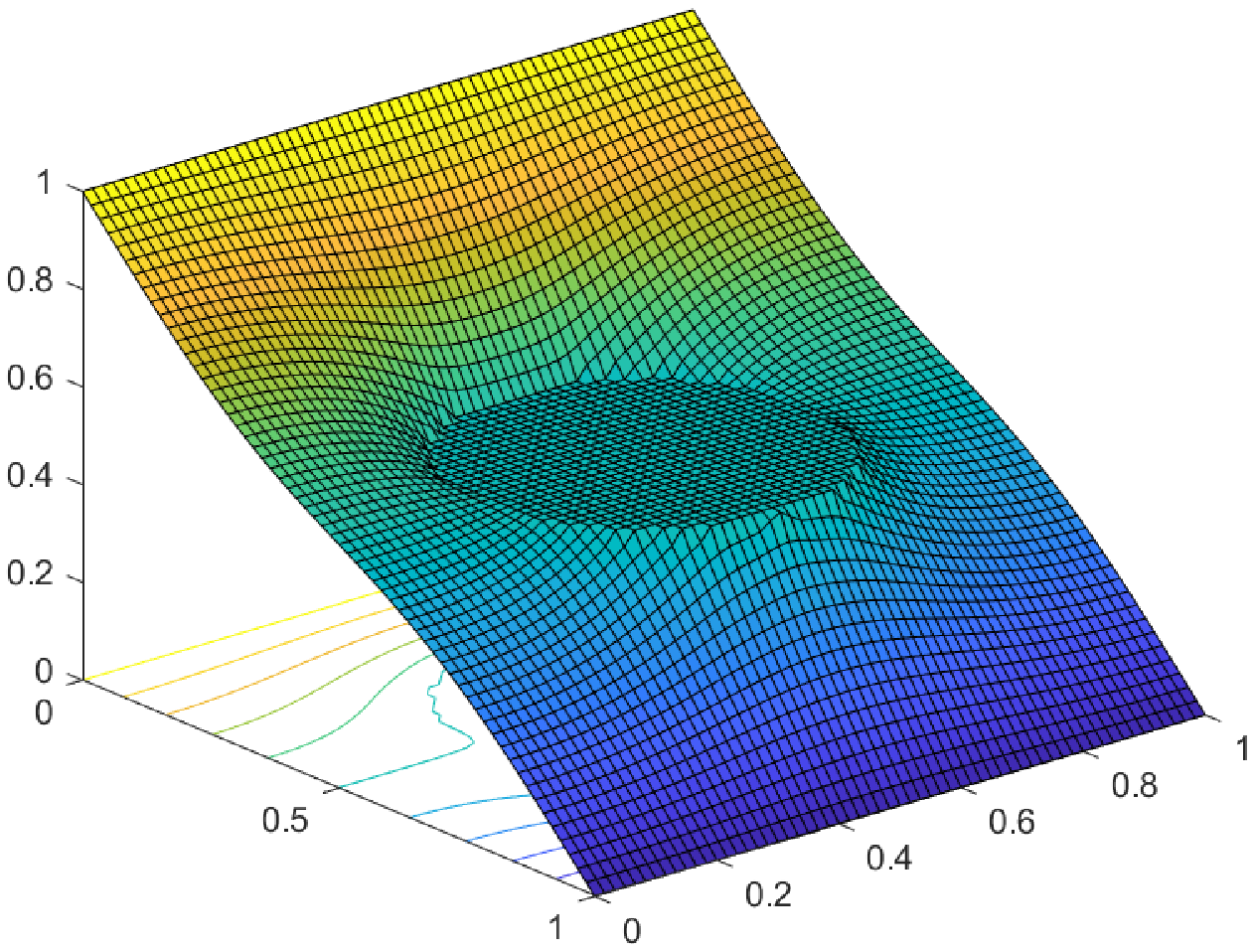}}
\subfigure[Solution on $60\times60$ mesh, fracture 2]{\includegraphics[width = 3in]{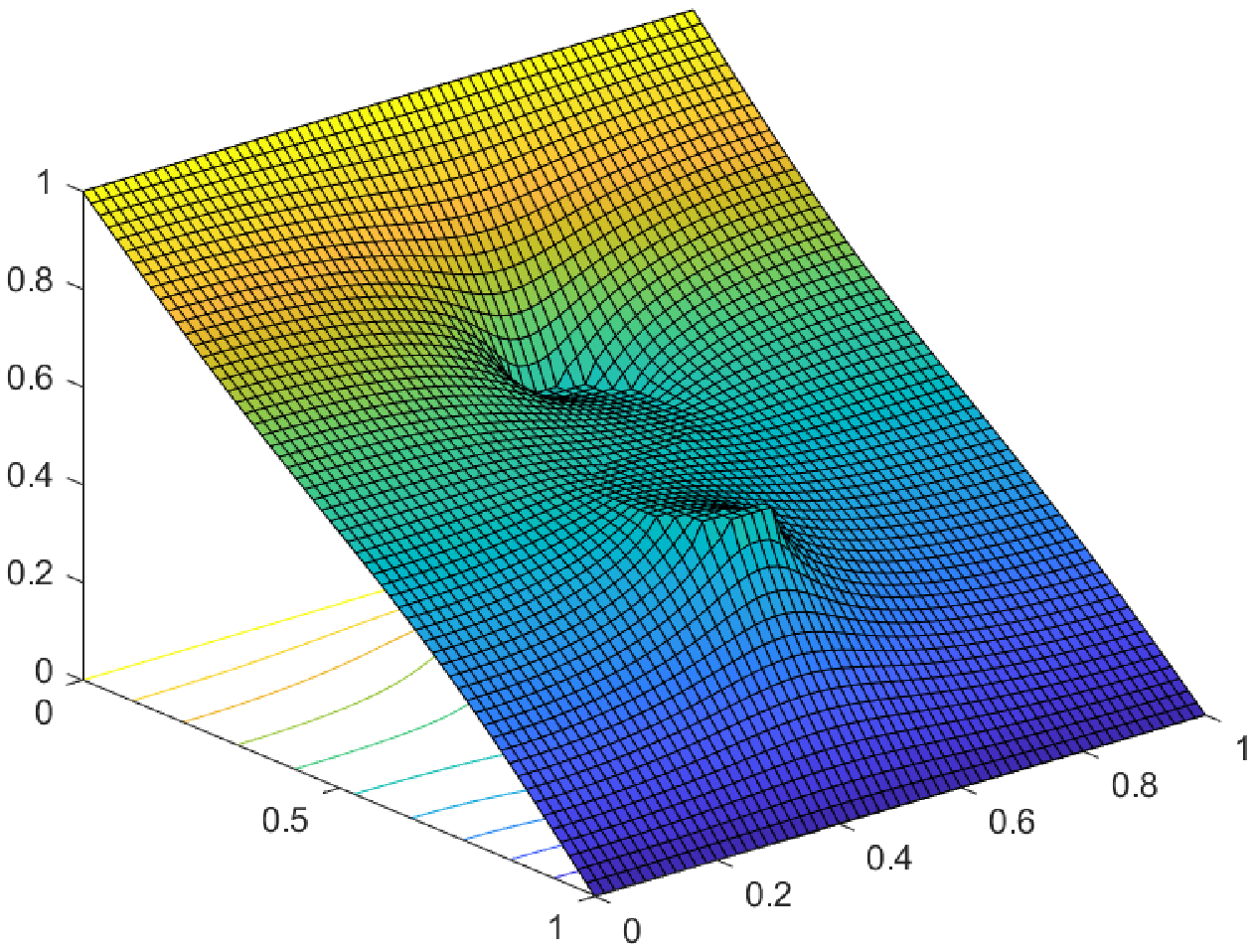}}\\
\caption{Solutions on different meshes of Example \ref{ex7}}\label{fig:ex7Solutions}
\end{figure}

\begin{figure}[!htbp]
\subfigure[Pressure along $y=0.5$, fracture 1]{\includegraphics[width = 3in]{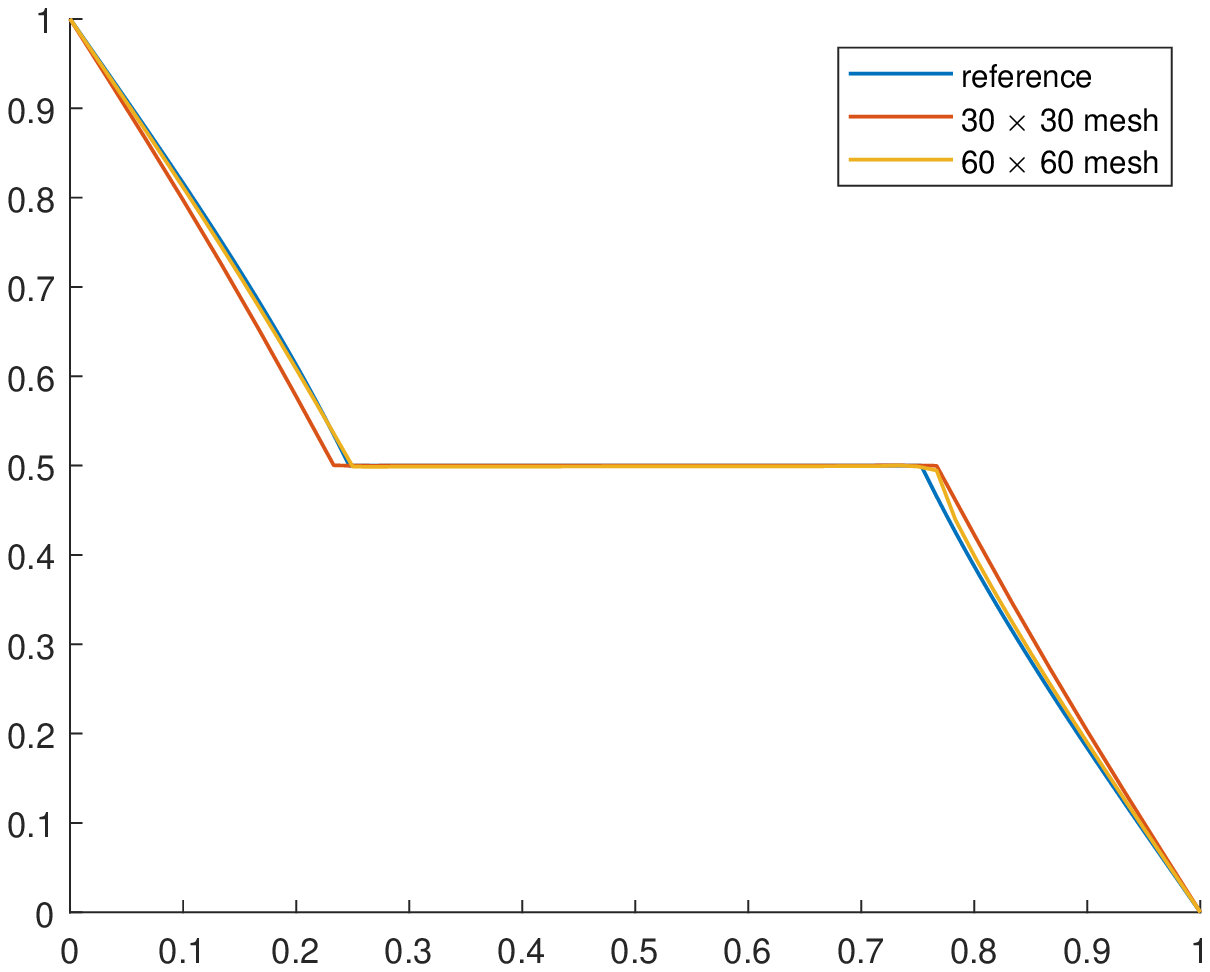}}
\subfigure[Pressure along $y=0.5$, fracture 2]{\includegraphics[width = 3in]{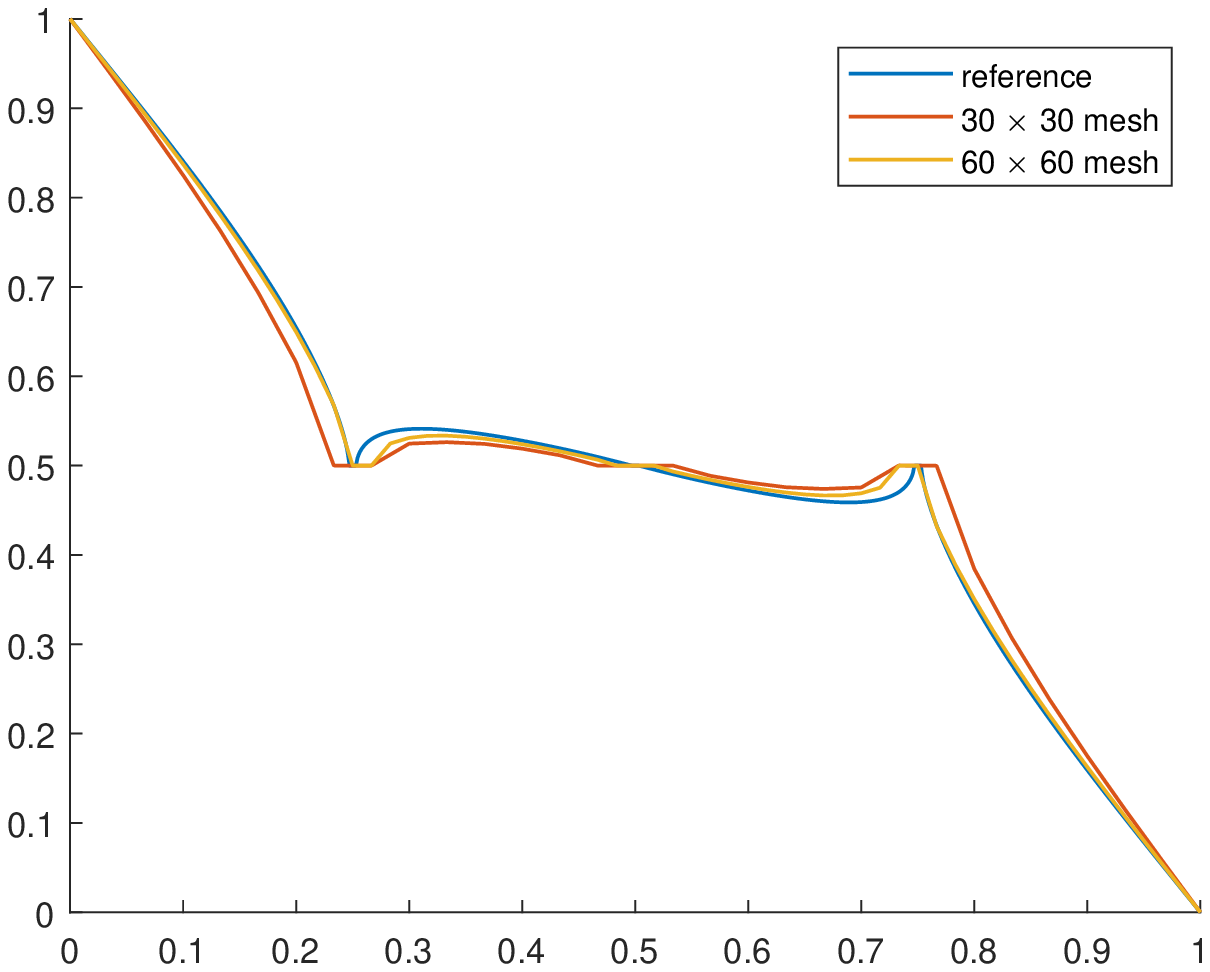}}\\
\subfigure[Pressure along $x=0.3$, fracture 1]{\includegraphics[width = 3in]{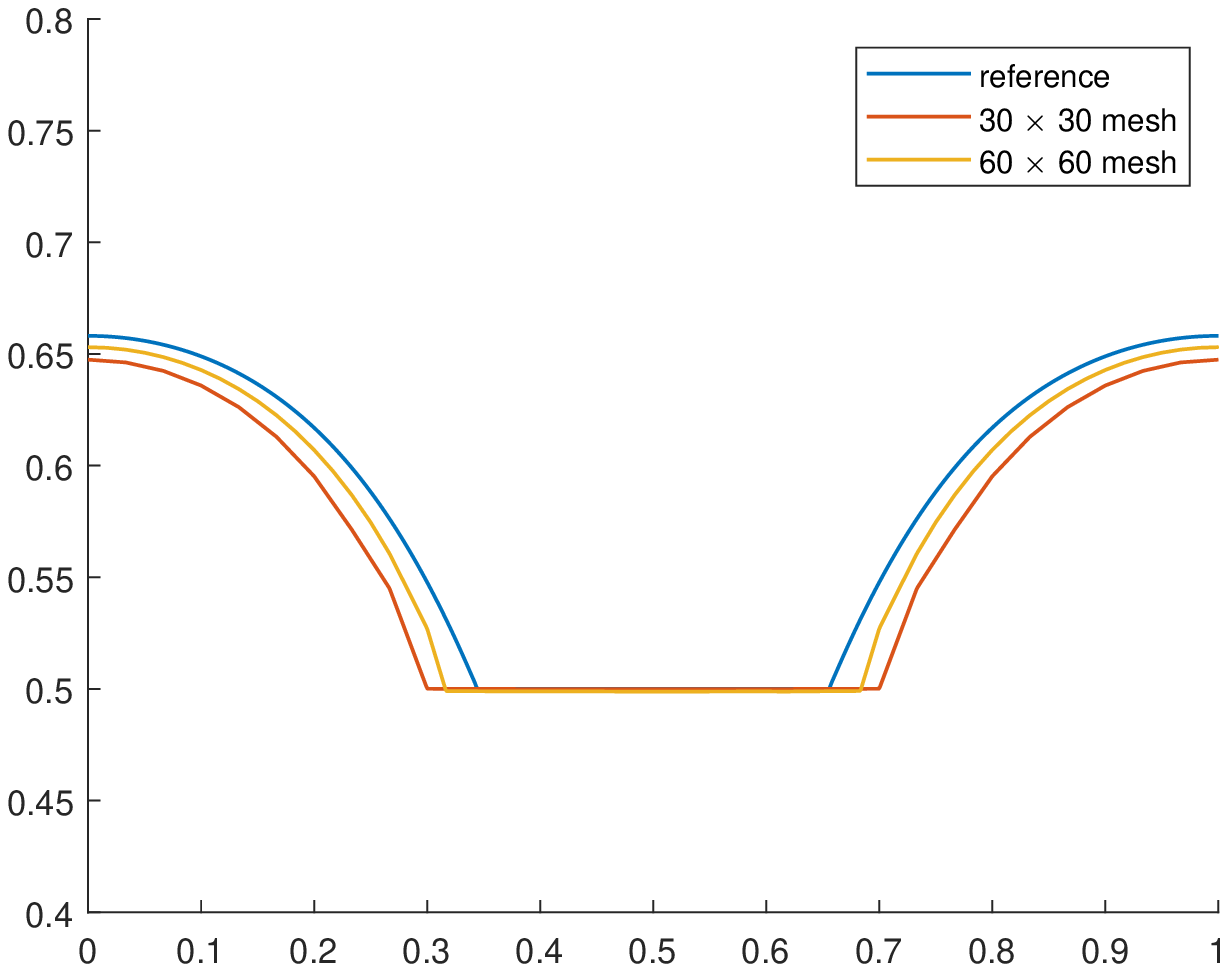}}
\subfigure[Pressure along $x=0.3$, fracture 2]{\includegraphics[width = 3in]{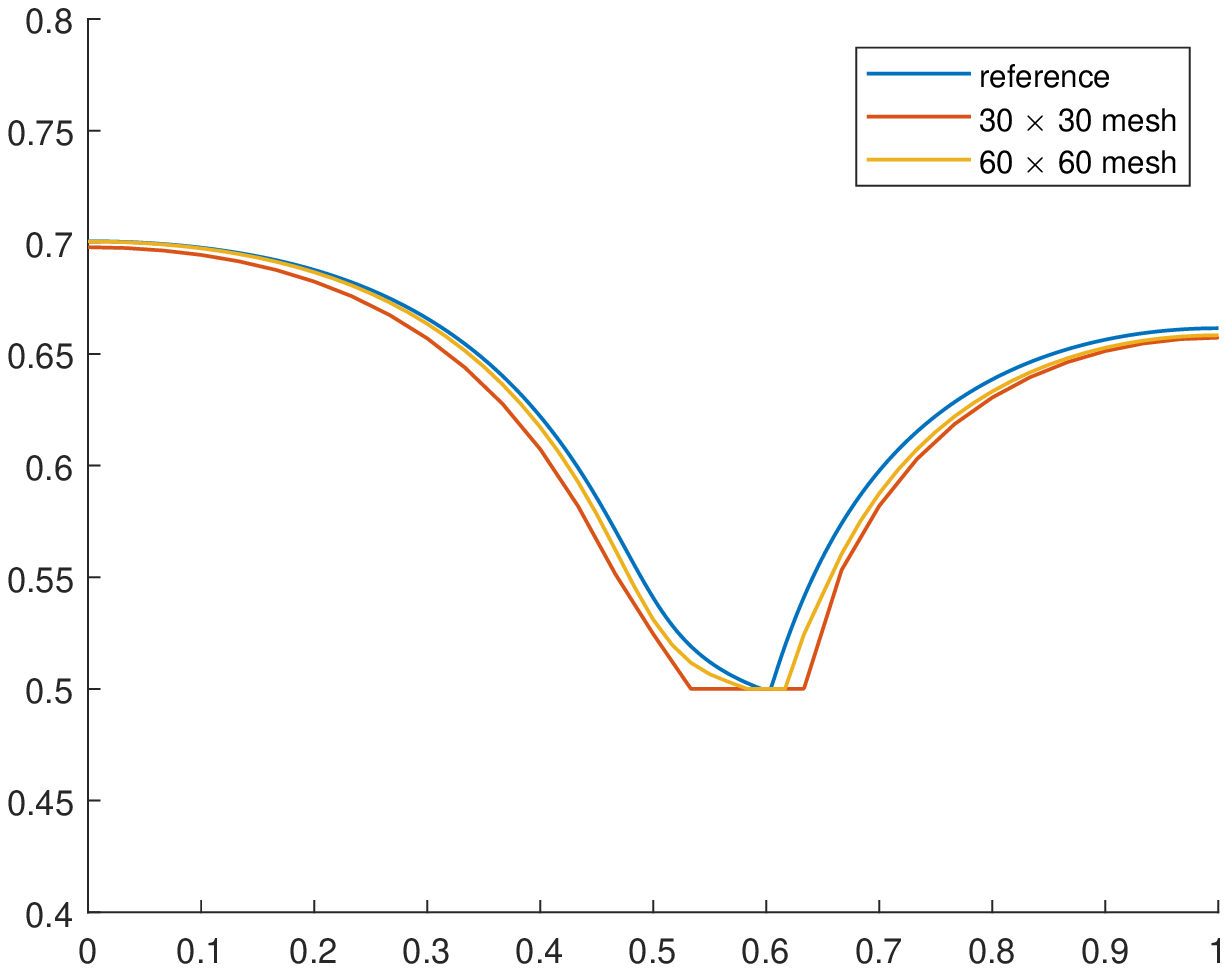}}\\
\subfigure[Pressure along $x=0.4$, fracture 1]{\includegraphics[width = 3in]{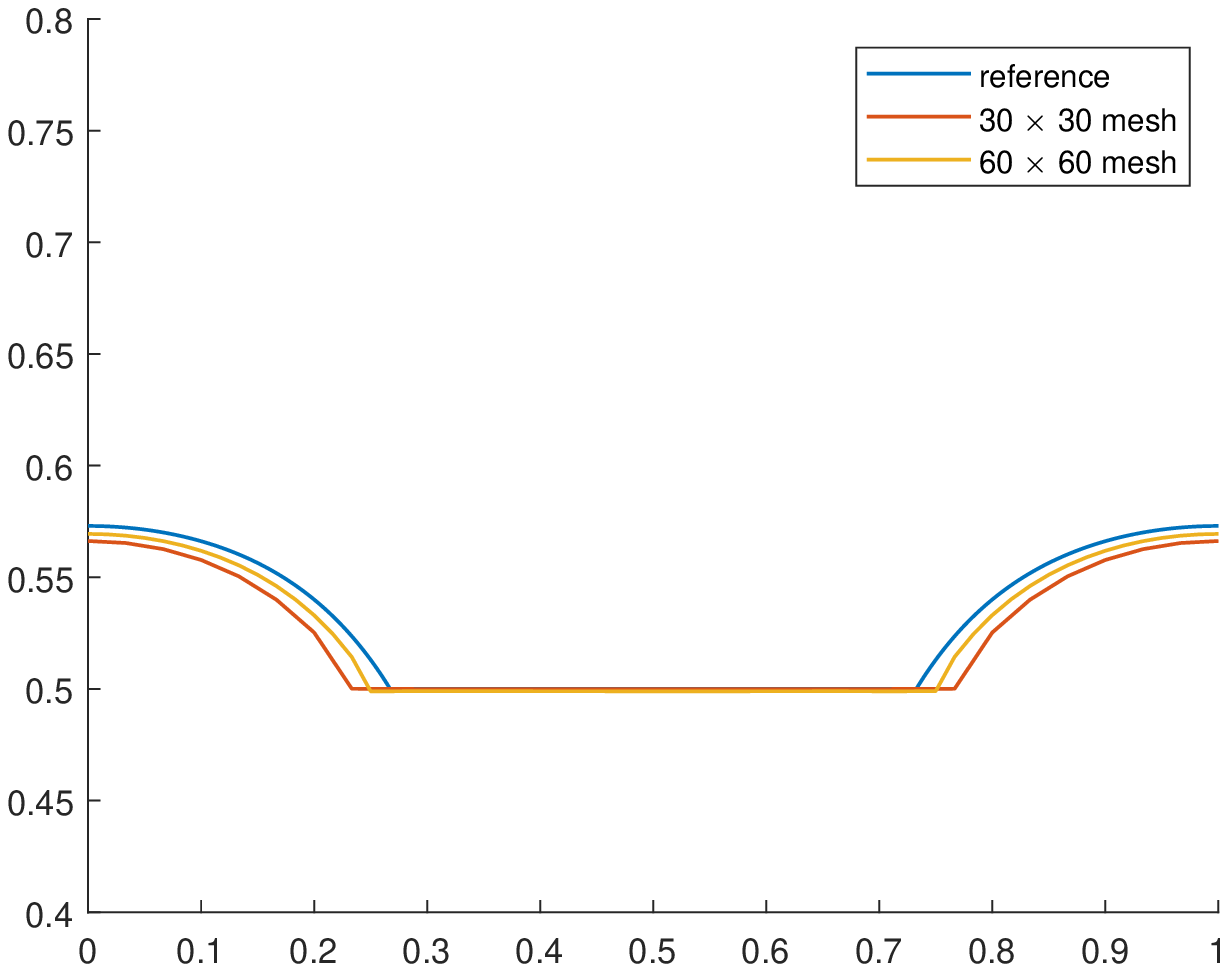}}
\subfigure[Pressure along $x=0.4$, fracture 2]{\includegraphics[width = 3in]{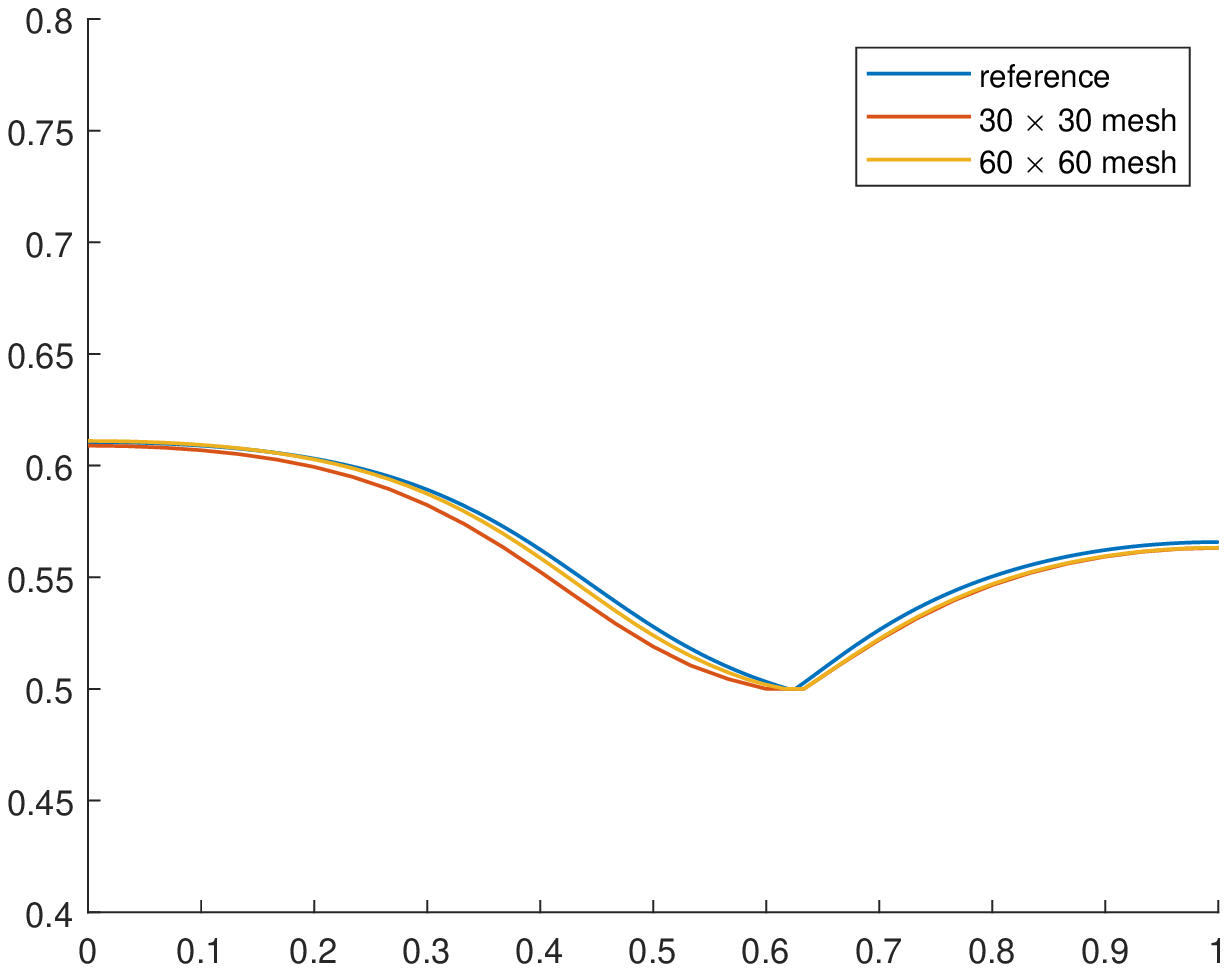}}
\caption{Pressure profiles along different lines of Example \ref{ex7}}\label{fig:ex7Slices}
\end{figure}

\begin{ex} \label{ex7}
\textbf{Curved Fractures}

In this example, we show the validity of the NDFM for curved fractures. Two types of curves are tested. The first one is a circle (closed) and the second one is a combination of semicircles (non-closed), See Figure \ref{fig:ex7references}. The computational domain is $[0,1]\times[0,1]$.The parametric equations of the curves are
$$x_1(t)=\frac{1}{2}+\frac{1}{4}\cos(t),\quad y_1(t)=\frac{1}{2}+\frac{1}{4}\sin(t),\quad 0\leq t\leq 2\pi,$$
and
$$ x_2(t) =
 \begin{cases}
      \frac{3}{8}-\frac{1}{8}\cos(t), & 0\leq t\leq \pi \\
      \frac{5}{8}+\frac{1}{8}\cos(t), & \pi\leq t\leq 2\pi  ,
   \end{cases},\quad y_2(t)=\frac{1}{2}+\frac{1}{8}\sin(t),\quad 0\leq t\leq 2\pi
,$$
respectively. The thickness of fractures is set to be $0.005$, and the permeability of porous matrix and fractures regions are $1$ and $10^{7}$, respectively. Moreover, the Dirichlet boundary conditions $p_D=1$ and $p_D=0$ are imposed on the left and right boundaries, respectively, and the top and bottom boundaries are impervious. See Figure \ref{fig:ex7references} for an illustration of the domain, fractures and boundary conditions setting.
\end{ex}
We employ the standard Galerkin finite element method for the equi-dimensional model on a fine mesh ($N_x=N_y=1000$) to give the reference solution. The surface and contour of the reference pressure are given in the Figure \ref{fig:ex7references}, and the results of NDFM on different rectangular meshes are shown in Figure \ref{fig:ex7Solutions}. Moreover, a comparison of pressure profiles sliced along $y=0.5, x=0.3$ and $x=0.4$ are shown in Figure \ref{fig:ex7Slices}, from which we can see that they match well.

\begin{figure}[!htbp]
\subfigure[Fracture setting 1]{\includegraphics[width = 3in]{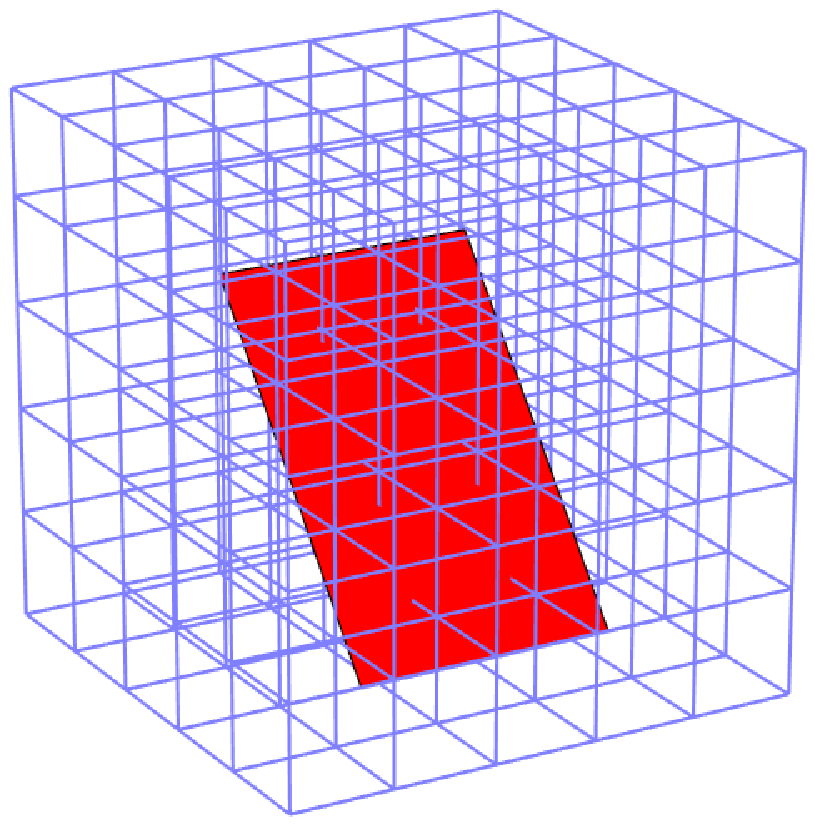}}
\subfigure[Fracture setting 2]{\includegraphics[width = 3in]{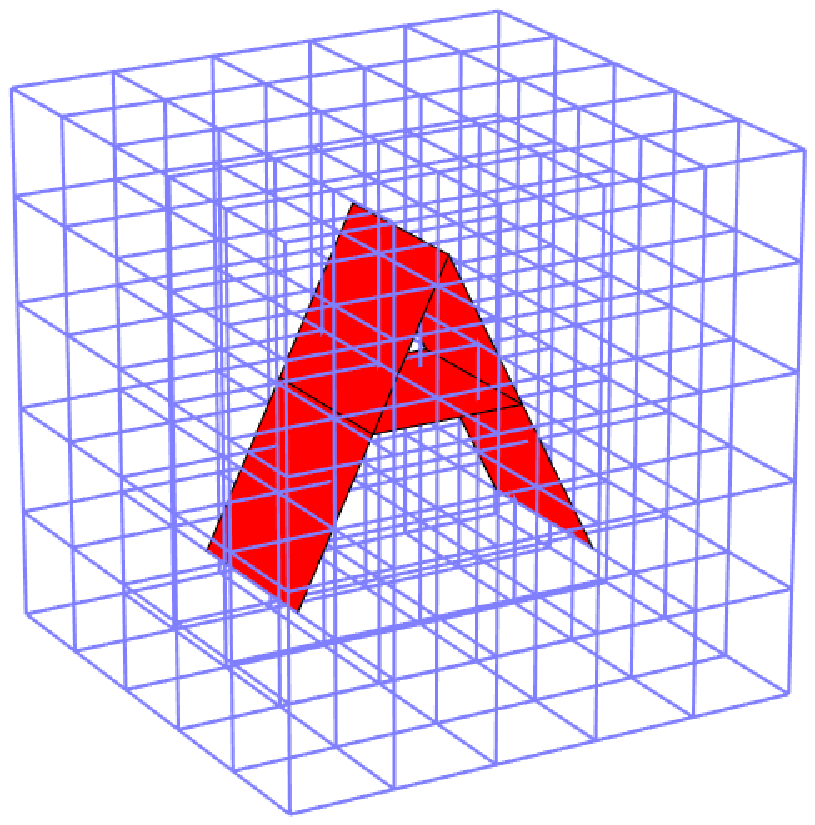}}\\
\caption{Fracture settings in $3$D space of Example \ref{ex8}, visualization algorithm from \cite{visualization}\label{fig:ex8fracturesetting}}
\end{figure}

\begin{ex} \label{ex8}
\textbf{3D cases}

In this last example, we give a brief demonstration on how to extend the NDFM to 3D cases and show the performance of the 3D algorithm qualitatively by presenting several 2D slices from the 3D space.

The extension of the model to $3D$ is almost trivial. By analogy, the hybrid-dimensional representation of the permeability tensor is
\begin{equation}
\displaystyle{\bf{K}} = \textbf{K}_m + \sum^{L}_{i=1} \epsilon_i k_{fi}\delta_i(\cdot)\mathbbm{1}_i(\cdot)\left( {\bf{I}} - \bm{\sigma}_i\bm{\sigma}_i^{T} \right),
\end{equation}
where ${\bf{I}}$ is the identity tensor and $\bm{\sigma}_i$ is the unit normal vector of the $i$th fracture. The expressions abbreviated by $(\cdot)$ in the functions $\delta_i$ and $\mathbbm{1}_i$ are even more complicated but the geometric meaning is still concise: the geometric information that determines a fracture. The corresponding variational formulation is
\begin{equation}
\displaystyle\int_{\Omega}\left(\textbf{K}_m\nabla p\right)\cdot\nabla v ~dx dy dz+ \sum^{L}_{i=1}\int_{{\bm S}_i}\epsilon_i k_{fi} \left( \nabla p\cdot \nabla v-\frac{\partial p}{\partial\sigma_i}\frac{\partial v}{\partial\sigma_i} \right)~dS = \int_{\Omega}f v~dx dy dz+\int_{\Gamma_N}q_N v~dS,
\end{equation}
where ${\bm S}_i$ denotes the $i$th fracture surface.

Unfortunately, due to the restriction of the memory of the computer and the computational cost, we are not able to give a fully resolved reference solution for the 3D problems, thus no quantitative measurement for the error is available here.

Two types of fractures are simulated. The first one is a single planar fracture and the second one is an "A-"shaped planar fractures network, see Figure \ref{fig:ex8fracturesetting}. The single fracture is a rectangle with four vertices $(\frac{1}{4}, \frac{1}{5}, \frac{1}{10}),$ $ (\frac{1}{4}, \frac{7}{10}, \frac{7}{10}), $ $(\frac{3}{4}, \frac{7}{10}, \frac{7}{10}),$ $ (\frac{3}{4}, \frac{1}{5}, \frac{1}{10})$. The fractures network is composed of 3 rectangles, whose vertices are $(\frac{1}{2}, \frac{1}{3}, \frac{4}{5}),$ $ (\frac{1}{5}, \frac{1}{3}, \frac{1}{5}),$ $ (\frac{1}{5}, \frac{2}{3}, \frac{1}{5}),$ $ (\frac{1}{2}, \frac{2}{3}, \frac{4}{5})$ and $(\frac{1}{2}, \frac{1}{3}, \frac{4}{5}),$ $ (\frac{4}{5}, \frac{1}{3}, \frac{1}{5}),$ $ (\frac{4}{5}, \frac{2}{3}, \frac{1}{5}), $ $(\frac{1}{2}, \frac{2}{3}, \frac{4}{5})$, and $(\frac{7}{20}, \frac{1}{3}, \frac{1}{2}), $ $(\frac{13}{20}, \frac{1}{3}, \frac{1}{2}),$ $ (\frac{13}{20}, \frac{2}{3}, \frac{1}{2}), (\frac{7}{20}, \frac{2}{3}, \frac{1}{2})$, respectively. The thickness of fractures is $10^{-3}$ and the permeability is $10^{8}$. The domain $\Omega=[0,1]\times[0,1]\times[0,1]$  and the permeability of porous matrix is $1$. We impose the Dirichlet boundary condition $p_D=1$ and $p_D=0$ on the left and right boundary of $\Omega$. The other boundaries are impervious.
\end{ex}

\begin{figure}[!htbp]
\subfigure[slice plans in $x$ direction, fracture 1]{\includegraphics[width = 3in]{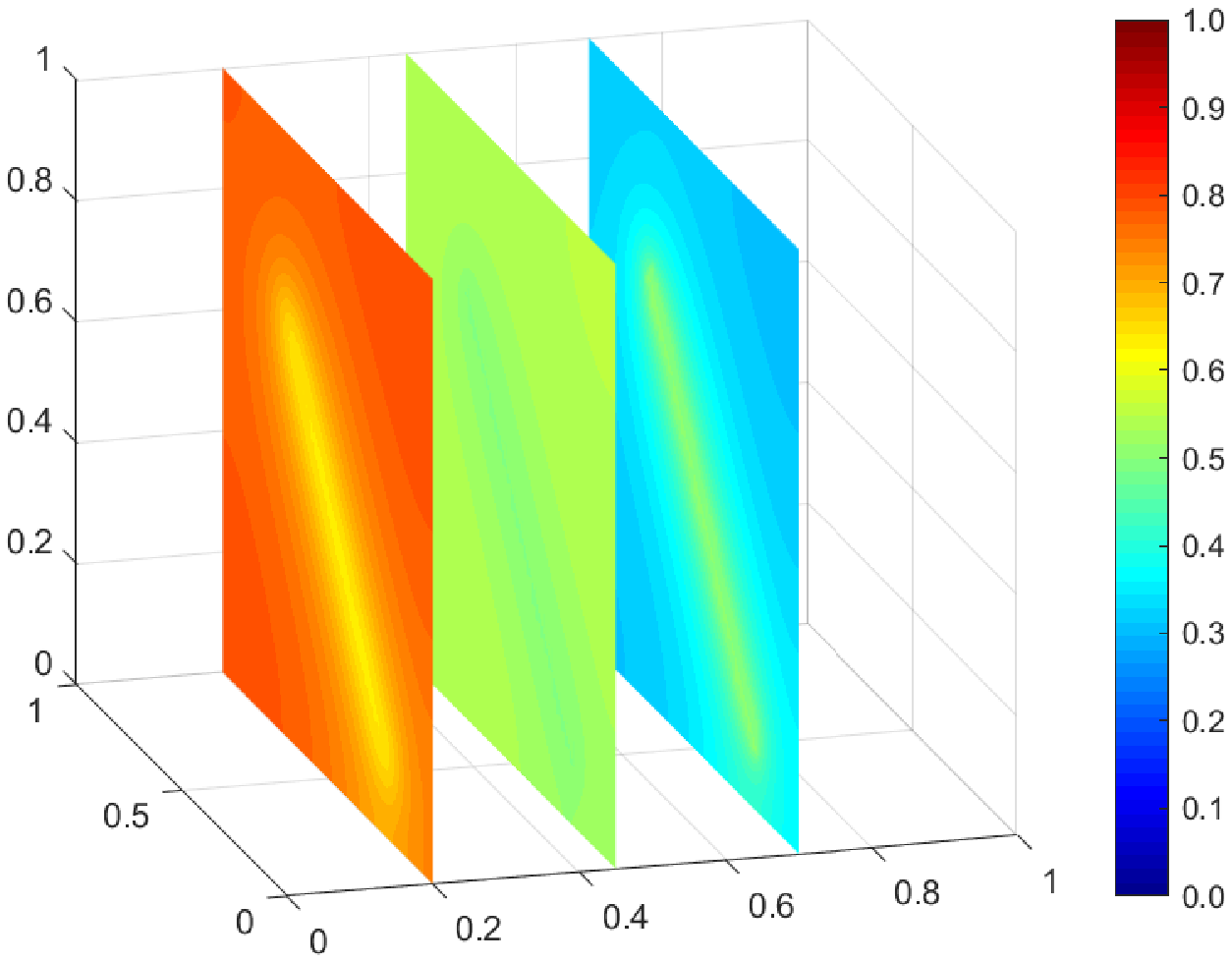}}
\subfigure[slice plans in $x$ direction, fracture 2]{\includegraphics[width = 3in]{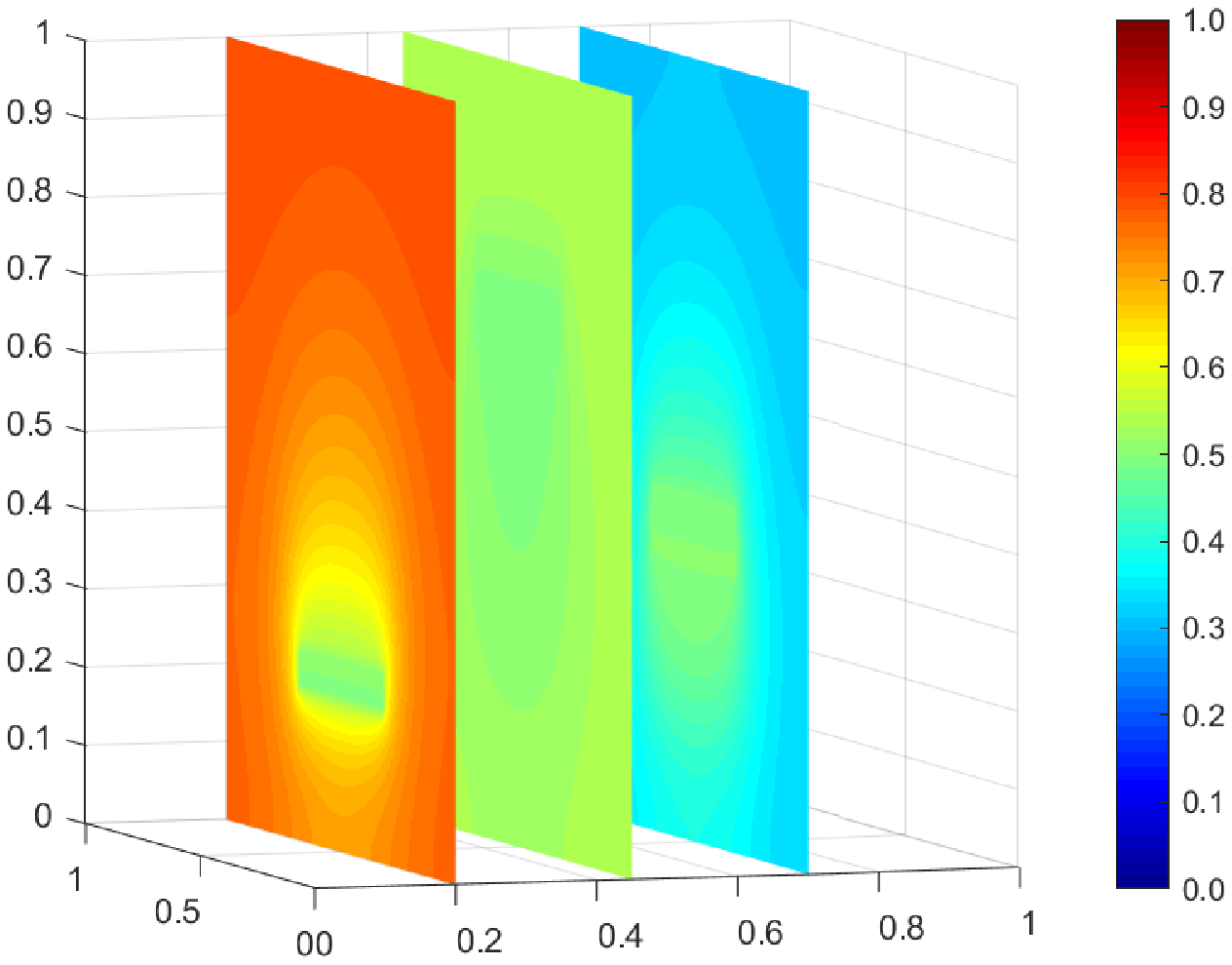}}\\
\subfigure[slice plans in $y$ direction, fracture 1]{\includegraphics[width = 3in]{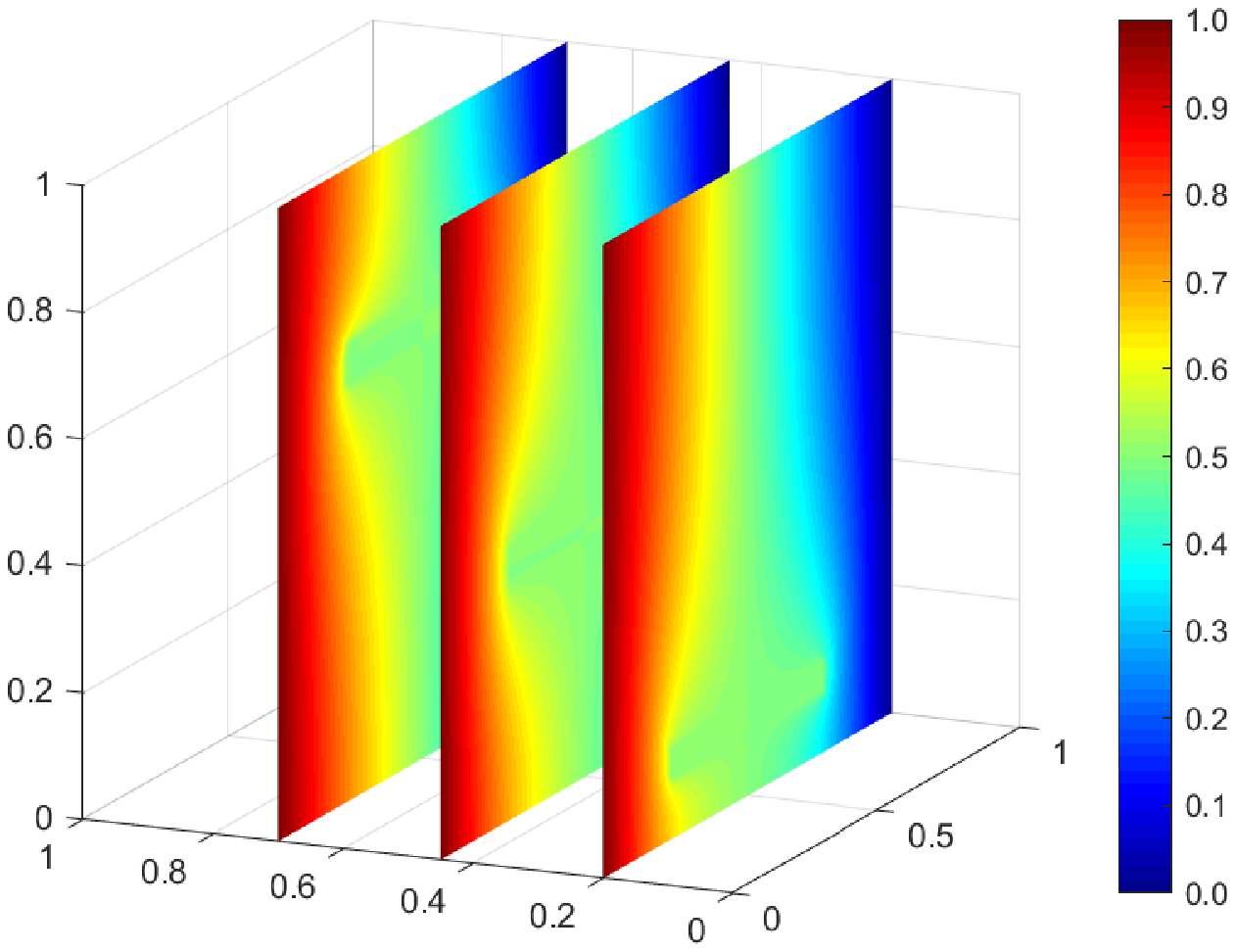}}
\subfigure[slice plans in $y$ direction, fracture 2]{\includegraphics[width = 3in]{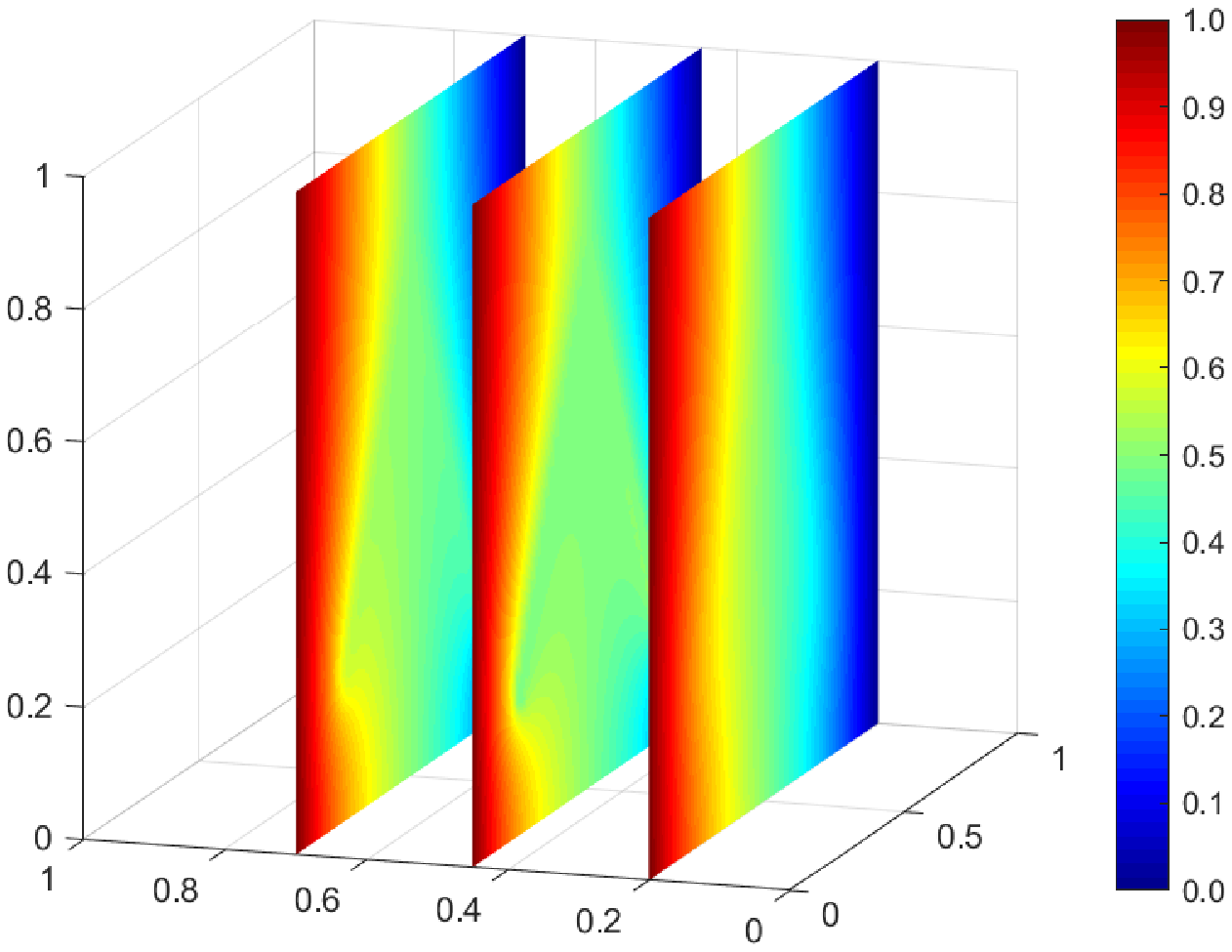}}\\
\subfigure[slice plans in $z$ direction, fracture 1]{\includegraphics[width = 3in]{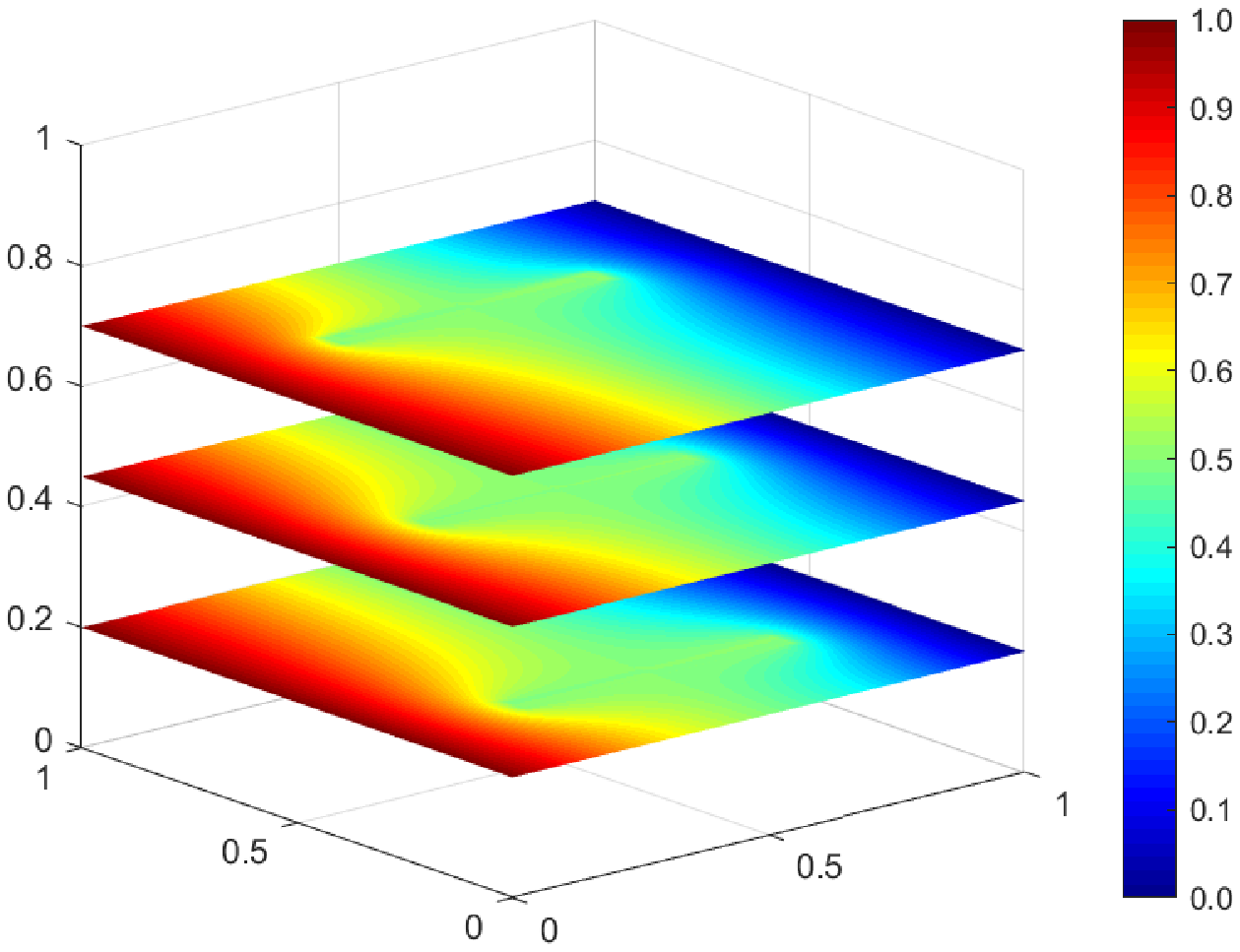}}
\subfigure[slice plans in $z$ direction, fracture 2]{\includegraphics[width = 3in]{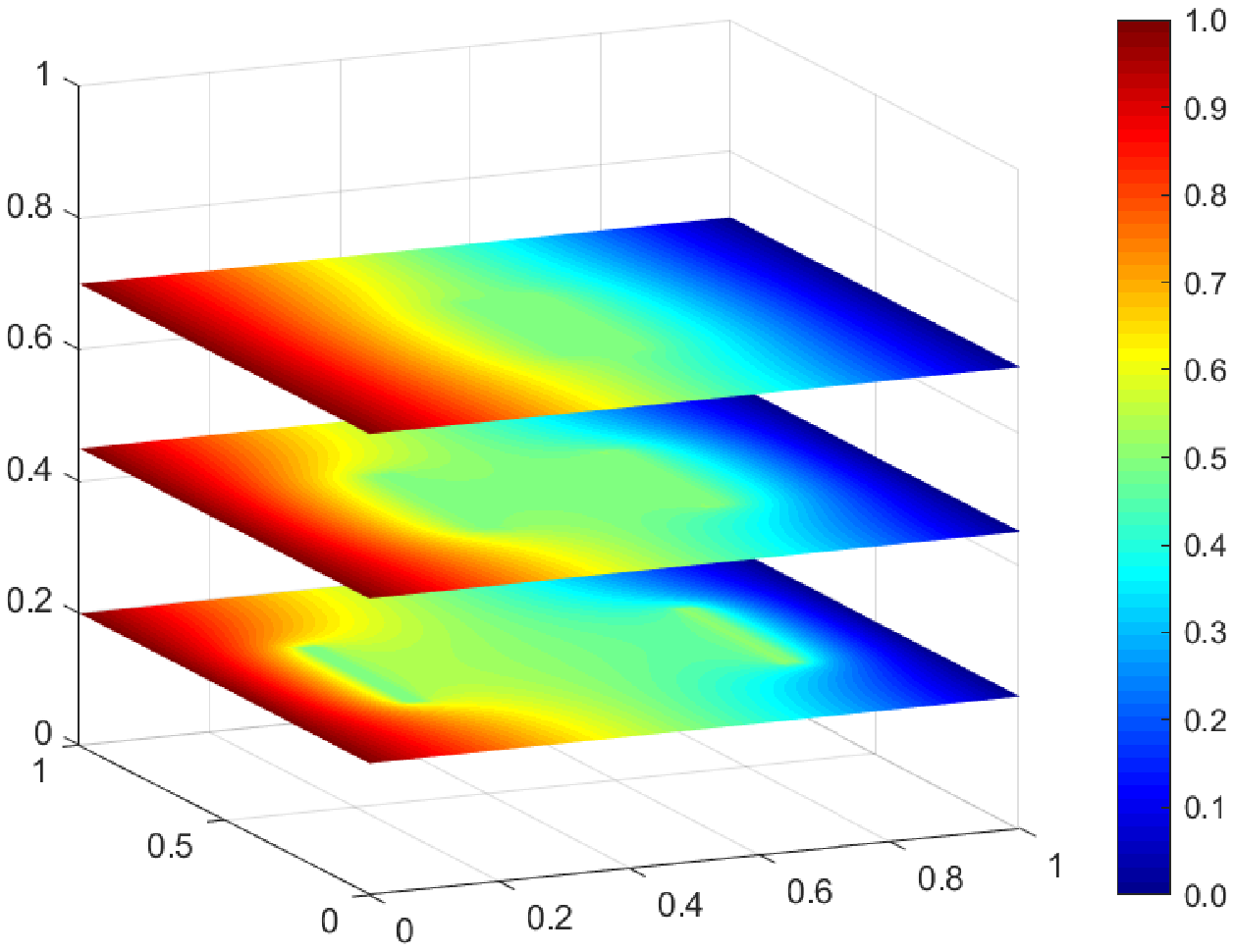}}\\
\caption{Slice planes of the pressure in $3$D space of Example \ref{ex8}  \label{fig:ex8slices}}
\end{figure}
We present the slice planes of the pressure along $x$ direction $x=0.2, 0.45, 0.7$, $y$ direction $y=0.2, 0.45, 0.7$, and $z$ direction $z=0.2, 0.45, 0.7$ in Figure \ref{fig:ex8slices}.
From the figures we can see that the highly conductive zone indicated by the slices coincide with the position of the fracture.

\section{Concluding remarks}
In this paper, we explored the hybrid-dimensional representation of permeability tensor of fractured media and constructed the Galerkin finite element discrete fracture model on non-conforming meshes based on it. Analytical analysis and numerical tests showed its consistency with the traditional discrete fracture model on conforming meshes and its accuracy and efficiency on non-conforming meshes.

There are some major limitations of the approach for further improvement. First, the model is established based on the assumption that the fracture has very tiny thickness and high permeability. Only under this condition can we reduce the fracture to a $1$D object and use Dirac-$\delta$ function to represent its permeability. If the thickness of the fracture is not small enough, the model will not be very suitable since the pressure jump across the fracture is not representable in the model but it's non-negligible in this case. Second, the model only works for conductive fracture. For the other type of fracture, the barrier element, the model is not applicable. Third, due to the property of finite element method, the scheme in the paper is not locally mass conservative. This shortcoming is not problematic in the steady-state single-phase flow problem, but might be more obvious when it is applied to two-phase or multi-phase fluid flow simulations.

Several possible improvements can be made regarding these limitations. For the barrier elements, the suitable adaptation for the model is our ongoing work. As for the mass conservation, some locally mass conservative methods, like the discontinuous Galerkin method \cite{DGconservation} or modified conservative Galerkin methods \cite{FEMconservation} could be employed to the hybrid-dimensional model problem \eqref{PoissonEq}, \eqref{HybridDimensionK} and \eqref{BVP} instead of the traditional Galerkin methods.

\begin{appendices}
\appendix

\section{Geometric data of fracture network in Example \ref{ex6}}

\begin{Verbatim}
NUMBER,START_X(m),START_Y(m),END_X(m),END_Y(m)
1,269.611206,152.05243,356.9240112,310.14123
2,249.5117187,514.990780001,272.218872,470.97082
3,258.3590698,515.574580001,271.9851684,490.9682
4,270.6622924,524.702640001,269.1347046,147.78143
5,355.8302002,348.479800001,337.5810733205,600
6,366.9730835,338.132990001,426.9185141723,600
7,198.237915,222.724420001,175.1561889,597.603030001
8,151.2785034,261.724610001,154.4623059774,600
9,29.5026855,300.724610001,96.3599853,514.82739
10,386.0808105,33.3621800002,440.585083,275.191830001
11,459.6350708,40.2413900001,461.751709,204.812620001
12,297.180603,237.62103,468.1018066,40.2413900001
13,312.5264892,272.01678,417.3016967,140.7832
14,330.5181884,298.47522,439.5266723,156.6582
15,340.5723877,320.70019,367.5598755,286.304380001
16,492.9725952,312.762820001,576.5811157,419.6546
17,505.6726684,309.05859,576.0520019,405.367190001
18,537.4227905,297.94598,623.3187866,376.68463
19,322.5338745,380.76941,521.8778076,593.552180001
20,344.9320678,481.56122,409.8867798,503.959410001
21,371.8098755,468.12219,510.6787109,383.009210001
22,432.2849731,510.678830001,642.8280029,374.04999
23,527.528634971,600,700,473.015615092
24,0,333.73321,441.2443847,0
25,13.4389038,342.692380001,347.171875,595.791990001
26,22.3981933,450.203790001,311.3347778,291.176630001
27,26.8778076,506.199220001,199.343811,400.92779
28,44.7963867,528.597410001,365.0905151,342.692380001
29,378.5294189,309.095210001,512.918518,116.470640001
30,461.4027099,253.099610001,530.8370971,134.38922
31,347.171875,374.04999,640.5881958,253.099610001
32,490.5203857,268.77844,564.4343872,145.58844
33,47.0361938,181.425410001,53.7556152,306.85541
34,382.4152832,424.151000001,447.8997192,371.76343
35,587.9967651,394.78222,549.1029663,362.635190001
36,589.9812011,393.59161,527.6716919,313.8194
37,597.125,378.90722,533.6248169,295.960200001
38,533.6248169,448.75738,453.8527832,326.91638
39,511.7966919,461.85419,489.5715942,395.17901
40,565.3748779,425.34161,483.6184692,315.40698
41,534.4185791,407.482240001,467.3466186,315.803830001
42,627.2874756,527.3388,574.8999023,498.763610001
43,644.3532104,519.00439,586.4093017,490.03241
44,655.8626098,502.335630001,602.6812133,476.53863
45,415.355896,585.679380001,391.9401855,561.47003
46,417.3402099,578.535580001,397.8933105,554.326230001
47,403.0526733,592.029420001,382.0183105,561.86682
48,495.1278686,505.113580001,468.1403198,481.30121
49,533.6248169,254.84381,420.9121093,159.196590001
50,508.6217041,221.10943,441.152771,159.59363
51,418.5308838,229.04681,312.961914,93.3154300004
52,362.5714111,174.6748,322.883789,120.69983
53,357.8088989,216.3468,295.102478,114.74658
54,402.2589111,283.41882,366.1433105,226.66559
55,337.5681762,253.256220001,374.4776001,211.18744
56,386.7808838,264.765620001,509.8123169,101.25281
57,473.2996826,278.65643,561.0092163,144.909240001
58,471.7122192,253.653200001,554.6593017,129.034240001
59,559.0249023,219.125,567.3593139,153.64044
60,567.7561035,214.759400001,573.7092895,162.37182
61,574.8999023,215.553040001,579.6624145,173.88104
62,557.0404663,285.006410001,600.6968994,325.48761
63,565.3748779,283.022030001,607.0468139,323.503230001
\end{Verbatim}

\end{appendices}

\end{document}